\documentclass[12pt]{article}
\usepackage{e-jc}

\usepackage{amsmath,amssymb}
\usepackage{subfigure}
\usepackage{tikz}
\tikzstyle{every node}=[circle, draw, fill=white, inner sep=0pt, minimum width=4pt]
\tikzstyle{nodelabel}=[rounded corners, fill=none, inner sep=5pt, draw=none]
\usepackage{enumerate}
\theoremstyle{plain}
\newtheorem{notation}[theorem]{Notation}
\theoremstyle{remark}
\newtheorem{assertion}[theorem]{Assertion}

\dateline{Mar 18, 2019}{Dec 20, 2019}{TBD}

\MSC{05C70, 05C75, 05C40}

\title{On essentially $4$-edge-connected cubic bricks}

\author{Nishad Kothari\thanks{Supported by Austrian Science Foundation FWF (START Y463)
            and FAPESP Brazil \mbox{(2018/04679-1)}.}\\
\small University of Campinas, Brazil\\[-0.8ex]
\small\tt nishadkothari@gmail.com\\
\and
Marcelo H. de Carvalho\thanks{Supported by Fundect-MS and CNPq.}\\
\small UFMS Campo Grande, Brazil\\[-0.8ex]
\small\tt mhc@facom.ufms.br\\
\and
Cl{\'a}udio L. Lucchesi\thanks{Supported by CNPq.}\\
\small University of Campinas, Brazil\\[-0.8ex]
\small\tt lucchesi@ic.unicamp.br\\
\and
Charles H. C. Little\\
\small Massey University, New Zealand\\[-0.8ex]
\small\tt c.little@massey.ac.nz
}

\begin{document}

\newcommand{\mc}{matching covered}
\newcommand{\mcg}{matching covered graph}
\newcommand{\efec}{essentially \mbox{$4$-edge-connected}}
\newcommand{\efeccg}{essentially \mbox{$4$-edge-connected} cubic
    graph}
\newcommand{\efeccb}{essentially \mbox{$4$-edge-connected} cubic
    brick}
\newcommand{\efeccnnbb}{essentially \mbox{$4$-edge-connected} cubic
    non-near-bipartite brick}
\newcommand{\binv}{\mbox{$b$-invariant}}
\newcommand{\qbinv}{\mbox{quasi-$b$-invariant}}
\newcommand{\iso}{\simeq}

\maketitle

% E-JC papers must include an abstract. The abstract should consist of a
% succinct statement of background followed by a listing of the
% principal new results that are to be found in the paper. The abstract
% should be informative, clear, and as complete as possible. Phrases
% like "we investigate..." or "we study..." should be kept to a minimum
% in favor of "we prove that..."  or "we show that...".  Do not
% include equation numbers, unexpanded citations (such as "[23]"), or
% any other references to things in the paper that are not defined in
% the abstract. The abstract may be distributed without the rest of the
% paper so it must be entirely self-contained.  Try to include all words
% and phrases that someone might search for when looking for your paper.

%\begin{center}
\centerline{Dedicated to Professor U. S. R. Murty on the year of his 80th birthday}
%\end{center}

\begin{abstract}
    
Lov{\'a}sz (1987) proved that every matching covered graph $G$ may be uniquely decomposed into
a list of bricks (nonbipartite) and braces (bipartite); we let $b(G)$ denote the number of bricks.
An edge $e$ is removable if $G-e$ is also matching covered; furthermore, $e$ is $b$-invariant
if $b(G-e)=1$, and $e$ is quasi-$b$-invariant if $b(G-e)=2$.
(Each edge of the Petersen graph is quasi-$b$-invariant.)

A brick $G$ is near-bipartite if it has a pair of edges $\{e,f\}$ so that $G-e-f$
is matching covered and bipartite; such a pair $\{e,f\}$ is a removable doubleton.
(Each of $K_4$ and the triangular prism $\overline{C_6}$ has three removable
doubletons.)
Carvalho, Lucchesi and Murty (2002) proved a conjecture of Lov{\'a}sz which states
that every brick, distinct from $K_4$, $\overline{C_6}$ and the Petersen graph,
has a $b$-invariant edge.

A cubic graph is essentially $4$-edge-connected if it is $2$-edge-connected and if its
only $3$-cuts are the trivial ones; it is well-known that
each such graph is either a brick or a brace;
we provide a graph-theoretical proof of this fact.

We prove that if $G$ is any essentially $4$-edge-connected cubic brick then its edge-set
may be partitioned into three (possibly empty) sets: (i) edges that participate in a removable
doubleton, (ii) $b$-invariant edges, and (iii) quasi-$b$-invariant edges;
our Main Theorem
states that if $G$ has two adjacent quasi-$b$-invariant edges, say $e_1$~and~$e_2$,
then either $G$ is the Petersen graph or the (near-bipartite) Cubeplex graph,
or otherwise, each edge of $G$ (distinct from $e_1$ and $e_2$) is $b$-invariant.
As a corollary, we deduce that each essentially $4$-edge-connected cubic
non-near-bipartite brick $G$, distinct from the Petersen graph,
has at least $|V(G)|$ $b$-invariant edges.

\end{abstract}

\section{Matching covered graphs}
\label{sec:mcg}

  A graph  is {\it  matchable} if  it has  a perfect  matching.  Tutte
  \cite{tutt47}    proved    his   celebrated    \mbox{$1$-factor}    Theorem
  characterizing matchable graphs, and deduced  as a corollary that in
  a  \mbox{$2$-edge-connected}  cubic  graph  each edge  lies  in  a  perfect
  matching.

  \smallskip
  Let $G$ be a matchable graph.  A nonempty subset~$S$ of its vertices
  is a  {\it barrier} if  it satisfies the equation~${\sf  odd}(G-S) =
  |S|$, where ${\sf odd}(G-S)$ denotes the number of odd components of
  $G-S$.   For  distinct vertices  $u$~and~$v$  of~$G$,  it is  easily
  deduced from Tutte's Theorem that  the graph~$G-u-v$ is matchable if
  and only if no barrier of~$G$  contains both $u$ and $v$.  A barrier
  is {\it trivial} if it has a single vertex.

  \smallskip
  An edge~$e$  of $G$  is {\it  admissible} if  there is  some perfect
  matching   of~$G$   that   contains~$e$;  otherwise   it   is   {\it
    inadmissible}. Clearly, an  edge $e$ is admissible if  and only if
  no barrier of~$G$ contains both ends of~$e$.

  \smallskip
  A  connected  graph with  two  or  more  vertices is  {\it  \mc} if each of its edges is admissible.  The observation made
  above implies that a matchable graph  $G$ is \mc\ if and
  only if every barrier of $G$ is stable.
The following fundamental theorem is due to Kotzig;
see \cite[p. 150]{lopl86}.
\begin{theorem}
\label{thm:canonical-partition}
The maximal barriers of a \mcg\ partition
its vertex set.
\end{theorem}

The aforementioned corollary of Tutte's Theorem
may be rephrased as follows.
\begin{theorem}
\label{cor:cubic-mcg-iff-2ec}
A cubic graph is \mc\ if and only if it is $2$-edge-connected.
\end{theorem}

\subsection{Tight cut decompositions}
\label{sec:tight-cut-decompositions}

For a nonempty proper subset~$X$ of  the vertices of a graph~$G$, we
  denote by  $\partial(X)$ the cut  associated with~$X$, that  is, the
  set of all edges  of~$G$ that have one end in~$X$  and the other end
  in~$\overline{X}:=V(G)-X$.   We refer  to $X$~and~$\overline{X}$  as
  the {\it shores} of~$\partial(X)$. A cut  is {\it trivial} if either of
  its shores  is a  singleton.  We  say that  $\partial(X)$ is  a {\it
    $k$-cut} if $|\partial(X)|=k$.

\begin{figure}[!htb]
    \centering
    \subfigure[A barrier cut]
{
    \begin{tikzpicture}[scale=1.3]
      \draw (-0.5,-0.5)  -- (0,0) -- (0.5,-0.5)  -- (-0.5,-0.5); 
      \draw (-0.5,-0.5)node{} to [out=115,in=210] (0,1.5); 
      \draw (0,0)node{} to [out=90,in=200]   (1.5,1.5);  
      \draw (0.5,-0.5)node{} to [out=80,in=190] (3,1.5); 
      \draw[ultra  thick] (-1,0.8) to [out=0, in=120] (1.3,-0.5); 
      \draw (1.5,0) -- (0,1.5); 
      \draw (1.5,0) -- (1.5,1.5);  
      \draw  (1.5,0)node{}  --  (3,1.5);  
      \draw  (3,0)  -- (0,1.5)node{}; 
      \draw (3,0) -- (1.5,1.5)node{}; 
      \draw (3,0)node{} -- (3,1.5)node{};
    \end{tikzpicture}
\label{fig:barrier-cut}
}
\hspace*{1in}
\subfigure[A $2$-separation cut]
{
\begin{tikzpicture}[scale=0.65]
      \draw[ultra thick] (4.5,1.8) -- (1.5,-1.8);
      \draw (0,0) -- (1.5,0); 
      \draw (0,0)  -- (3,2); 
      \draw (1.5,0) -- (3,2);
      \draw  (0,0)node{} --  (3,-2);  
      \draw (1.5,0)node{}  -- (3,-2);  
      \draw (4.5,0)  -- (6,0);  
      \draw  (4.5,0) --  (3,2);  
      \draw (4.5,0)node{}  -- (3,-2); 
      \draw (6,0) -- (3,2)node{}; 
      \draw (6,0)node{} -- (3,-2)node{};
    \end{tikzpicture}
\label{fig:2-separation-cut}
}
    \caption{Nontrivial tight cuts}
  \end{figure}
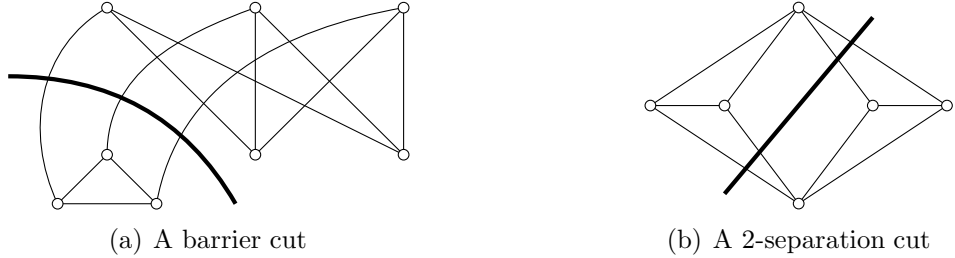

  \smallskip
  For a cut $\partial(X)$, we denote the graph obtained by contracting
  the  shore~$X$  to   a  single  vertex~$x$  by
  $G/(X \rightarrow  x)$.
The graph $G / (\overline{X} \rightarrow \overline{x})$ is defined analogously.
In case the  label of
  the contraction vertex $x$ or $\overline{x}$ is irrelevant, we simply write
  $G/ X$ or $G / \overline{X}$, respectively.
The  two graphs $G/X$ and  $G/ \overline{X}$ are
  called  the  \mbox{{\it   $\partial(X)$-contractions}}  of~$G$.   In
  Figure~\ref{fig:barrier-cut}, the  three edges crossing the  bold line
  constitute  a   nontrivial  cut,  say~$\partial(X)$,  and   the  two
  $\partial(X)$-contractions are $K_4$ and $K_{3,3}$.

%Let $G$ be a \mcg.
A cut $\partial(X)$ of a \mcg~$G$ is a {\it separating cut}
if each $\partial(X)$-contraction of $G$ is also \mc.
Clearly, each trivial cut is a separating cut.
The triangular prism~$\overline{C_6}$
has a (unique) nontrivial $3$-cut $\partial(X)$ that is a separating cut.
More generally, if $G$ is cubic, and if $\partial(X)$ is any $3$-cut,
then each $\partial(X)$-contraction of $G$ is a cubic graph that is $2$-edge-connected;
whence, by Theorem~\ref{cor:cubic-mcg-iff-2ec},
we have the following.

\begin{proposition}
\label{prp:cubic-mcg-3-cut-separating-cut}
In a cubic \mcg,
each $3$-cut is a separating cut. \qed
\end{proposition}

However, a cubic \mcg\ may have a separating cut
that is not a $3$-cut.
For instance, the Petersen graph has
nontrivial separating cuts, each of which is a $5$-cut; for any such cut $\partial(X)$,
each of the $\partial(X)$-contractions is isomorphic to the odd wheel~$W_5$.

\smallskip
A cut  $\partial(X)$ is a {\it tight cut} if $|M \cap  \partial(X)|=1$
for every perfect matching $M$  of~$G$. It is easily verified that every
tight cut is a separating cut.
The converse is not true. For instance, as noted earlier,
$\overline{C_6}$ has a nontrivial separating cut;
however, $\overline{C_6}$ is free of nontrivial
tight cuts.

  \smallskip
  For  instance,  if $B$  is  a  barrier of~$G$,  and  $K$  is an  odd
  component   of~$G-B$,  then   $\partial(V(K))$   is   a  tight   cut
  of~$G$. Such a  tight cut is called a {\it  barrier cut associated
  with the barrier} $B$, or simply a {\it barrier cut}.  The graph
  in Figure~\ref{fig:barrier-cut} has a barrier cut depicted by the bold
  line, and each of its contractions (that is, $K_4$ and $K_{3,3}$) is
  free of nontrivial tight cuts.

  \smallskip
  By a {\it $2$-separation}, we mean a $2$-vertex-cut.
  Now suppose that $\{u,v\}$ is a  $2$-separation of~$G$ that is not a
  barrier; that is,  each component of $G-u-v$ is even.   Let $K$ be a
  subgraph that is formed by the  union of some, but not all, components of $G-u-v$.
  Then each of the  sets $V(K) \cup \{u\}$ and $V(K)  \cup \{v\}$ is a
  shore of a nontrivial tight cut of~$G$. Such a tight cut is called a
  {\it     \mbox{$2$-separation}      cut associated with the $2$-separation}
  $\{u,v\}$, or simply a {\it \mbox{$2$-separation} cut}.
The      graph     in
  Figure~\ref{fig:2-separation-cut}  has  a  \mbox{$2$-separation}
  cut, and each of its contractions is $K_4$ with multiple edges.

\smallskip
Let $G$ be a \mcg.
If  $\partial(X)$ is  a
  nontrivial tight cut of~$G$,  then each $\partial(X)$-contraction is
  a \mcg\ that  has strictly fewer vertices than~$G$.
  If either  of the $\partial(X)$-contractions has  a nontrivial tight
  cut, then  that graph  can be further  decomposed into  even smaller
  {\mcg}s.  We  can repeat  this  procedure until  we
  obtain a list  of {\mcg}s, each of  which is free of
  nontrivial tight cuts.  This procedure is  known as a {\it tight cut
    decomposition} of~$G$.

  \smallskip
  Let  $G$  be a  \mcg\  free of  nontrivial  tight
  cuts. If $G$ is bipartite then it  is a {\it brace}; otherwise it is
  a {\it brick}.  Thus, a tight  cut decomposition of~$G$ results in a
  list of bricks and braces.

  \smallskip
  In general,  a \mcg\ may  admit several  tight cut
  decompositions.  However, Lov{\'a}sz  \cite{lova87} proved the following
  remarkable result.
  \begin{theorem}
    \label{thm:lovasz-tight-cut-decomposition}
    Any two tight cut decompositions of a \mcg\ yield
    the  same  list   of  bricks  and  braces   (except  possibly  for
    multiplicities of edges).
  \end{theorem}

  In  particular,  any two  tight  cut  decompositions of  a  \mcg~$G$
yield the same  number of bricks; this  number is
  denoted by $b(G)$.   We remark that $G$ is bipartite  if and only if
  $b(G) =  0$.

  \smallskip
  A  graph~$G$, with  four or  more vertices,  is {\it  bicritical} if
  $G-u-v$ is  matchable for  every pair of  distinct vertices  $u$ and
  $v$.       For      instance,      the      graph      shown      in
  Figure~\ref{fig:2-separation-cut}  is  bicritical,
whereas the graph shown in Figure~\ref{fig:barrier-cut} is not.   It  follows
  from Tutte's Theorem that a matchable graph $G$ is bicritical if and
  only if every barrier of $G$ is trivial.

  \smallskip
  Since a brick is a nonbipartite \mcg\ which is free
  of nontrivial  tight cuts,  it follows  from the  above observations
  that every  brick is $3$-connected and  bicritical.  Edmonds, Lov{\'a}sz
  and Pulleyblank \cite{elp82} established the converse.
  \begin{theorem}
    \label{thm:elp-brick-characterization}
    A  graph  is a  brick  if  and only  if  it  is $3$-connected  and
    bicritical.
  \end{theorem}

  In  fact, the  difficult  direction of  their  theorem statement  is
  equivalent to the following.
  \begin{theorem}
    \label{thm:existence-of-elp-cut}
    If a  \mcg\ has  a nontrivial tight cut,  then it
    has  a nontrivial  tight cut  that is  either a  barrier cut  or a
    $2$-separation cut.
  \end{theorem}

For cubic graphs, one may easily deduce the following
strengthening.
(See Lemma~\ref{lem:2-vertex-cut-with-cubic-vertex-implies-barrier-cut}.)
\begin{corollary}
\label{cor:existence-of-barrier-cut-in-cubic-mcgs}
If a cubic \mcg\ has a nontrivial tight cut,
then it has a nontrivial tight cut that is a barrier cut.
\end{corollary}

In general, a cubic \mcg\ need
  not be a brick  or a brace.  For instance, the  graph shown in
  Figure~\ref{fig:barrier-cut}  has a  nontrivial  tight  cut, and  this
  particular cut happens to be a  $3$-cut.
  In fact, this is not a coincidence.
  
  %Our first result shows that
  %this is not a coincidence.
  \begin{theorem}
    \label{thm:cubic-tight}
    In a cubic \mcg, each tight cut is a $3$-cut.
  \end{theorem}

  The above theorem may be proved easily using Edmond's characterization
  of the perfect matching polytope by considering the vector that
  assigns $\frac{1}{3}$ to each edge (see~\cite{kss09}).
  A     graph-theoretical proof     of     the      above     theorem     appears     in
  Section~\ref{sec:cubic-graphs-and-tight-cuts};
it  is   rather
  straightforward,  and  was  already  known   to  \mbox{C. N. Campos}  and
  \mbox{C. L. Lucchesi} in 1999.

\smallskip
We  say  that   a  cubic graph  is  {\it
    essentially $4$-edge-connected} if it is $2$-edge-connected
and if it is free of nontrivial  $3$-cuts.
(It is easy to see that such a graph is necessarily $3$-edge-connected
unless it is isomorphic to $C_4$ with multiple edges,
and that it is triangle-free unless it is isomorphic to $K_4$.)
The following is  an immediate
  consequence of Theorem~\ref{thm:cubic-tight}.
  \begin{corollary}
    \label{cor:efeccg-brick-brace}
    Every \efeccg\ is either a brick
    or a brace.
  \end{corollary}

\begin{figure}[!htb]
    \centering
    \subfigure[The Tricorn]
{
    \begin{tikzpicture}[scale=0.7]
      \draw[ultra thick] (70:3) -- (110:3);
      \draw[ultra thick] (190:3) -- (230:3);
      \draw[ultra thick] (310:3) -- (350:3);
      \draw (110:3) -- (190:3);
      \draw (230:3) -- (310:3);
      \draw (350:3) -- (70:3);
      \draw (0:0) -- (90:1.5);
      \draw (0:0) -- (210:1.5);
      \draw (0:0) node{} -- (330:1.5);
      \draw (90:1.5) -- (70:3)node{};
      \draw (110:3)node{} -- (90:1.5)node{};
      \draw (210:1.5) -- (190:3)node{};
      \draw (230:3)node{} -- (210:1.5)node{};
      \draw (330:1.5) -- (310:3)node{};
      \draw (350:3)node{} -- (330:1.5)node{};
    \end{tikzpicture}
\label{fig:Tricorn}
}
\hspace*{0.2in}
\vline
\hspace*{0.2in}
\subfigure[Staircases]
{
\begin{tikzpicture}[scale=0.7]
\draw (0,0) -- (1.5,1.5) -- (0,3) -- (0,0);
\draw (1.5,1.5)node{} -- (4.5,1.5)node{};
\draw (0,3)node{} -- (6, 3);
\draw (0,0)node{} -- (6,0);
\draw (6,0) -- (4.5,1.5)node{} -- (6,3)node{} -- (6,0)node{};
\draw[ultra thick] (3,1.5) -- (3,3);
\draw (3,1.5)node{};
\draw (3,3)node{};
\end{tikzpicture}
\hspace{0.3in}
\begin{tikzpicture}[scale=0.7]
\draw (0,0) -- (1.5,1.5) -- (0,3) -- (0,0);
\draw (1.5,1.5)node{} -- (4.5,1.5)node{};
\draw (0,3)node{} -- (6, 3);
\draw (0,0)node{} -- (6,0);
\draw (6,0) -- (4.5,1.5)node{} -- (6,3)node{} -- (6,0)node{};
\draw[ultra thick] (2.5,1.5) -- (2.5,3);
\draw[ultra thick] (3.5,1.5) -- (3.5,3);
\draw (2.5,1.5)node{};
\draw (2.5,3)node{};
\draw (3.5,1.5)node{};
\draw (3.5,3)node{};
\end{tikzpicture}
\label{fig:staircases}
}
    \caption{Some cubic bricks that are not \efec}
\label{fig:cubic-bricks-not-efec}
  \end{figure}
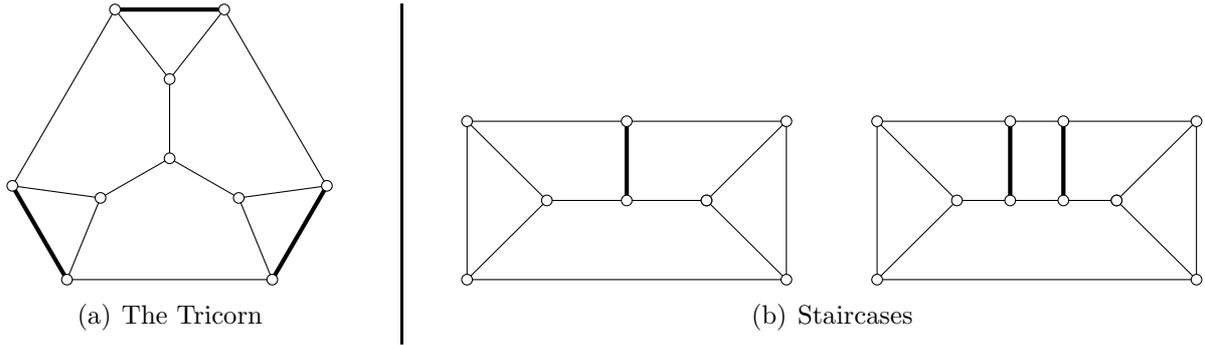

However, there exist cubic bricks that are not \efec.
For instance,
the `staircases' form one such infinite family; they play an important
role in \cite{noth08,komu16}.
The smallest staircase is $\overline{C_6}$, and the next two members
are shown in Figure~\ref{fig:staircases}.
Another example is the Tricorn; see Figure~\ref{fig:Tricorn}.

\smallskip
On the other hand, every cubic brace is in fact \efec.
This is due to Proposition~\ref{prp:cubic-mcg-3-cut-separating-cut},
and the fact that, in a bipartite \mcg,
every separating cut is also a tight cut; see~\cite[Corollary 2.22]{clm02}.

\smallskip
A \mcg~$G$ is {\it solid} if every separating cut of $G$
is a tight cut; otherwise $G$ is nonsolid.
The class of solid graphs is thus a generalization of bipartite graphs,
and it has played an important role in the theory of {\mcg}s;
see \cite{clm04,clm06,clm12,lckm18}. In particular, solid bricks
are precisely those bricks that are free of nontrivial separating cuts.
By Proposition~\ref{prp:cubic-mcg-3-cut-separating-cut},
every cubic solid brick is \efec.

\smallskip
Two infinite families of {\efeccg}s, worth mentioning here,
are the `prisms' and the `M{\"o}bius ladders'. See \cite{komu16} for definitions.
Each bipartite member of these families is a brace.
The nonbipartite M{\"o}bius ladders are solid bricks,
whereas the nonbipartite prisms are nonsolid bricks.

\subsection{$b$-invariant edges}
An edge  $e$ of a matching  covered graph $G$ is  {\it removable} if
  $G-e$ is also  \mc.  Furthermore, a  removable edge $e$
  is {\it \binv}  if $b(G-e) = b(G)$.  Note that,  if $G$ is bipartite
  then  any  removable  edge  $e$  is  \binv\  since  $b(G-e)=b(G)=0$.
  Furthermore, it can be easily shown that  if $G$ is a brace of order
  six or more,  then each edge is removable and  thus \binv.  However,
  the  notions  of  removability   and  $b$-invariance  are  far  more
  interesting, and nontrivial, in the case of bricks.

  \smallskip
  For  a matching  covered graph  $G$, an  edge $e$  {\it depends  on}
  another edge  $f$ if every  perfect matching that contains  $e$ also
  contains $f$, and  in this case the edge $f$  is not removable.  Two
  edges $e$ and $f$ are {\it mutually dependent} if $e$ depends on $f$
  and  $f$  depends on  $e$.   Lov{\'a}sz  \cite{lova87} proved the following.
\begin{proposition}
\label{prp:mutually-dependent}
If $\{e,f\}$ is a pair of mutually dependent edges in a brick~$G$
then $G-e-f$ is a matchable bipartite graph.
\end{proposition}

In particular, for a brick~$G$,
the complement of a pair of mutually dependent edges is a cut
of~$G$. Thus, if $\{e,f\}$ and $\{e,f'\}$ are distinct pairs of mutually
dependent edges then the symmetric difference of their complements
is a cut of~$G$. That is, $\{f,f'\}$ is a cut of~$G$. This is absurd since
bricks are $3$-edge-connected. This proves the following.
\begin{corollary}
\label{cor:mutually-dependent}
In a brick, any two distinct pairs of mutually dependent edges
are disjoint. \qed
\end{corollary}

In general, for a brick~$G$ and a pair of mutually dependent edges $\{e,f\}$,
the bipartite graph~$G-e-f$ need  not be  \mc.
We say that $R:=\{e,f\}$  is a
{\it removable doubleton} if  $G-R$ is \mc\ (and bipartite).
A brick $G$ is {\it  near-bipartite} if it has a removable
  doubleton; otherwise $G$ is {\it non-near-bipartite}.  For instance,
  the Petersen graph  is non-near-bipartite.  On the  other hand, each
  of  $K_4$   and  the triangular prism~$\overline{C_6}$
  has  three   distinct  removable
  doubletons; furthermore, each of them  is devoid of removable edges.
  Lov{\'a}sz \cite{lova87} proved that every brick distinct from
$K_4$~and~$\overline{C_6}$ has  a removable edge.

\smallskip
If $G$ is a brick and $e$ is a removable edge then $b(G-e) \geq 1$, and in general,
$b(G-e)$ can be arbitrarily large. A removable edge $e$ of a brick $G$
is \binv\ if and only if $b(G-e)=1$.
For instance, the Tricorn brick, shown in Figure~\ref{fig:Tricorn},
has precisely three removable edges (indicated by bold lines), each of which is \binv.
  
\smallskip
On the other hand, the Petersen graph is devoid of \binv\ edges despite
the fact that each edge is removable.
Confirming a  conjecture of
  Lov{\'a}sz, the  following result was  proved by Carvalho,  Lucchesi and
  Murty \cite{clm02a}.
  \begin{theorem}
    Every brick distinct from $K_4$, $\overline{C_6}$ and the Petersen
    graph has a \binv\ edge.
  \end{theorem}

As an application of the above theorem,
Carvalho, Lucchesi and Murty \cite{clm02b} gave an alternative
proof for the characterization of the matching lattice, and proved
other deep results concerning `ear decompositions' of {\mcg}s.
Since then, the existence of \binv\ edges in bricks, as well
as the existence of special types of \binv\ edges (such as `thin' and `strictly thin' edges)
in bricks and in braces,
have found many applications in matching theory;
see \cite{mccu04,clm06,clm12,komu16}.

\smallskip
In certain applications, it is helpful to have the presence of ``many'' \binv\ edges
--- for instance, one incident with each vertex.
An immediate consequence of Theorem 3.3 and Corollary 6.12 in \cite{clm12}
is the following.
\begin{theorem}
In a solid brick $G$, distinct from $K_4$, each vertex is incident with at
least one \binv\ edge; consequently, $G$ has at least $\frac{|V(G)|}{2}$
\binv\ edges.
\end{theorem}

In this paper, we establish a similar lower bound for another rich class of bricks.
\begin{theorem}
\label{thm:efeccnnbb-lower-bound}
%{\sc [Main Result]}
Every \efeccnnbb\
$G$, distinct from the Petersen graph,
has at least  $|V(G)|$ \binv\ edges.
\end{theorem}

We will prove a much stronger result that will immediately
imply Theorem~\ref{thm:efeccnnbb-lower-bound};
on the way there, we will discover other
noteworthy facts concerning the edges of an {\efeccb}. These are
discussed in the next section.

%new stuff added by Nishad starts here

\begin{figure}[!htb]
\centering
\subfigure[]{
\begin{tikzpicture}[scale=0.65]
	\draw (0,-6) -- (-1.5,-4);
	\draw[ultra thick] (0,-6) -- (0,-4);
	\draw (0,-6) -- (1.5,-4);
	\draw[ultra thick] (-1.5,-4) to [out=180,in=270] (-4,0);
	\draw (0,-4) to [out=135,in=270] (-2,0);
	\draw[ultra thick] (1.5,-4) to [out=0,in=270] (4,0);
	\draw (0,-4) to [out=45,in=270] (2,0)node{};
	\draw (-4,0) -- (-2,0);
	\draw (-1,4) to [out=180,in=90] (-4,0)node{};
	\draw[ultra thick] (1,4) to [out=210,in=90] (-2,0);
	\draw (-2,0) node{};
	\draw (4,0) -- (2,0);
	\draw (1,4) to [out=0,in=90] (4,0)node{};
	\draw[ultra thick] (-1,4) to [out=330,in=90] (2,0);
	\draw (2,0) node{};
	\draw[ultra thick] (0,2) -- (0,-2);
	\draw (0,2) -- (-1,4)node{};
	\draw (0,2) -- (1,4)node{};
	\draw (0,-2) -- (-1.5,-4)node{};
	\draw (0,-2) -- (1.5,-4)node{};
	\draw (0,2)node{};
	\draw (0,-2)node{};
	\draw (0,-4)node{};
	\draw (0,-6)node{};
      \end{tikzpicture}
\label{fig:Cubeplex-Mate-qbinv}
}
\hspace*{0.6in}
\subfigure[]{      
\begin{tikzpicture}[scale=0.75]
            
      \draw[ultra thick] (0,4) -- (0,1.5);
      \draw[ultra thick] (4,0) -- (1.5,0);
      \draw[ultra thick] (0,-4) -- (0,-1.5);
      \draw[ultra thick] (-4,0) -- (-1.5,0);
      
      \draw (0,4) -- (4,0) -- (0,-4) -- (-4,0) -- (0,4);
      \draw (0,1.5) -- (1.5,0) -- (0,-1.5) -- (-1.5,0) -- (0,1.5);
      
      \draw[ultra thick] (2,2) -- (-2,-2);
      \draw[ultra thick] (-2,2) -- (2,-2);
      
      \draw (0,4)node{};
      \draw (4,0)node{};
      \draw (0,-4)node{};
      \draw (-4,0)node{};
      
      \draw (0,1.5)node{};
      \draw (1.5,0)node{};
      \draw (0,-1.5)node{};
      \draw (-1.5,0)node{};
      
      \draw (2,2)node{};
      \draw (-2,2)node{};
      \draw (2,-2)node{};
      \draw (-2,-2)node{};
      
      \end{tikzpicture}
\label{fig:Cubeplex-Mate-Cube}
}
\caption{Two drawings of an \efeccnnbb~$G$ that has precisely $|V(G)|$ \binv\ edges}
\label{fig:Cubeplex-Mate}
\end{figure}
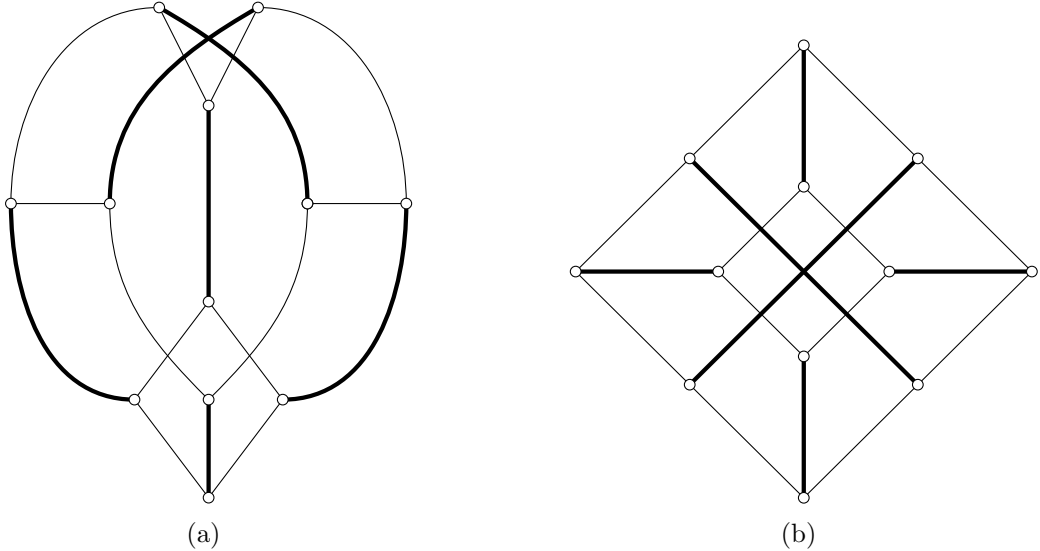

\smallskip
Figure~\ref{fig:Cubeplex-Mate} shows (two drawings of) an \efeccnnbb, of order~$12$,
that has precisely $12$ \binv\ edges
(and $6$ \qbinv\ edges that are indicated by bold lines);
however, apart from this graph, we do not know of any other graph
that meets the lower bound of Theorem~\ref{thm:efeccnnbb-lower-bound}.
In fact, as per the results of our recent computations,
each \efeccnnbb, of order at most $20$, has no more than $6$ \qbinv\ edges.
(Is this a coincidence?)

%new stuff added by Nishad ends here

\smallskip
It should be noted that
the lower bound of Theorem~\ref{thm:efeccnnbb-lower-bound}
does not hold for cubic bricks, in general.
For instance,
each staircase $G$, shown in Figure~\ref{fig:staircases},
has exactly $\frac{|V(G)|-6}{2}$ \binv\ edges;
these are near-bipartite but not \efec.
The Tricorn, shown in Figure~\ref{fig:Tricorn},
is neither \efec\ nor near-bipartite,
and it has exactly
three \binv\ edges.
(In Figure~\ref{fig:cubic-bricks-not-efec},
the \binv\ edges are indicated by bold lines.)

\smallskip
On the other hand, we
conjecture
that a weaker lower bound
holds for cubic bricks that are \efec\ and near-bipartite.
\begin{conjecture}
\label{con:efeccnbb-lower-bound}
Every \efec\ cubic near-bipartite brick $G$,
distinct from $K_4$, has at least $\frac{|V(G)|}{2}$
\binv\ edges.
\end{conjecture}

A proof of the above conjecture has already been
announced by Lu, Feng and Wang \cite{lfw19}.
Furthermore, they prove that prisms of order~$4k+2$,
and M{\"o}bius ladder of order~$4k$ (where $k \geq 2$)
are the only graphs that attain the conjectured lower bound.

%they prove that the only infinite families
%that attain the conjectured lower bound are:
%prisms of order $4k+2$, and M{\"o}bius ladders of order $4k$, where $k \geq 2$.
%See \cite{komu16} for definitions.

\subsection{Edges of an \efeccb}
The triangular prism $\overline{C_6}$ has a nontrivial $3$-cut $C$,
and for each $e \in C$,
the edge $e$ is neither removable nor does it participate in a removable doubleton.
We prove that this phenomenon cannot occur in {\efeccb}s.
\begin{theorem}
    \label{thm:efeccb-edges}
    In  an  \efeccb,  each  edge  is  either  removable  or  otherwise
    participates in a removable doubleton.
  \end{theorem}

It should be noted that, in a brick,
any edge can participate in at most one removable
doubleton. See Corollary~\ref{cor:mutually-dependent}.

\smallskip
As mentioned earlier,
for a removable edge $e$ of a brick $G$, the quantity $b(G-e)$
may be arbitrarily large. This is also true for cubic bricks, in general.
In Section~\ref{sec:infinite-family}, we describe how one may
construct such cubic bricks.
However, our next result shows that $b(G-e) \in \{1,2\}$
if the brick under consideration is cubic as well as \efec.

\smallskip
A removable edge $e$ of a brick is {\it \qbinv} if $b(G-e)=2$. 
For instance, each edge of the Petersen graph is \qbinv.
  \begin{theorem}
    \label{thm:efeccb-removable-edges}
    In an \efeccb,  each removable edge is either  \binv\ or otherwise
    \qbinv.
  \end{theorem}

The above theorem is reminiscent of the following result of Carvalho, Lucchesi and
Murty~\cite{clm02}.
\begin{theorem}
\label{thm:solid-removable-implies-binv}
In a solid brick, each removable edge is \binv.
\end{theorem}

It follows from
Theorems~\ref{thm:efeccb-edges}~and~\ref{thm:efeccb-removable-edges}
that the edge set of an \efeccb\ may be partitioned into three
disjoint (possibly empty) sets: (i) edges that participate in a removable
doubleton --- these come in pairs, (ii) edges that are \binv,
and (iii) edges that are \qbinv.

\smallskip
As mentioned earlier, it is often helpful to have the presence of ``many''
\binv\ edges. On a related note, a brick being near-bipartite is most likely
``good'' news. For instance, while there has been no significant progress
in characterizing `Pfaffian' bricks, Fischer and Little \cite{fili01}
were able to characterize Pfaffian near-bipartite bricks.
(See Section~\ref{sec:Pfaffian-and-conformal} for definition
and further discussion on this topic.)
In this sense,
the only ``unpleasant'' outcome for any particular edge $e$
of an \efeccb~$G$ is that it happens to be \qbinv.
With this in mind, we started wondering whether we could prove
an upper bound on the number of \qbinv\ edges
in an \efeccb. Such a result would yield a lower bound on the number
of \binv\ edges --- in the case that $G$ is non-near-bipartite.

\smallskip
As noted earlier, in the Petersen graph, every edge is \qbinv.
In particular, if $v$ is any vertex, then all three edges incident with $v$
are \qbinv. Our next result shows that, among the {\efeccb}s, this
phenomenon is unique to the Petersen graph. In fact,
we prove the following stronger result that immediately implies
Theorem~\ref{thm:efeccnnbb-lower-bound}
(see Section~\ref{sec:main-theorem-consequences}).
  \begin{theorem}
    \label{thm:efeccb-two-qbinv-at-vertex}
{\sc [Main Theorem]}
Let $G$ be an \efeccb\ that has two adjacent \qbinv\ edges $e_1$ and $e_2$.
Then the following statements hold:
    \begin{enumerate}[(i)]
    \item
    \label{itm:bricks-J-K4}
             For $i \in \{1,2\}$, both bricks of $G-e_i$ are isomorphic to $K_4$ (up to multiple edges).
    \item
    \label{itm:G-nonplanar}
            The graph $G$ is nonsolid, nonplanar and non-Pfaffian.
    \item
    \label{itm:near-bipartite-implies-Cubeplex}
            If $G$ is near-bipartite then $G$ is the Cubeplex.
    \item
    \label{itm:any-more-qbinv-implies-Petersen}
            If $G$ has a \qbinv\ edge, distinct from $e_1$ and $e_2$,
            then $G$ is the Petersen graph.
    \end{enumerate}
Consequently, every edge of $G$, distinct from $e_1$ and $e_2$,
is \binv, \underline{unless} $G$ is either the Cubeplex or the Petersen graph.
  \end{theorem}

\begin{figure}[!htb]
    \centering
    \begin{tikzpicture}[scale=0.65]

      \draw (-5.5,0) -- (-3,-4);
      \draw (-3.5,0) -- (-1,-4);

      \draw (5.5,0) -- (1,-4);
      \draw (3.5,0) -- (3,-4);

      \draw (-1,-6) -- (-1,-4);
      \draw (-1,-6) -- (-3,-4);
      \draw (-1,-6) -- (3,-4);

      \draw (1,-6) -- (1,-4);
      \draw (1,-6) -- (-3,-4);
      \draw (1,-6) -- (3,-4);

      \draw (-1,-6)node{};
      \draw (1,-6)node{};
      \draw (3,-4)node{};
      \draw (-3,-4)node{};

      \draw (-5.5,0) -- (-3.5,0);
      \draw (-1,4) to [out=180,in=90] (-5.5,0)node{};

      \draw (1,4) to [out=210,in=90] (-3.5,0)node{};

      \draw (5.5,0) -- (3.5,0);
      \draw (1,4) to [out=0,in=90] (5.5,0)node{};

      \draw (-1,4) to [out=330,in=90] (3.5,0)node{};
      
      \draw (0,2) -- (0,-2);
      
      \draw (0,2) -- (-1,4)node{};
      \draw (0,2) -- (1,4)node{};

      \draw (0,-2) -- (-1,-4)node{};
      
      \draw (0,-2) -- (1,-4)node{};
      
      \draw (-1,4.5)node[nodelabel]{$u_2$};
      \draw (1,4.5)node[nodelabel]{$u_3$};

      \draw (0,2)node{};
      \draw (0.3,1.8)node[nodelabel]{$v$};
      \draw (0,-2)node{};
      \draw (0.4,-1.8)node[nodelabel]{$u_1$};

    \end{tikzpicture}
    \caption{The edges $vu_1$ and $vu_2$ are \qbinv\ whereas $vu_3$ is
      \binv}
    \label{fig:two-qbinv-one-binv-smallest}
  \end{figure}
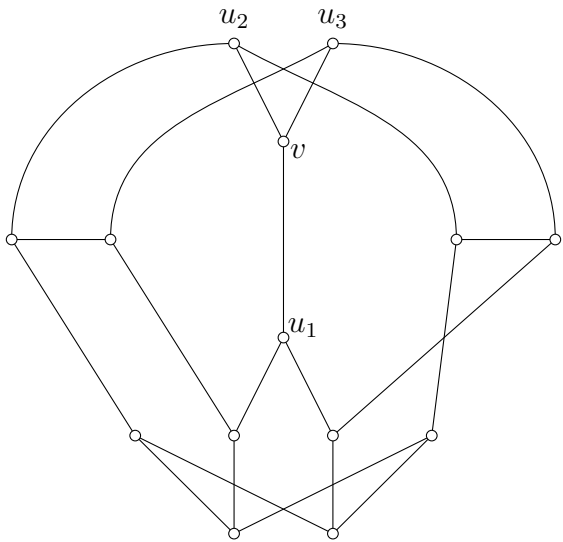

\begin{figure}[!htb]
\centering
\subfigure[]{
\begin{tikzpicture}[scale=0.65]

	%edges incident with y'
	\draw (0,-6) -- (-1.5,-4);
	\draw (0,-6) -- (0,-4);
	\draw (0,-6) -- (1.5,-4);

	%edges d' and f'
	\draw (-2,0) to [out=270,in=150] (-1.5,-4);
	\draw (-1.6,-1.2)node[nodelabel]{$d'_1$};
	\draw (-4,0) to [out=270,in=135] (0,-4);
	\draw (-4,-1.6)node[nodelabel]{$f'_1$};

	%edges g' and h'
	\draw (1.5,-4) to [out=0,in=270] (4,0);
	\draw (4.1,-2)node[nodelabel]{$g'_1$};
	\draw (0,-4) to [out=45,in=270] (2,0)node{};
	\draw (2.1,-1.9)node[nodelabel]{$h'_1$};

	%edges d and f
	\draw (-4,0) -- (-2,0);
	\draw (-1,4) to [out=180,in=90] (-4,0)node{};
	\draw (-4.5,0)node[nodelabel]{$s_2$};
	\draw (-3.9,2.3)node[nodelabel]{$d_1$};
	\draw (1,4) to [out=210,in=90] (-2,0)node[fill=black]{};
	\draw (-1.5,0)node[nodelabel]{$s_3$};
	\draw (-2,1.9)node[nodelabel]{$f_1$};

	%edges g and h
	\draw (4,0) -- (2,0);
	\draw (1,4) to [out=0,in=90] (4,0)node[fill=black]{};
	\draw (4.5,0)node[nodelabel]{$t_3$};
	\draw (4,2.3)node[nodelabel]{$h_1$};
	\draw (-1,4) to [out=330,in=90] (2,0)node{};
	\draw (1.5,0)node[nodelabel]{$t_2$};
	\draw (2,1.9)node[nodelabel]{$g_1$};

	%edge e=uv
	\draw (0,2) -- (0,-2);
	\draw (0.4,0)node[nodelabel]{$e_1$};

	%other edges incident at v
	\draw (0,2) -- (-1,4)node[fill=black]{};
	\draw (0,2) -- (1,4)node{};

	%other edges incident at u
	\draw (0,-2) -- (-1.5,-4)node{};
	\draw (-1.7,-4.5)node[nodelabel]{$s_1$};
	\draw (0,-2) -- (1.5,-4)node{};
	\draw (1.7,-4.5)node[nodelabel]{$t_1$};

	%barrier B
	\draw (-1,4.5)node[nodelabel]{$u_2$};
	\draw (1,4.5)node[nodelabel]{$u_3$};

	%ends u and v of edge e
	\draw (0,2)node{};
	\draw (0.3,1.8)node[nodelabel]{$v$};
	\draw (0,-2)node[fill=black]{};
	\draw (0.4,-1.8)node[nodelabel]{$u_1$};

	%vertices y and y'
	\draw (0,-4)node[fill=black]{};
	\draw (0,-3.5)node[nodelabel]{$y_1$};
	\draw (0,-6)node[fill=black]{};
	\draw (0,-6.5)node[nodelabel]{$y_2$};

      \end{tikzpicture}
\label{fig:Cubeplex-qbinv}
}
\hspace*{0.2in}
\subfigure[]{
      \begin{tikzpicture}[scale=0.75]
      
      \draw (0,0) to [out=75,in=105] (2.75,0);
      \draw (0,0) to [out=105,in=75] (-2.75,0);
      \draw (0,0) to [out=195,in=165] (0,-2.75);
      
      \draw (0,4) -- (0,1.5);
      \draw (4,0) -- (1.5,0);
      \draw (0,-4) -- (0,-1.5);
      \draw (-4,0) -- (-1.5,0);
      
      \draw (0,4) -- (4,0) -- (0,-4) -- (-4,0) -- (0,4);
      \draw (0,1.5) -- (1.5,0) -- (0,-1.5) -- (-1.5,0) -- (0,1.5);
      
      \draw (0,4)node[fill=black]{};
      \draw (0.4,4.1)node[nodelabel]{$y_1$};%label
      \draw (4,0)node{};
      \draw (4,-0.4)node[nodelabel]{$t_2$};%label
      \draw (0,-4)node[fill=black]{};
      \draw (0.4,-4.1)node[nodelabel]{$u_2$};%label
      \draw (-4,0)node{};
      \draw (-4,-0.4)node[nodelabel]{$s_2$};%label
      
      \draw (0,1.5)node[fill=black]{};
      \draw (0.4,1.6)node[nodelabel]{$y_2$};
      \draw (1.5,0)node{};
      \draw (1.5,-0.4)node[nodelabel]{$t_1$};%label
      \draw (0,-1.5)node[fill=black]{};
      \draw (0.4,-1.6)node[nodelabel]{$u_1$};%label
      \draw (-1.5,0)node{};
      \draw (-1.5,-0.4)node[nodelabel]{$s_1$};%label

      \draw (2.75,0)node[fill=black]{};
      \draw (2.75,-0.4)node[nodelabel]{$t_3$};%label
      \draw (-2.75,0)node[fill=black]{};
      \draw (-2.75,-0.4)node[nodelabel]{$s_3$};%label
      \draw (0,-2.75)node{};
      \draw (0.4,-2.75)node[nodelabel]{$v$};%label
      
      \draw (0,0)node{};
      \draw (0.4,-0.1)node[nodelabel]{$u_3$};%label
      
      \end{tikzpicture}
\label{fig:Cubeplex-Cube}
}
    \caption{Two drawings of the Cubeplex}
\label{fig:Cubeplex}
\end{figure}
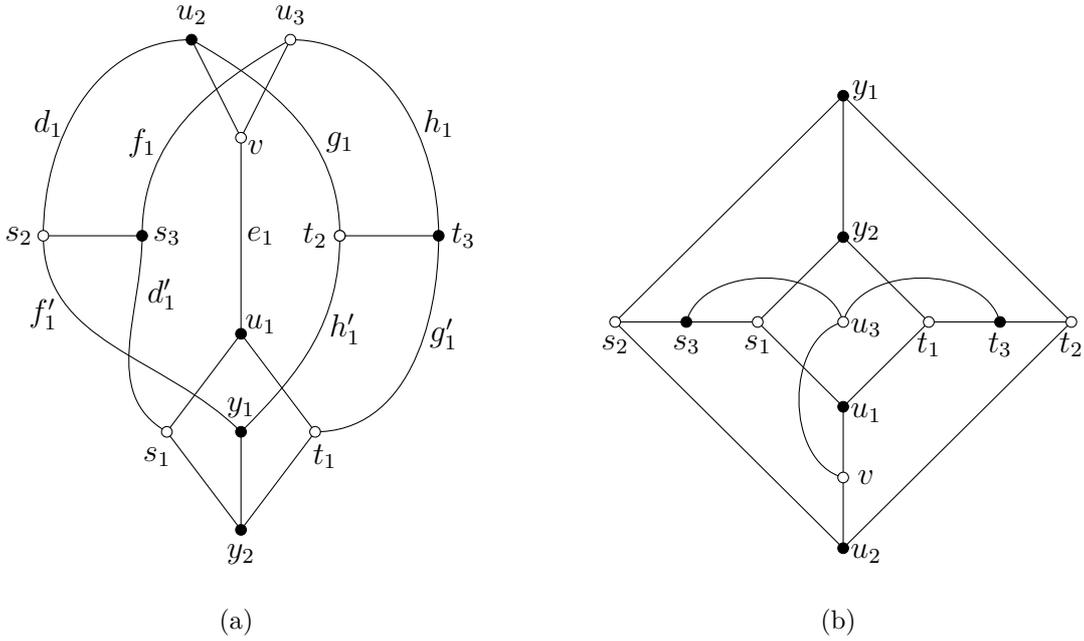

The proof appears     in
  Section~\ref{sec:main-theorem-proof}.
  Figures~\ref{fig:Petersen}~and~\ref{fig:Cubeplex}
  show  the two  special graphs
  that  appear  in Theorem~\ref{thm:efeccb-two-qbinv-at-vertex}.   The
  Cubeplex  first  appeared  in  the   work  of  Fischer  and  Little
  \cite{fili01}  where they  showed that  it is  one of  two minimally
  non-Pfaffian near-bipartite graphs; they  referred to the two graphs
  as $\Gamma_1$  and $\Gamma_2$.   The names Cubeplex~for~$\Gamma_1$,
  and  Twinplex~for~$\Gamma_2$,  are  due  to   Norine  and  Thomas
  \cite{noth08}.
(See Section~\ref{sec:Pfaffian-and-conformal} for further discussion.)
Figure~\ref{fig:two-qbinv-one-binv-smallest}  shows
  the smallest \efeccb\ that has  a vertex $v$ incident with precisely
  two \qbinv\  edges and one \binv\  edge. The choice of  the specific
  drawings in
  Figures~\ref{fig:Cubeplex-Mate-qbinv}, \ref{fig:two-qbinv-one-binv-smallest}, \ref{fig:Petersen}
  and \ref{fig:Cubeplex-qbinv}
  is  by  no means  coincidental; it is  related  to the  proofs  of
  Theorems~\ref{thm:efeccb-removable-edges}~and~\ref{thm:efeccb-two-qbinv-at-vertex}.

\smallskip
Throughout this research, we made extensive use of computations ---
especially, in discovering the statement and proof of the
Main Theorem~(\ref{thm:efeccb-two-qbinv-at-vertex}).
To this end, we downloaded the exhaustive lists of cubic graphs from
the House of Graphs~\cite{bcgm13}, filtered the essentially $4$-edge-connected bricks,
and performed various computations using SageMath~\cite{sagemath}.

\smallskip
It should be noted that the class of {\efeccg}s has been studied in
the literature; however, perhaps not from the point of view
of tight cuts, bricks and braces, and \binv\ edges.
We mention a particular work of Wormald \cite{worm79},
wherein he proves a generation theorem for 
this class of graphs --- that he refers to as
`cyclically $4$-connected cubic graphs'.
Also, see Bondy and Murty \cite[Exercise 9.4.7]{bomu08}.
It would be worth investigating whether
Wormald's generation theorem can be used to obtain
simpler proofs for any of our results --- especially, for the Main
Theorem (\ref{thm:efeccb-two-qbinv-at-vertex}).

\subsection{Organization and summary of this paper}

In Section~\ref{sec:cubic-graphs-and-tight-cuts},
we prove Theorem~\ref{thm:cubic-tight}, which implies that
every \efeccg~$G$ is either a brick or a brace. The rest of this paper
deals with the case in which $G$ is a brick.
In Section~\ref{sec:removability},
we prove Theorem~\ref{thm:efeccb-edges} which states that
each edge is either removable or otherwise participates in a removable
doubleton.
In Section~\ref{sec:binv-and-qbinv},
we prove Theorem~\ref{thm:efeccb-removable-edges} which states that
if $e$ is a removable edge, then $b(G-e) \in \{1,2\}$;
and in Section~\ref{sec:infinite-family},
we demonstrate why such a result does not hold for cubic bricks, in general.
The Main Theorem (\ref{thm:efeccb-two-qbinv-at-vertex})
is proved in Section~\ref{sec:main-theorem-proof},
and it immediately implies Theorem~\ref{thm:efeccnnbb-lower-bound}
which states that if $G$ is non-near-bipartite
then at least two-third of its edges are \binv\ \underline{unless}
$G$ is the Petersen graph.

\smallskip
Our Main Theorem~(\ref{thm:efeccb-two-qbinv-at-vertex})
has some consequences pertaining to `Pfaffian orientations'
and the related notion of `conformal minors'. These
topics are discussed in Section~\ref{sec:Pfaffian-and-conformal},
and the relevant consequences of the Main Theorem
are stated in Section~\ref{sec:main-theorem-consequences}.

\section{Cubic graphs and tight cuts}
\label{sec:cubic-graphs-and-tight-cuts}
In this section, our goal is to provide a graph-theoretical proof of
Theorem~\ref{thm:cubic-tight}.
First, we need some easy facts pertaining to tight cuts.

\smallskip
Two cuts $C:=\partial(X)$ and $D:=\partial(Y)$ are said to be {\it crossing} if all four sets $X \cap Y$, $\overline{X} \cap Y$,
$X \cap \overline{Y}$ and $\overline{X} \cap \overline{Y}$ are nonempty; otherwise $C$ and $D$ are said to be {\it laminar}.
The following lemma is useful in proving theorems concerning {\mcg}s, and was used by Lov{\'a}sz \cite{lova87} in
the proof of Theorem~\ref{thm:lovasz-tight-cut-decomposition}.
\begin{lemma}
\label{lem:uncrossing-tight-cuts}
Let $G$ be a \mcg, and let $C:=\partial(X)$ and $D:=\partial(Y)$ be crossing tight cuts such that $|X \cap Y|$ is odd.
Then:
\begin{enumerate}[(i)]
\item $I:=\partial(X \cap Y)$ and $U:=\partial(\overline{X} \cap \overline{Y})$ are both tight cuts,
\item there are no edges between $\overline{X} \cap Y$ and $X \cap \overline{Y}$, and
\item $|C| + |D| = |I| + |U|$. \qed
\end{enumerate}
\end{lemma}

The following lemma immediately implies
Corollary~\ref{cor:existence-of-barrier-cut-in-cubic-mcgs}.
\begin{lemma}
\label{lem:2-vertex-cut-with-cubic-vertex-implies-barrier-cut}
Let $G$ be a \mcg\ that has a $2$-separation $\{u,v\}$.
If either of $u$ and $v$ is a cubic vertex then $G$ has a
barrier of cardinality two.
\end{lemma}
\begin{proof}
As in the statement,
let $\{u,v\}$ be a $2$-separation of a \mcg~$G$,
and assume that $v$ is a cubic vertex.
If $G-u-v$ has an odd component then $\{u,v\}$ is a barrier.
Now suppose that each component of $G-u-v$ is even.
Since $G$ is $2$-connected, and since $v$ is cubic,
there exists a component $L$ of $G-u-v$
such that $v$ has exactly one neighbour, say $w$, that lies in $L$.
Observe that $\{u,w\}$ is a barrier of $G$.
\end{proof}

We state an immediate consequence that will be useful later.
(See Theorem~\ref{thm:elp-brick-characterization}.)
\begin{corollary}
\label{cor:bicritical-not-brick}
Let $G$ be a bicritical graph. Then, for every $2$-separation $\{u,v\}$,
each of $u$~and~$v$ is a non-cubic vertex. (In particular, if $G$ is not a brick
then $G$ has at least two non-cubic vertices.) \qed
\end{corollary}

The graph, shown in Figure~\ref{fig:2-separation-cut}, is the smallest bicritical graph
that is not a brick, and it has exactly two noncubic vertices.

\smallskip
The next fact is easily verified.
\begin{lemma}
\label{lem:parity}
Let $G$ be a cubic graph, and let $S \subseteq V(G)$.
Then $|\partial(S)| \equiv |S| {\rm ~mod~} 2$.
Furthermore, if $G$ is $2$-edge-connected,
then $|\partial(S)| \geq 2$ whenever $S \neq \emptyset$,
and $|\partial(S)| \geq 3$ whenever $|S|$ is odd. \qed
\end{lemma}

\medskip
\noindent
\textbf{Theorem~\ref{thm:cubic-tight}}
\textit{In a cubic \mcg, each tight cut is a $3$-cut.}

\smallskip
\noindent
\begin{proof}
    Let  $G$  be  a  cubic  \mcg.   We  proceed  by
    induction on $|V(G)|$.   Every trivial cut  is a $3$-cut.
    Now    let    $C$    denote   a    nontrivial    tight    cut.

    \smallskip
    By Corollary~\ref{cor:existence-of-barrier-cut-in-cubic-mcgs},
    $G$ has a nontrivial tight cut $D$ that is a barrier cut
    associated with some nontrivial barrier, say~$B$.
    Since $G$ is cubic  and since
    $B$ is a  stable set, $|\partial(B)|=3|B|$.  Also,  for each (odd)
    component~$J$  of $G-B$, we  have that $|\partial(V(J))|  \geq 3$.
    We  infer  that $|\partial(V(J))|=3$  for  each  component $J$  of
    $G-B$. Thus $|D|=3$.
    In particular, every barrier cut is indeed a $3$-cut.
    We may thus assume that $C$ and $D$ are distinct cuts.

    \smallskip
    It remains
    to deduce that $|C|=3$.
    We consider two cases depending on whether $C$~and~$D$ are laminar
    cuts, or whether they are crossing cuts.

    \smallskip
    First suppose  that $C$ and  $D$ are  laminar cuts. Let  $G_1$ and
    $G_2$ be the  two $D$-contractions of $G$, and  adjust notation so
    that $C$ is  a tight cut of $G_1$. Since  $|D|=3$, it follows that
    $G_1$  is a  cubic  \mcg,  whence by  induction
    hypothesis $|C|=3$.

    \smallskip
    Now  suppose  that $C$  and  $D$  are  crossing  cuts. As  in  the
    statement  of   Lemma~\ref{lem:uncrossing-tight-cuts},  we  adjust
    notation so  that $C:=\partial(X)$  and $D:=\partial(Y)$  and that
    $|X \cap  Y|$ is  odd. Then  each of  $I:=\partial(X \cap  Y)$ and
    $U:=\partial(\overline{X} \cap \overline{Y})$ is  a tight cut that
    is  laminar with  the cut  $D$.  From  the preceding  paragraph we
    infer that $|I|=|U|=3$.  By Lemma~\ref{lem:uncrossing-tight-cuts},
    we have that $|C|+|D|=|I|+|U|$, and thus $|C|=3$.
\end{proof}

\section{Essentially $4$-edge-connected cubic bricks}
In the last section, we established that every \efeccg\
is either a brick or a brace. Now we focus on the nonbipartite members
of this class: bricks.

\subsection{Removability}
  \label{sec:removability}
In this section, we present a proof of Theorem~\ref{thm:efeccb-edges}.
We will use the following fact that is easily verified.
\begin{lemma}
\label{lem:3-conn-3-cut}
In a $3$-connected graph, every nontrivial $3$-cut is a matching. \qed
\end{lemma}

In our attempts to prove
the Main Theorem (\ref{thm:efeccb-two-qbinv-at-vertex}),
we ran into a problem wherein
we had to deal with a slightly more general class of cubic bricks
(that includes the class of {\efeccb}s).
It is for this reason that we    prove    the    following   generalization    of
Theorem~\ref{thm:efeccb-edges}.
\begin{theorem}
\label{thm:efeccb-edges-generalization}
Let $G$ be a cubic brick, and let $e$ denote an edge that participates
in every nontrivial $3$-cut of $G$ (if such cuts exist). Let $f$ denote any edge
that is adjacent with~$e$. Then:
\begin{enumerate}[(i)]
\item either $f$ is removable, or
\item or otherwise there exists a unique edge $f'$ that depends on $f$,
and $\{f,f'\}$ is a removable doubleton of $G$.
\end{enumerate}
\end{theorem}

In order to see why the above theorem
implies Theorem~\ref{thm:efeccb-edges},
it suffices to note that if $G$ is \efec\ then any edge may play the role of $e$
in the statement of the above theorem.

\smallskip
Observe that, a staircase, as shown in Figure~\ref{fig:staircases},
is a cubic brick
that is not \efec. However, each staircase has an edge that participates
in every nontrivial $3$-cut. On the other hand, the Tricorn,
shown in Figure~\ref{fig:Tricorn}, has
no such edge.

\begin{proof}[Proof of Theorem~$\ref{thm:efeccb-edges-generalization}$]
Let $v$ denote the common end of edges $e$ and $f$.    
Suppose that $f$ is not removable in $G$.
That is, some edge depends on $f$.
Our goal is to deduce that $f$ participates in a removable doubleton.

\smallskip
Let $f'$ denote any edge that depends on~$f$. In what follows,
we will show that $\{f,f'\}$
is a pair of mutually dependent edges;
thus, by Corollary~\ref{cor:mutually-dependent},
$f'$ is the unique edge that depends on~$f$.
Furthermore, since the dependence relation is transitive,
no edge in the set $E(G)-f-f'$
depends on either of $f$ and $f'$.
Hence, $G-f-f'$ is \mc\ and $\{f,f'\}$
is a removable doubleton of $G$.
(Now, it remains to show that $f$ depends on $f'$.)

\smallskip
Since $f'$ depends on $f$, the edge~$f'$ is inadmissible in the matchable graph $G-f$.
Applying Tutte's 1-factor Theorem,
$G-f$ has a barrier $B$ that contains both ends of $f'$.
Let $\mathcal{J}$ denote the set of odd components of $G-f-B$.
Since $G$ itself is free of nontrivial barriers, $f$ must have its
ends in distinct members of $\mathcal{J}$, say $J_1$ and $J_2$.
Adjust notation so that $v \in V(J_1)$.
Observe that, for each $J \in \mathcal{J}-J_1$, the edge $e$
does not lie in the cut $\partial(V(J))$.

\smallskip
Since $G$ is $3$-edge-connected, for each $J \in \mathcal{J}$,
we have $|\partial(V(J))| \geq 3$.
Since $f$ has its ends in $J_1$ and in $J_2$,
there are at least $3|B|-2$ edges
that have one end in a member of~$\mathcal{J}$, and the other end in $B$.
Since $G$ is cubic, and since $f'$ has both ends in $B$,
we infer that $|\partial(B)| \leq 3|B|-2$.
Thus there are exactly $3|B|-2$ edges that have one end in a member of $\mathcal{J}$,
and the other end in $B$. Consequently,
for each $J \in \mathcal{J}$, we have $|\partial(V(J))|=3$.
Furthermore, $G-f-B$ has no even components.
Since $e$ participates in every nontrivial $3$-cut of $G$,
we conclude that every component of $G-f-B$, except
perhaps $J_1$, is trivial. We will now argue that $J_1$ is also trivial.

\smallskip
Suppose to the contrary that $J_1$ is nontrivial,
whence $\partial(V(J_1))$ is a nontrivial $3$-cut of~$G$.
Thus $e \in \partial(V(J_1))$.
Note that $f \in \partial(V(J_1))$.
Since $e$ and $f$ are adjacent,
$\partial(V(J_1))$ is a nontrivial $3$-cut that is not a matching.
This contradicts
Lemma~\ref{lem:3-conn-3-cut}.

\smallskip
In summary, each member of $\mathcal{J}$ is trivial.
Thus $G-f-f'$ is bipartite; one of its color classes is $B$ which contains
both ends of $f'$; the other color class contains both ends of $f$.
Clearly, $f$ depends on $f'$.
As discussed earlier, this completes
the proof of Theorem~\ref{thm:efeccb-edges-generalization}.
\end{proof}

The following is an
easy application of Theorem~\ref{thm:efeccb-edges-generalization},
and it will be useful to us
in the proof of the Main Theorem~(\ref{thm:efeccb-two-qbinv-at-vertex}).

\begin{corollary}
\label{cor:dependence-implies-doubleton}
Let $J$ be a cubic brick, and let $xx'$ denote an edge
that participates in every nontrivial $3$-cut of $J$ (if such cuts exist).
Let $\partial(x):=\{e,d,f\}$ and let $\partial(x'):=\{e,d',f'\}$.
If either of $d$ and $d'$ depends on the other,
or if either of $f$ and $f'$ depends on the other,
then each of $\{d,d'\}$ and $\{f,f'\}$ is a removable doubleton of $J$. \qed
\end{corollary}

\subsection{$b$-invariance and quasi-$b$-invariance}
\label{sec:binv-and-qbinv}

  In this section, we will prove Theorem~\ref{thm:efeccb-removable-edges}
  which states that every removable edge of an \efeccb~$G$ is either \binv\
  or \qbinv.
  In fact, we will prove a stronger result that also describes
  the structure of $G$ with respect to any given \qbinv\ edge.
  This will help us in proving the Main Theorem (\ref{thm:efeccb-two-qbinv-at-vertex}).
  
\smallskip
Before that, we need a  few   preliminary results.
The following two statements are self-evident.
\begin{lemma}
\label{lem:6-cut}
Let $G$ be an \efeccg, and let $\partial(X)$ be a $6$-cut.
If $G[X]$ is disconnected then $|X|=2$.~\qed
\end{lemma}
\begin{lemma}
    \label{lem:4-cut}
    Let  $\partial(X)$ be a $4$-cut of an \efeccg~$G$.
    If there exist two edges in $\partial(X)$ that have a  common end in
    $X$ then $G[X]$ is isomorphic to $K_2$.~\qed
\end{lemma}

Recall that the tight cut decomposition of a bipartite \mcg\
yields only braces. In the same spirit, the following result from
Lov{\'a}sz and Plummer \cite[chapter 5]{lopl86},
implies that the tight cut decomposition
of a bicritical graph yields only bricks.
\begin{proposition}
\label{prp:bicritical-tight-cut}
Let $G$ be a bicritical graph, and let $C$ denote a $2$-separation cut.
Then each $C$-contraction of $G$ is also bicritical. \qed
\end{proposition}

A barrier $B$ of a \mcg~$G$ is
  \emph{special} if $G-B$ has precisely one nontrivial component.
\begin{lemma}
    \label{lem:every-barrier-special}
    Let $e$ be a removable edge  of an \efeccb\ $G$, and let $B$ denote
    a nontrivial barrier of $G-e$.  Then $B$ is special.
    Furthermore, $|\partial(V(J))|=5$,
    where $J$ denotes the unique nontrivial component of $G-e-B$.
\end{lemma}
\begin{proof}
Since $e$ is removable, $G-e$ is \mc.
Thus $B$ is stable, and $G-e-B$ has no even components.
Since $G$ is cubic, $|\partial(B)|=3|B|$.

\smallskip
We let $\mathcal{J}$ denote the set of components of $G-e-B$.
For each $J \in \mathcal{J}$,
it follows from Lemma~\ref{lem:parity} that
$|\partial(V(J))|$ is odd, and is at least $3$.
Also, since $G$ itself is free of nontrivial barriers, $e$ has its ends
in distinct members of $\mathcal{J}$.
From these facts, we infer that there exists a unique $J' \in \mathcal{J}$
such that $|\partial(V(J'))| =5$, and for each $J \in \mathcal{J} - J'$,
we have $|\partial(V(J))|=3$. Since $G$ is \efec, $J'$
is the only nontrivial component of $G-e-B$.
In particular, $B$ is special, and this completes the proof
of Lemma~\ref{lem:every-barrier-special}.
\end{proof}

For  a special barrier $B$  of $G-e$, we let  $I$ denote
  the  set of  isolated vertices  of  $G-e-B$, and
  we let $X:=B \cup I$.
  Note that, since $B$ is special, $|I| = |B|-1$.
Often we  may use subscripts  or superscripts,  or both, to  denote a
  special barrier  --- for instance, $B'_1$. In this  case, the
  corresponding set  of isolated vertices
  will  be decorated similarly --- that is, $I'_1$ --- and
  likewise $X'_1 := B'_1 \cup  I'_1$.
\begin{lemma}
\label{lem:bipartite-shore-properties}
Let $G$ be a bicritical graph, let $e:=uv$ be a removable edge of $G$
and let $B$ denote a nontrivial special barrier of $G-e$.
Assume that $e$ has exactly one end, say $v$, in $I$.
Then, for any distinct $y,z \in B$, the bipartite graph $G[X-v-y-z]$
is matchable.
Furthermore, if $|B| \geq 3$, then the bipartite graph $G[X-v]$
is connected and is free of nontrivial $1$-cuts.
\end{lemma}
\begin{proof}
Note that the end $u$ of edge $e$ lies in the unique nontrivial
component of $G-e-B$.

\smallskip
Let $y$ and $z$ denote distinct vertices in $B$.
Observe that $G[X-v-y-z]$ has equicardinal color classes:
$B-y-z$ and $I-v$.
Since $G$ is bicritical, $G-y-z$ has a perfect matching, say $M$.
Since the neighbourhood of $I-v$ is a subset of $B$,
the restriction of $M$ to $G[X-v]$ is in fact a perfect matching
of $G[X-v-y-z]$.

\smallskip
Now assume that $|B| \geq 3$, and let $H:=G[X-v]$.
Note that $H$ is bipartite with color classes $I-v$ and $B$;
furthermore, $|B|=|I-v|+2$.
Observe that, in a bicritical graph, the neighbourhood of any
(nonempty) stable set $S$ has cardinality at least $|S| + 2$;
in particular, this holds for each nonempty subset of $I-v$.
These facts imply that $H$ is connected. It remains
to show that $H$ is free of nontrivial $1$-cuts.

\smallskip
Suppose, to the contrary, that $\{ab\}$ is a nontrivial $1$-cut of $H$,
where $ab$ is an edge with $a \in I-v$ and $b \in B$.
We let $H_1$ and $H_2$ denote the two (nontrivial) components of $H-ab$
such that $a \in V(H_1)$ and $b \in V(H_2)$.
We let $k_1 := |V(H_1) \cap (I-v)|$ and $k_2 := |V(H_2) \cap (I-v)|$.
It follows from the observations in the preceding paragraph
that $|V(H_1) \cap B| \geq k_1+1$ and $|V(H_2) \cap B| \geq k_2 + 2$.
Consequently,
$|B| \geq k_1 + k_2 + 3 = |I-v|+3$.
This is a contradiction.
Thus $H$ is indeed free of nontrivial $1$-cuts.
This completes the proof of Lemma~\ref{lem:bipartite-shore-properties}.
\end{proof}

In order to  prove Theorem~\ref{thm:efeccb-removable-edges}, we need
  the  Three Case  Lemma, a  result  of Carvalho,  Lucchesi and  Murty
  \cite{clm02} that  plays an important role  in a few of  their works
  \cite{clm02a, clm06, clm12}.

\smallskip
As we are  dealing with
  cubic bricks, one  of the cases of  the lemma does not  apply.
So, in fact, the version of the lemma stated below has only two cases.
\begin{lemma}
    \label{lem:three-case}
    Let $G$ be a  bicritical graph and let $e$ be  a removable edge of
    $G$ such that every barrier of  $G-e$ is special.  If $G-e$ is not
    bicritical then one of the following holds:
    \begin{enumerate}[(i)]
    \item
      The graph $G-e$ has only one maximal nontrivial barrier, say $B$.
      The graph $(G-e)/X$ is bicritical,
      and $e$ has at least one end in $I$.
    \item
      The graph $G-e$ has two maximal nontrivial barriers, say $B$ and $B^*$.
      The set $B' := B^* - I$ is the unique maximal nontrivial barrier of $(G-e)/X$;
      furthermore, $B'$ is a barrier of $G-e$ as well, and $I'=I^*-B$.
      The graph $((G-e)/X)/X'$ is bicritical,
      and the edge $e$ has one end in $I$ and the other end in $I'$. \qed
    \end{enumerate}
\end{lemma}

Now          we         are          ready         to          prove
  the aforementioned strengthening of Theorem~\ref{thm:efeccb-removable-edges}.
\begin{theorem}
    \label{thm:efeccb-qbinv-edge}
    Let $G$ be an \efeccb, and let $e:=uv$ denote a removable edge that
    is not  {\binv}.  Then,  $e$ is  \qbinv.  Moreover,  the following
    properties hold:
    \begin{enumerate}[(i)]
    \item
      \label{itm:two-barriers}
      The graph $G-e$ has two nontrivial special barriers, $B$ and $B'$, such
      that at least one of them is a maximal barrier,
      the sets $X$ and $X'$ are disjoint, $v \in I$ and $u \in I'$,
      and the graph
      $H:=(G-e) / (X \rightarrow x) / (X' \rightarrow x')$ is bicritical.
      (See Figure~\ref{fig:structure-of-G-wrt-qbinv-edge}.)
    \item
      \label{itm:x-xp}
      In $H$, each contraction vertex, $x$  and $x'$,
      has degree exactly four; whereas every other vertex is cubic.
      Furthermore, $b(H)=b(G-e)$, and $\{x,x'\}$ is the unique 2-separation of $H$.
    \item
      \label{itm:L-Lp}
      The graph $H-x-x'$ has precisely  two (even) components,
      $L$ and $L'$.
      Each of $x$~and~$x'$ has exactly two distinct neighbours in $L$,
      and likewise, in $L'$. Consequently, $H$ is a simple graph,
      and each of $B$ and $B'$ is a maximal barrier of $G-e$.
      Furthermore, in~$G$, the two edges joining $L$~and~$X$ are nonadjacent;
      an analogous statement holds for $L$~and~$X'$, for $L'$~and~$X$,
      and for $L'$~and~$X'$.
    \item
      \label{itm:qbinv}
      The graph $H$ has two bricks; that is, $b(H)=2$.
      Furthermore, the cubic graphs  $J:=H-L'+xx'$ and
      $J':=H-L+xx'$ are isomorphic to the underlying simple
      graphs of the two bricks of $H$.
    \item
      \label{itm:nonbipartite-subgraphs}
      Each of the four graphs $G[V(L) \cup X], G[V(L') \cup X'], G[V(L') \cup X]$
      and $G[V(L) \cup X']$ is nonbipartite.
    \item
      \label{itm:3-cut-x-xp}
      Every nontrivial  3-cut of  $J$ contains the edge $xx'$ (if such cuts exist).
      An analogous statement holds for $J'$.
    \item
      \label{itm:4-cut-matching}
      If $L$ is not  isomorphic to $K_2$, then
      $\partial(V(L))$ is a matching in $G$, and $L$
      is 2-connected.  An analogous statement holds for $L'$.
    \end{enumerate}
\end{theorem}
\begin{proof}
We will prove statements (i) to (vii) in order.

\medskip
\noindent (\ref{itm:two-barriers}) Since $G-e$ has vertices of
degree two, it is not bicritical.
Also, by Lemma~\ref{lem:every-barrier-special}, every barrier of $G-e$ is special.
We may thus invoke the Three Case Lemma~(\ref{lem:three-case}).

\smallskip
First suppose that $G-e$ has only one maximal nontrivial barrier $B$,
whence $H':=(G-e)/(X \rightarrow x)$ is bicritical, and the edge $e$
has at least one end in $I$.
Since $H'$ is bicritical, each of its vertices has degree at least three,
whence the edge $e$ in fact has both ends in $I$.
Note that $\partial(X)$ is a nontrivial tight cut of $G-e$,
and the graph $H'$ is one of the \mbox{$\partial(X)$-contractions} of $G-e$.
The other \mbox{$\partial(X)$-contraction} of $G-e$ is bipartite, whence it does not
yield any bricks.
Thus $b(G-e)=b(H')$.
Observe that each vertex of $H'$, except perhaps the contraction vertex $x$,
is cubic. By Corollary~\ref{cor:bicritical-not-brick}, $H'$ is in fact a brick.
Thus $b(G-e)=b(H')=1$. In other words, the edge $e$ is \binv,
contrary to our hypothesis.

\smallskip
It follows from the Three Case Lemma (\ref{lem:three-case})
that $G-e$ has two maximal nontrivial special barriers, say $B$ and $B^*$.
By Theorem~\ref{thm:canonical-partition}, $B$ and $B^*$ are disjoint,
whence $I$ and $I^*$ are disjoint.
Furthermore, $B' = B^* - I$ is also a nontrivial special barrier of $G-e$,
and $I'=I^*-B$.
Consequently, $X$ and $X'$ are disjoint.
Also, the graph
\mbox{$H:=(G-e)/(X \rightarrow x)/(X' \rightarrow x')$}
is bicritical, and the edge $e$ has one end in $I$ and the other end in $I'$.
We may adjust notation so that $v \in I$ and $u \in I'$.

\medskip
\noindent (\ref{itm:x-xp})
Observe that $\overline{X}$ is precisely the vertex set of the unique
nontrivial component of $G-e-B$.
We let $C:=\partial(X)$.
By Lemma~\ref{lem:every-barrier-special},
$|C|=5$. Likewise, $|C'|=5$, where $C':=\partial(X')$.
Since $e$ has one end in $I$ and the other end in $I'$,
and since $X$ and $X'$ are disjoint,
$e \in C \cap C'$.
Observe that, in $G-e$, the cuts $C-e$ and $C'-e$ are laminar nontrivial tight cuts,
and $H$ is obtained by shrinking their disjoint shores, $X$ and $X'$,
to single vertices $x$ and $x'$, respectively.
Since $|C-e|=|C'-e|=4$, each of the contraction vertices
$x$ and $x'$ has degree exactly four in~$H$,
Clearly, every other vertex of $H$ is cubic.
Furthermore, each of $(G-e)/\overline{X}$
and $(G-e)/\overline{X'}$ is bipartite; whence they do not yield any bricks.
Thus, $b(H)=b(G-e)$.

\smallskip
Since $e$ is not \binv, $H$ is not a brick. However, $H$ is bicritical.
By invoking Theorem~\ref{thm:elp-brick-characterization} and
Corollary~\ref{cor:bicritical-not-brick}, we conclude
that $\{x,x'\}$ is the unique $2$-separation of $H$.

\medskip
\noindent(\ref{itm:L-Lp})
Since $H$ is bicritical, each component of $H-x-x'$ is even.
Let $L$ and $L'$ denote two distinct components of $H-x-x'$.
As $H$ is $2$-connected, $x$ has at least one neighbour in $L$,
say $w$. If $w$ is the only neighbour of $x$ in $L$
then $\{w,x'\}$ is a barrier in~$H$, contradicting
the fact that $H$ is bicritical.
Thus, $x$ has at least two neighbours in $L$.
Likewise, $x$ has at least two neighbours in $L'$.
Since $x$ has degree exactly four in $H$, we infer that $L$ and $L'$
are the only components of $H-x-x'$, and that $x$
has exactly two neighbours in $L$, and likewise,
$x$ has exactly two neighbours in $L'$.
A similar conclusion holds for the vertex $x'$.
These facts imply that $H$ is in fact a simple graph.

\smallskip
As noted earlier, in the proof of (\ref{itm:two-barriers}),
$B$ is a maximal nontrivial barrier of $G-e$, and $B'$
is a subset of the maximal nontrivial barrier $B^*$.
Furthermore, $B^* \subseteq B' \cup I$ and $I^* \subseteq I' \cup B$,
whence $X^* \subseteq X \cup X'$.
Suppose that $B'$ is a proper subset of $B^*$, whence $|B^*| \geq 3$.
By Lemma~\ref{lem:bipartite-shore-properties}, the subgraph~$G[X^*-u]$
is connected; its vertex set
meets each of the sets $X'-u$ and $X$; however, $G$ has no edges
joining these sets. This is absurd.
We thus conclude that $B'=B^*$. In other words,
$B'$ is in fact a maximal barrier of $G-e$.

\smallskip
Now, consider the two edges joining $X$ and $L$, say $d$ and $f$.
We have already established that $H$ is simple;
in particular, $d$ and $f$ do not share
a common end in $L$. Observe that $\partial(V(L))$ is a $4$-cut in $G$.
Also, the shore $\overline{V(L)}$ is clearly not isomorphic to $K_2$.
By Lemma~\ref{lem:4-cut}, $d$~and $f$~do not share a common end in $X$.
In other words, $d$ and $f$ are nonadjacent.

\medskip
\noindent(\ref{itm:qbinv})
We consider a $2$-separation cut of $H$
associated with its only $2$-separation $\{x,x'\}$,
say $D:=\partial_H(V(L) \cup \{x\})$. Observe
that the underlying simple graph of one of the two \mbox{$D$-contractions} of $H$
is isomorphic to $J:=H-L'+xx'$, and that the underlying simple
graph of the other $D$-contraction is isomorphic to $J':=H-L+xx'$.
Thus, $b(G-e)=b(H)=b(J)+b(J')$.
Since $H$ is bicritical, Proposition~\ref{prp:bicritical-tight-cut}
implies that each of $J$ and $J'$ is bicritical;
observe that each of them is cubic.
Corollary~\ref{cor:bicritical-not-brick} implies that
each of $J$ and $J'$ is in fact a brick.
Thus $b(G-e)=b(H)=2$.
In particular, $e$ is a \qbinv\ edge of $G$.

\medskip
\noindent(\ref{itm:nonbipartite-subgraphs})
As noted above, $D:=\partial_H(V(L) \cup \{x\})$
is a nontrivial tight cut of $H$. Thus, $D$ is a nontrivial tight cut
of $G-e$ as well, and its shores are $V(L) \cup X$ and $V(L') \cup X'$.
Each $D$-contraction of $G-e$ is a \mcg\ that
has exactly one brick;
in particular, each $D$-contraction is nonbipartite.
Thus, each shore of $D$ induces a nonbipartite subgraph.
An analogous argument shows that each of $G[V(L') \cup X]$
and $G[V(L) \cup X']$ is also nonbipartite.

\medskip
\noindent(\ref{itm:3-cut-x-xp})
Let $F$ denote any cut of $J$ such that $xx' \notin F$.
Then $F$ has a shore that is disjoint with $\{x,x'\}$,
whence $F$ is a cut of $G$ as well.
Since $G$ is \efec, $F$ is a $3$-cut if and only if $F$ is trivial.
Consequently, every nontrivial $3$-cut of $J$ contains the edge $xx'$.

  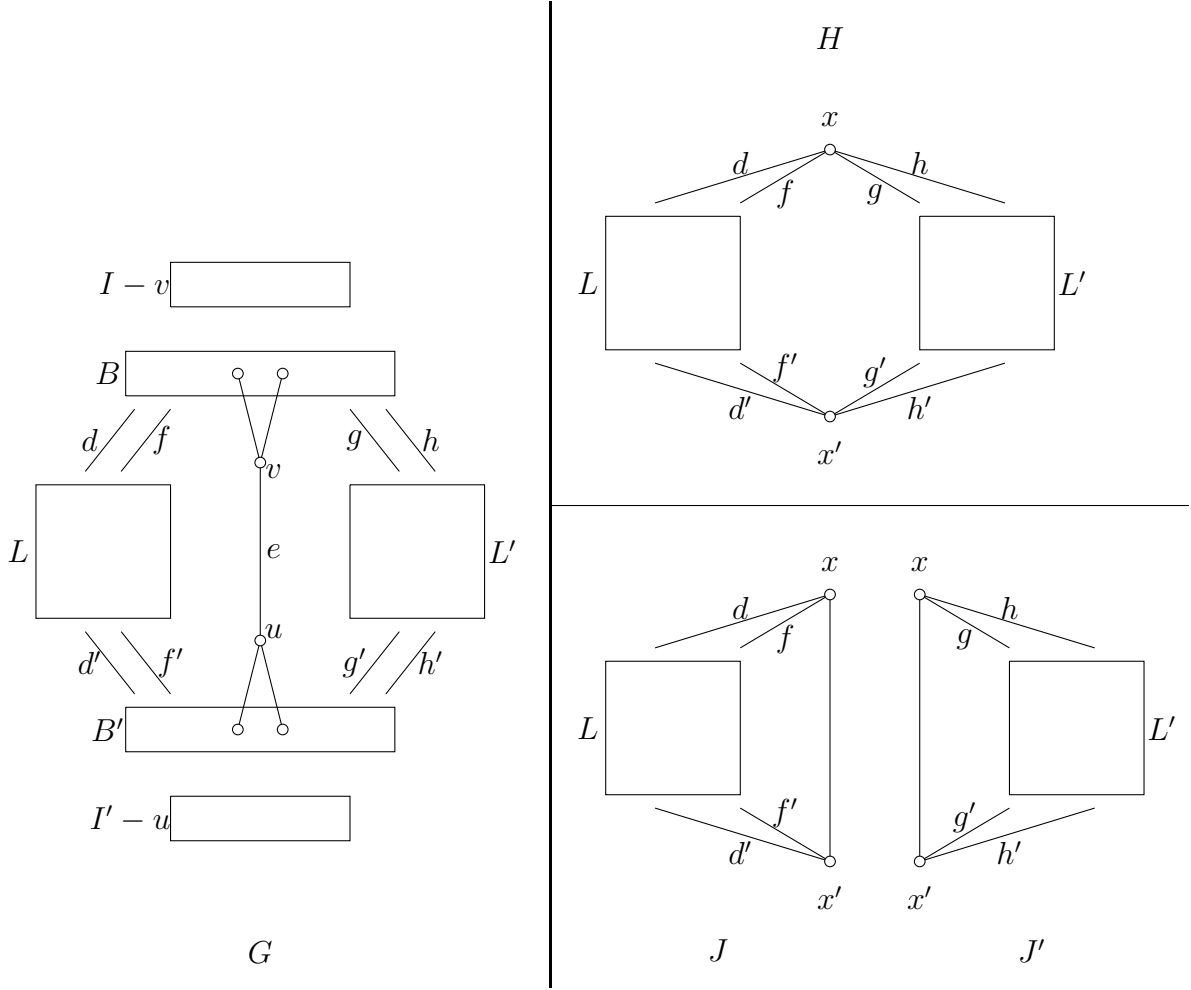
\begin{figure}[!htb]
    \centering
    \begin{tikzpicture}[scale=0.59]

       \draw (0,2) -- (0,-2);
      \draw (0.3,0)node[nodelabel]{$e$};

      \draw (0,2) -- (-0.5,4)node{};
      \draw (0,2) -- (0.5,4)node{};

      \draw (0,-2) -- (-0.5,-4)node{};
      \draw (0,-2) -- (0.5,-4)node{};

      \draw (-3,3.5) -- (3,3.5) -- (3,4.5) -- (-3,4.5) -- (-3,3.5);
      \draw (-3.4,4)node[nodelabel]{$B$};

      \draw (-2,5.5) -- (2,5.5) -- (2,6.5) -- (-2,6.5) -- (-2,5.5);
      \draw (-2.8,6)node[nodelabel]{$I-v$};

      \draw (-3,-3.5) -- (3,-3.5) -- (3,-4.5) -- (-3,-4.5) -- (-3,-3.5);
      \draw (-3.4,-4)node[nodelabel]{$B'$};

      \draw (-2,-5.5) -- (2,-5.5) -- (2,-6.5) -- (-2,-6.5) -- (-2,-5.5);
      \draw (-2.9,-6)node[nodelabel]{$I'-u$};

      \draw (-2,1.5) -- (-5,1.5) -- (-5,-1.5) -- (-2,-1.5) -- (-2,1.5);
      \draw (-5.4,0)node[nodelabel]{$L$};

      \draw (2,1.5) -- (5,1.5) -- (5,-1.5) -- (2,-1.5) -- (2,1.5);
      \draw (5.4,0)node[nodelabel]{$L'$};

      \draw (-2.8,3.2) -- (-3.9,1.8);
      \draw (-3.8,2.5)node[nodelabel]{$d$};
      \draw (-2,3.2) -- (-3.1,1.8);
      \draw (-2.2,2.5)node[nodelabel]{$f$};

      \draw (2.8,3.2) -- (3.9,1.8);
      \draw (3.8,2.5)node[nodelabel]{$h$};
      \draw (2,3.2) -- (3.1,1.8);
      \draw (2.1,2.5)node[nodelabel]{$g$};

      \draw (-2.8,-3.2) -- (-3.9,-1.8);
      \draw (-3.8,-2.5)node[nodelabel]{$d'$};
      \draw (-2,-3.2) -- (-3.1,-1.8);
      \draw (-2,-2.5)node[nodelabel]{$f'$};

      \draw (2.8,-3.2) -- (3.9,-1.8);
      \draw (3.8,-2.5)node[nodelabel]{$h'$};
      \draw (2,-3.2) -- (3.1,-1.8);
      \draw (2.1,-2.5)node[nodelabel]{$g'$};

      \draw (0,2)node{};
      \draw (0.3,1.8)node[nodelabel]{$v$};
      \draw (0,-2)node{};
      \draw (0.3,-1.8)node[nodelabel]{$u$};

      \draw (0,-9)node[nodelabel]{$G$};

    \end{tikzpicture}
    \vline
    \begin{tikzpicture}[scale=0.59]

      \draw (-2,1.5) -- (-5,1.5) -- (-5,-1.5) -- (-2,-1.5) -- (-2,1.5);
      \draw (-5.4,0)node[nodelabel]{$L$};

      \draw (2,1.5) -- (5,1.5) -- (5,-1.5) -- (2,-1.5) -- (2,1.5);
      \draw (5.4,0)node[nodelabel]{$L'$};

      \draw (0,3) -- (-3.9,1.8);
      \draw (-2,2.7)node[nodelabel]{$d$};
      \draw (0,3) -- (-2,1.8);
      \draw (-1,2)node[nodelabel]{$f$};

      \draw (0,3) -- (3.9,1.8);
      \draw (2,2.7)node[nodelabel]{$h$};
      \draw (0,3) -- (2,1.8);
      \draw (1,2)node[nodelabel]{$g$};

      \draw (0,-3) -- (-3.9,-1.8);
      \draw (-2,-2.8)node[nodelabel]{$d'$};
      \draw (0,-3) -- (-2,-1.8);
      \draw (-1,-1.9)node[nodelabel]{$f'$};

      \draw (0,-3) -- (3.9,-1.8);
      \draw (2,-2.8)node[nodelabel]{$h'$};
      \draw (0,-3) -- (2,-1.8);
      \draw (1,-2)node[nodelabel]{$g'$};

      \draw (0,3)node{}node[above,nodelabel]{$x$};
      \draw (0,-3)node{}node[below,nodelabel]{$x'$};

      \draw (-6.2,-5) -- (8,-5);

      \draw (0,5.5)node[nodelabel]{$H$};

      \draw (-2,1.5-10) -- (-5,1.5-10) -- (-5,-1.5-10) -- (-2,-1.5-10) -- (-2,1.5-10);
      \draw (-5.4,0-10)node[nodelabel]{$L$};
      \draw (0,3-10) -- (-3.9,1.8-10);
      \draw (-2,2.7-10)node[nodelabel]{$d$};
      \draw (0,3-10) -- (-2,1.8-10);
      \draw (-1,2-10)node[nodelabel]{$f$};
      \draw (0,-3-10) -- (-3.9,-1.8-10);
      \draw (-2,-2.8-10)node[nodelabel]{$d'$};
      \draw (0,-3-10) -- (-2,-1.8-10);
      \draw (-1,-1.9-10)node[nodelabel]{$f'$};
      \draw (0,3-10) -- (0,-3-10);
      \draw (0,3-10)node{}node[above,nodelabel]{$x$};
      \draw (0,-3-10)node{}node[below,nodelabel]{$x'$};
      \draw (-2.5,-15)node[nodelabel]{$J$};

      \draw (2+2,1.5-10) -- (5+2,1.5-10) -- (5+2,-1.5-10) -- (2+2,-1.5-10) -- (2+2,1.5-10);
      \draw (5.4+2,0-10)node[nodelabel]{$L'$};
      \draw (0+2,3-10) -- (3.9+2,1.8-10);
      \draw (2+2,2.7-10)node[nodelabel]{$h$};
      \draw (0+2,3-10) -- (2+2,1.8-10);
      \draw (1+2,2-10)node[nodelabel]{$g$};

      \draw (0+2,-3-10) -- (3.9+2,-1.8-10);
      \draw (2+2,-2.8-10)node[nodelabel]{$h'$};
      \draw (0+2,-3-10) -- (2+2,-1.8-10);
      \draw (1+2,-2-10)node[nodelabel]{$g'$};

      \draw (0+2,3-10) -- (0+2,-3-10);

      \draw (0+2,3-10)node{}node[above,nodelabel]{$x$};
      \draw (0+2,-3-10)node{}node[below,nodelabel]{$x'$};

      \draw (4.5,-15)node[nodelabel]{$J'$};

    \end{tikzpicture}
    \caption{(left) An \efeccb\ $G$ with a \qbinv\ edge $e:=vu$ ; (top right) bicritical graph $H$ obtained from $G-e$ ;
      (bottom right) cubic bricks $J$ and $J'$ obtained from $H$}
    \label{fig:structure-of-G-wrt-qbinv-edge}
  \end{figure}
\medskip
\noindent(\ref{itm:4-cut-matching})
As noted earlier, $\partial(V(L))$
is a $4$-cut of $G$. It follows from statement (\ref{itm:L-Lp})
that, for any two edges in $\partial(V(L))$, their ends in $\overline{V(L)}$
are distinct.
If $L$ is not isomorphic to $K_2$,
then Lemma~\ref{lem:4-cut} implies that, for any two edges in $\partial(V(L))$,
their ends in $V(L)$ are also distinct, whence $\partial(V(L))$ is indeed a matching in $G$.

\smallskip
Now, suppose that $L$ is not isomorphic to $K_2$.
Assume, to the contrary, that $L$ is not $2$-connected,
and let $w$ denote a cut-vertex of $L$.
Let $L_1$ and $L_2$ denote distinct components of $L-w$.
For $i \in \{1,2\}$, let $F_i:=\partial(V(L_i))$.
Then $F_1 \cup F_2 \subseteq \partial(w) \cup \partial(V(L))$.
Since $F_1$~and~$F_2$ are disjoint,
$|F_1| + |F_2| \leq |\partial(w)| + |\partial(V(L))| = 7$.
Thus, at least one of $F_1$~and~$F_2$ is a $3$-cut in $G$.
Adjust notation so that $F_1$ is a $3$-cut.
Clearly, the edge $xx'$ does not lie in $F_1$.
Thus $F_1$ is a $3$-cut in $J$.
It follows from statement (\ref{itm:3-cut-x-xp})
that $F_1$ is a trivial cut of $J$.
However, since $F_1 \subseteq \partial(w) \cup \partial(V(L))$,
we infer that two edges of $\partial(V(L))$ share a common end in $V(L)$,
contradicting what we have established in the preceding paragraph.
We thus conclude that if $L$ is not isomorphic to $K_2$
then $L$ is indeed $2$-connected.

\medskip
This completes the proof of Theorem~\ref{thm:efeccb-qbinv-edge}.
\end{proof}

\smallskip
We conclude this section with an easy lemma --- that will be very useful
in proving the Main Theorem (\ref{thm:efeccb-two-qbinv-at-vertex}).
It appears implicitly and crucially in the proof of the main result
of Carvalho, Lucchesi and Murty \cite{clm06};
we include their proof for the sake of completeness.
\begin{lemma}
\label{lem:intersection-of-barriers}
Let $G$ be a bicritical graph, and let $e_1$ and $e_2$ be two adjacent edges.
If $B_1$ is a barrier of $G-e_1$, and if $B_2$ is a barrier of $G-e_2$, then $|B_1 \cap B_2| \leq 1$.
\end{lemma}
\begin{proof}
Suppose to the contrary that $B_1 \cap B_2$ contains two distinct vertices $s$ and $t$.
Since $G$ is bicritical, let $M$ denote a perfect matching of $G-s-t$. Since $s,t \in B_1$, the graph $G-e_1-s-t$ has no perfect matching,
whence $e_1 \in M$. Likewise, we infer that $e_2 \in M$. This is absurd since $e_1$ and $e_2$ are adjacent edges.
Thus, $|B_1 \cap B_2| \leq 1$.
\end{proof}

\subsection{Matchable subgraphs}

In the proof of the Main Theorem~(\ref{thm:efeccb-two-qbinv-at-vertex}),
we will often need to construct perfect matchings of certain subgraphs
of the brick under consideration. In this section,
we prove some technical lemmas that will help us in constructing
the desired perfect matchings.

\smallskip
Before that, we state two results concerning
bipartite matchable graphs that may be easily derived from
Hall's Theorem.
These will be invoked only in the proof of
Lemma~\ref{lem:near-bipartite-J-admissibility-of-xx'}.
\begin{proposition}
\label{prp:bip-mc-characterization}
    Let  $H[A,B]$ denote  a  bipartite matchable  graph.  Then $H$  is
    \mc\ if and only if $H-a-b$ is matchable for each pair
    of vertices $a \in A$ and $b \in B$. \qed
\end{proposition}
\begin{proposition}
\label{prp:bip-matchable-not-mc}
    Let $H[A,B]$ denote  a bipartite matchable graph. If  an edge $ab$
    is inadmissible then there exist partitions $[A_1,A_2]$ of $A$ and
    $[B_1,B_2]$ of  $B$ such that $a  \in A_2$, $b \in  B_1$, $|A_1| =
    |B_1|$ and that there are no edges from $A_1$ to $B_2$. \qed
\end{proposition}

Now let $e:=uv$ denote a \qbinv\ edge of an \efeccb~$G$.
We invoke Theorem 3.7, and adopt all of the notation therein
(pertaining to the structure of $G$ with respect to~$e$),
as well as the following notation.
(See Figure~\ref{fig:structure-of-G-wrt-qbinv-edge}.)
\begin{notation}
\label{Not:8-edges}
Let $\partial(V(L)):=\{d,f,d',f'\}$, where each of $d$ and $f$
has one end in $X$, and each of $d'$ and $f'$ has one end in $X'$.
Likewise, let $\partial(V(L')):=\{g,h,g',h'\}$, where each of $g$ and $h$
has one end in $X$, and each of $g'$ and $h'$ has one end in $X'$.
\end{notation}

The next three lemmas will help us in constructing perfect
matchings of certain subgraphs with the additional property
that a specified edge is included.
\begin{lemma}
\label{lem:L-p-q-matchable-e'-admissible-in-G-p-q}
Let $p$ and $q$ denote distinct vertices of $L$,
and let $e' \in \partial(v)-e$.
If $L-p-q$ is matchable then $e'$ is admissible in $G-p-q$.
\end{lemma}
\begin{proof}
Let $y$ denote the end of $e'$ in $B$. So $e'=vy$.
By Theorem~\ref{thm:efeccb-qbinv-edge}(\ref{itm:L-Lp}),
$g$ and $h$ are nonadjacent edges,
whence at least one of them, say $g$, is not in $\partial(y)$.
We let $z$ denote the end of $g$ in $B$. (Thus $y \neq z$.)
Now, let $M'$ be a perfect matching of the brick $J'$
such that $g \in M'$.
Clearly, $M'$ contains exactly one of $g'$ and $h'$.
Adjust notation so that $g' \in M'$,
and let $z'$ denote the end of $g'$ in $B'$.
The vertex $u$ has two distinct neighbours in $B'$,
and we let $y' \in B'$ denote a neighbour of $u$ that is distinct from $z'$.

\smallskip
We now invoke Lemma~\ref{lem:bipartite-shore-properties} twice.
The graph $G[X-v-y-z]$ has a perfect matching, say $N$.
Likewise, $G[X'-u-y'-z']$ has a perfect matching, say~$N'$.
By hypothesis, $L-p-q$ has a perfect matching, say $M$.
Observe that $M \cup M' \cup N \cup N' \cup \{vy,uy'\}$
is a perfect matching of $G-p-q$ that contains the edge $e'$, as desired.
\end{proof}
\begin{lemma}
\label{lem:extend-M-union-Mp-union-e}
Let $p$ and $q$ denote distinct vertices of $L$, let $M$ be a perfect
matching of $J-p-q$, and let $M'$ be a perfect matching of $J'$.
Suppose that $M \cup M'$ does not contain the edge $xx'$,
and that $M \cup M'$ is a matching in $G$.
Then $M \cup M' \cup \{e\}$ may be extended to a perfect
matching of $G-p-q$.
\end{lemma}
\begin{proof}
Note that $xx' \notin M \cup M'$.
Thus $M$ contains precisely one edge incident
with $x$, and it contains precisely one edge incident with $x'$.
An analogous statement holds for $M'$.
By hypothesis, $M \cup M'$ is a matching in $G$.
We let $y$ and $z$ denote the two distinct vertices of $X$ that are incident with edges
in $M \cup M'$.
Likewise, we let $y'$ and $z'$ denote the two distinct vertices of $X'$ that are incident
with edges in $M \cup M'$.

\smallskip
We now invoke Lemma~\ref{lem:bipartite-shore-properties} twice.
The graph $G[X-v-y-z]$ has a perfect matching, say $N$.
Likewise, $G[X'-u-y'-z']$ has a perfect matching, say $N'$.
Observe that $M \cup M' \cup \{e\} \cup N \cup N'$ is a perfect
matching of $G-p-q$. This completes the proof.
\end{proof}
\begin{lemma}
\label{lem:near-bipartite-J-admissibility-of-xx'}
Suppose that $\{d,d'\}$ is a removable doubleton of $J$.
Let $d:=xy$, and let $[T,T']$ denote the bipartition of $J-d-d'$
such that $x,y \in T$. Then, for each $p \in T'-x'$,
the edge $xx'$ is admissible in the bipartite graph $J-d-d'-p-y$.
\end{lemma}
\begin{proof}
Since  $J-d-d'$ is bipartite and  \mc,
    Proposition~\ref{prp:bip-mc-characterization}     implies     that
    $J-d-d'-p-y$ is matchable, and it  is bipartite with color classes
    $T-y$  and   $T'-p$.

\smallskip
Suppose  that  $xx'$   is  inadmissible  in
    $J-d-d'-p-y$.     By   Proposition~\ref{prp:bip-matchable-not-mc},
    there exist  partitions $[T_1,T_2]$ of $T-y$  and $[T'_1,T'_2]$ of
    $T'-p$ such that  $x \in T_1$, $x' \in  T'_2$, $|T_1|=|T'_1|$, and
    that there are no edges between $T'_1$ and $T_2$. Since $d=xy$ and
    $x  \in T_1$,  no end  of  $d$ is  in  $T_2$. Thus,  in $J$,  each
    neighbour of  $T_2$ lies in  $T'_2 \cup  \{p\}$, whence $J$  has a
    nontrivial  barrier, contradiction.   We  conclude  that $xx'$  is
    indeed admissible in $J-d-d'-p-y$.
\end{proof}

\subsection{Bricks of order $10$}

In this section, we describe the {\efeccb}s, up to order~$10$,
that have a \qbinv\ edge. It is an easy application
of Theorem~\ref{thm:efeccb-qbinv-edge}.

\smallskip
Let $G$ be an \efeccb\ that has a \qbinv\ edge $e:=\{u,v\}$.
We invoke Theorem~\ref{thm:efeccb-qbinv-edge}, and adopt all of the notation
therein regarding the structure of $G$
with respect to edge $e$. See Figure~\ref{fig:structure-of-G-wrt-qbinv-edge}.
Note that each of the sets $B, B', V(L)$ and $V(L')$ has at least two vertices.
Since these sets are pairwise-disjoint, and since none of them meets $\{u,v\}$,
the graph $G$ has at least $10$ vertices.

\begin{figure}[!htb]
\centering
\subfigure [A subgraph of each]
{
\begin{tikzpicture}[scale=0.55]

	\draw (1,-4) to [out=0,in=270] (4,0);
	\draw (-1,-4) to [out=30,in=270] (2,0)node{};

	\draw (-4,0) -- (-2,0);
	\draw (-1,4) to [out=180,in=90] (-4,0)node{};

	\draw (1,4) to [out=210,in=90] (-2,0)node{};

	\draw (4,0) -- (2,0);
	\draw (1,4) to [out=0,in=90] (4,0)node{};

	\draw (-1,4) to [out=330,in=90] (2,0)node{};

	\draw (0,2) -- (0,-2);
	\draw (0.4,0)node[nodelabel]{$e$};

	\draw (0,2) -- (-1,4)node{};
	\draw (0,2) -- (1,4)node{};

	\draw (0,-2) -- (-1,-4)node{};

	\draw (0,-2) -- (1,-4)node{};

	\draw (0,2)node{};
	\draw (0.4,1.8)node[nodelabel]{$v$};
	\draw (0,-2)node{};
	\draw (0.4,-1.8)node[nodelabel]{$u$};

      \end{tikzpicture}
\label{fig:subgraph-of-order-10-brick}
}
\hspace*{0.2in}
\subfigure [The Petersen graph]
{
\begin{tikzpicture}[scale=0.55]

	\draw (-2,0) to [out=270,in=150] (-1,-4);

	\draw (-4,0) to [out=270,in=135] (1,-4);

	\draw (1,-4) to [out=0,in=270] (4,0);

	\draw (-1,-4) to [out=30,in=270] (2,0)node{};

	\draw (-4,0) -- (-2,0);
	\draw (-1,4) to [out=180,in=90] (-4,0)node{};

	\draw (1,4) to [out=210,in=90] (-2,0)node{};

	\draw (4,0) -- (2,0);
	\draw (1,4) to [out=0,in=90] (4,0)node{};

	\draw (-1,4) to [out=330,in=90] (2,0)node{};

	\draw (0,2) -- (0,-2);
	\draw (0.4,0)node[nodelabel]{$e$};

	\draw (0,2) -- (-1,4)node{};
	\draw (0,2) -- (1,4)node{};

	\draw (0,-2) -- (-1,-4)node{};

	\draw (0,-2) -- (1,-4)node{};

	\draw (0,2)node{};
	\draw (0.4,1.8)node[nodelabel]{$v$};
	\draw (0,-2)node{};
	\draw (0.4,-1.8)node[nodelabel]{$u$};

      \end{tikzpicture}
\label{fig:Petersen}
}
\hspace*{0.2in}
\subfigure [Petersen's Mate]
{
\begin{tikzpicture}[scale=0.55]

       \draw (-1,-4) to [out=180,in=270] (-4,0);
       \draw (1,-4) to [out=150,in=270] (-2,0);

	\draw (1,-4) to [out=0,in=270] (4,0);

	\draw (-1,-4) to [out=30,in=270] (2,0)node{};

	\draw (-4,0) -- (-2,0);
	\draw (-1,4) to [out=180,in=90] (-4,0)node{};

	\draw (1,4) to [out=210,in=90] (-2,0)node{};

	\draw (4,0) -- (2,0);
	\draw (1,4) to [out=0,in=90] (4,0)node{};

	\draw (-1,4) to [out=330,in=90] (2,0)node{};

	\draw (0,2) -- (0,-2);
	\draw (0.4,0)node[nodelabel]{$e$};

	\draw (0,2) -- (-1,4)node{};
	\draw (0,2) -- (1,4)node{};

	\draw (0,-2) -- (-1,-4)node{};
	\draw (0,-2) -- (1,-4)node{};

	\draw (0,2)node{};
	\draw (0.4,1.8)node[nodelabel]{$v$};
	\draw (0,-2)node{};
	\draw (0.4,-1.8)node[nodelabel]{$u$};

      \end{tikzpicture}
\label{fig:Petersen-Mate}
}
\caption{Bricks of order $10$}
\label{fig:order-10}
\end{figure}
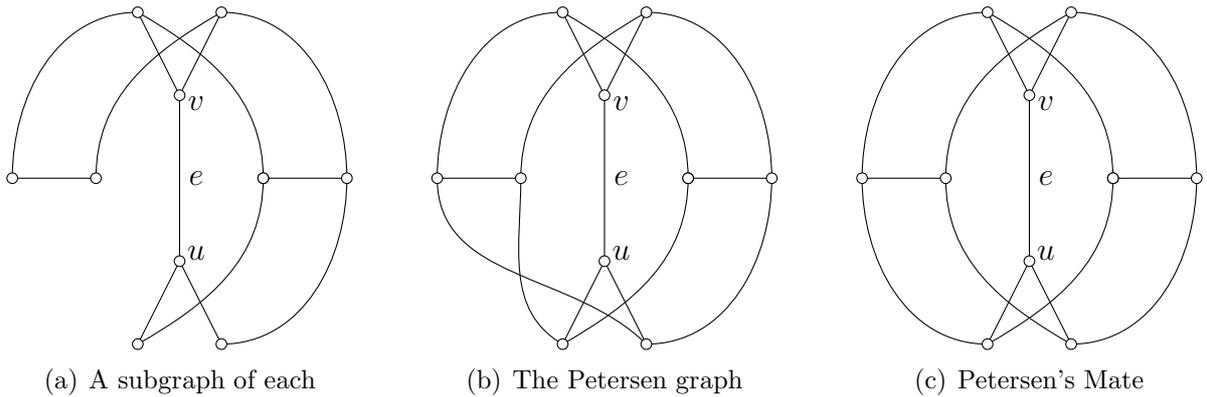

\smallskip
Now suppose that $|V(G)|=10$, whence each of the sets $B,B',V(L)$ and $V(L')$
is a doubleton. Furthermore, each of $L$ and $L'$ is isomorphic to $K_2$,
and each of $G[X]$ and $G[X']$ is a path of order $3$.
By Theorem~\ref{thm:efeccb-qbinv-edge}(\ref{itm:L-Lp}),
each of $\{d,f\}$, $\{d',f'\}$, $\{g,h\}$ and $\{g',h'\}$
is a pair of nonadjacent edges.
Consequently, the graph shown
in Figure~\ref{fig:subgraph-of-order-10-brick}
is a subgraph of~$G$.

\smallskip
Clearly, there are exactly two possibilities now. Either $G$
is the Petersen graph, that has girth $5$,
as shown in Figure~\ref{fig:Petersen}.
Otherwise, $G$ is the graph shown in Figure~\ref{fig:Petersen-Mate},
that has girth $4$;
we shall refer to this graph as {\it Petersen's Mate}.

Since the Petersen graph is edge-transitive, each of its edges
is \qbinv.
On the other hand, the reader may verify that Petersen's Mate
has exactly four edge-orbits: the edge~$e$
is the only \qbinv\ edge, each edge adjacent with $e$
participates in a removable doubleton, and members of
the remaining two edge-orbits are \binv.
We will find these facts, summarized below, useful in the proof of
the Main Theorem~(\ref{thm:efeccb-two-qbinv-at-vertex}).

\begin{proposition}
\label{prp:order-10}
There exist precisely two {\efeccb}s, of order at most $10$,
that have a \qbinv\ edge. These are the Petersen graph and
Petersen's Mate. Furthermore, each edge of the Petersen graph
is \qbinv; whereas, Petersen's Mate has a unique \qbinv\ edge. \qed
\end{proposition}

\section{Pfaffian graphs and conformal minors}
\label{sec:Pfaffian-and-conformal}
In this section, we discuss two related concepts: `Pfaffian orientations'
and `conformal minors'. We will use
Lemma~\ref{lem:rigid-bisubdivision-K23-implies-non-Pfaffian}
in the proof of the Main Theorem~(\ref{thm:efeccb-two-qbinv-at-vertex});
the reader may postpone reading the rest until
before Section~\ref{sec:main-theorem-consequences}.

\smallskip
Let $D$ be an orientation of an undirected graph~$G$.
For an even cycle~$C$ of~$G$,
we abuse notation and use $C$ to also refer to the
corresponding set of arcs in~$D$.
Note that regardless of the sense of traversal of~$C$,
the number of forward arcs and the number of reverse arcs
have the same parity.
We say that $C$ is {\it evenly-oriented} if the number
of forward arcs is even, and {\it oddly-oriented} otherwise.

\smallskip
Let $G$ be a matchable graph.
A cycle~$C$ of a graph~$G$ is {\it conformal} if $G-V(C)$ is matchable.
An orientation~$D$ (of~$G$) is {\it Pfaffian} if each conformal cycle
is oddly-oriented. Furthermore, we say that $G$ is a {\it Pfaffian} graph
if $G$ admits a Pfaffian orientation; otherwise
$G$ is {\it non-Pfaffian}.

\smallskip
The significance of Pfaffian orientations arises from the fact that
if a graph is Pfaffian then the number of its perfect matchings
may be computed in polynomial-time.
We may thus restrict ourselves to {\mcg}s.
Kasteleyn~\cite{kast63}
showed that all planar graphs are Pfaffian.
However, the
Pfaffian graph recognition problem remains unsolved for nonplanar graphs.

\begin{problem}
\label{prb:Pfaffian}
Characterize Pfaffian nonplanar {\mcg}s. (Is the problem of deciding whether
a given graph is Pfaffian in the complexity class $\mathcal{NP}$?
Is it in~$\mathcal{P}$?)
\end{problem}

The smallest non-Pfaffian graph is $K_{3,3}$.
We now describe a certificate that may be used to prove
that certain graphs are non-Pfaffian.
It should be noted that this certificate
does not exist in every non-Pfaffian graph. For instance,
the Petersen graph is non-Pfaffian; see~\cite{clm12}. However,
it does not contain the certificate we are about to describe.

\smallskip
To {\it bi-subdivide an edge} $e$ means to subdivide $e$ by inserting an
even number of vertices; or equivalently, to replace~$e$ by an odd path.
A {\it bi-subdivision} of a graph~$J$ is a graph~$H$ obtained from~$J$
by means of bi-subdividing a subset of its edges.
For a matchable graph~$G$,
we say that a subgraph~$H$ is a {\it rigid bi-subdivision of $K_{2,3}$}
if it satisfies the following properties:
\begin{itemize}
\item the subgraph~$H$ is a bi-subdivision of $K_{2,3}$ and
\item each cycle of~$H$ is conformal in~$G$.
\end{itemize}
\begin{lemma}
\label{lem:rigid-bisubdivision-K23-implies-non-Pfaffian}
Every matchable graph that has a rigid bi-subdivision of $K_{2,3}$ is non-Pfaffian.
\end{lemma}
\begin{proof}
Let $G$ be a matchable graph, and let $H$ denote a subgraph that is
a rigid bi-subdivision of $K_{2,3}$.
Note that $H$ has two cubic vertices, say~$u$ and $v$;
and it has three edge-disjoint \mbox{$uv$-paths}, each of even length,
say, $P_1, P_2$ and $P_3$. Let $D$ denote any orientation of~$G$.
We will argue that $D$ is not a Pfaffian orientation.

\smallskip
For $i,j \in \{1,2,3\}$,
such that $i < j$,
let $C_{i,j}$ denote the (even) cycle $P_i \cup P_j$, and
let $l_{i,j}$ denote the number of forward arcs in $C_{i,j}$ ---
traversing $P_i$ from $u$~to~$v$, and $P_j$ from $v$~to~$u$.
Since $H$ is rigid, each of $C_{1,2}, C_{2,3}$ and $C_{1,3}$
is a conformal cycle of~$G$.

\smallskip
Observe that each edge of the path~$P_2$ is a forward arc
for exactly one of $C_{1,2}$ and $C_{2,3}$.
Consequently, $l_{1,3} = l_{1,2} + l_{2,3} - |P_2|$.
Since $P_2$ is of even length, $l_{1,3} \equiv l_{1,2} + l_{2,3} {\rm ~(mod~2)}$.
It follows that at least one of the (conformal) cycles $C_{1,2}, C_{2,3}$ and $C_{1,3}$
is evenly-oriented. Thus $D$ is not a Pfaffian orientation, whence $G$
is non-Pfaffian.
\end{proof}

We will find the above lemma useful in the proof of the
Main Theorem (\ref{thm:efeccb-two-qbinv-at-vertex}).
Now we discuss another concept that is intrinsically related
to Pfaffian graphs.

\smallskip
Let $G$ be a \mcg.
A subgraph~$H$ of~$G$ is {\it conformal} if $G-V(H)$ is matchable.
Observe that $G$ is Pfaffian if and only if each conformal \mc\ subgraph of~$G$
is Pfaffian.

\smallskip
Now, let $J$ be a cubic \mcg.
We say that $J$ is a {\it conformal minor} of a \mcg~$G$ if 
the latter has a conformal subgraph~$H$ that is a bi-subdivision of~$J$.
For the sake of brevity,
we say that $G$ is {\it $J$-based} if the latter is a conformal minor of
the former; otherwise we say that $G$ is {\it $J$-free}.
The following fundamental result was proved by Little~\cite{litt75}.

\begin{theorem}
\label{thm:bipartite-Pfaffian-iff-K33free}
A bipartite \mcg\ is Pfaffian if and only if it is $K_{3,3}$-free.
\end{theorem}

It was shown by Little and Rendl \cite{lire91} that a \mcg\ is Pfaffian
if and only if each of its bricks and braces is Pfaffian. Consequently,
it suffices to solve Problem~\ref{prb:Pfaffian} for bricks and braces.
In particular, for the case of bipartite graphs, it suffices
to characterize the $K_{3,3}$-free braces; this feat was accomplished
by Robertson, Seymour and Thomas \cite{rst99}, and independently by
McCuaig \cite{mccu04}; their works yield a polynomial-time
algorithm to decide whether a given bipartite graph is Pfaffian or not.

\smallskip
Thus, in order to solve Problem~\ref{prb:Pfaffian},
one may restrict their attention to nonplanar and nonbipartite {\mcg}s,
and in particular to nonplanar bricks.

\smallskip
A result similar to Theorem~\ref{thm:bipartite-Pfaffian-iff-K33free}
was proved for near-bipartite graphs by Fischer and Little \cite{fili01}.
In particular, they proved that a near-bipartite \mcg\ is Pfaffian
if and only if it does not contain any of seven cubic graphs as a conformal minor;
three of these graphs are $K_{3,3}$, Cubeplex and Twinplex;
the remaining four are obtained from these three by replacing at most one
or two (specific) vertices by triangles.
Also, Miranda and Lucchesi \cite{milu08}
gave a polynomial-time algorithm to decide whether a given near-bipartite
graph is Pfaffian or not; their algorithm does not rely on the result
of Fischer and Little. Thus one may further restrict attention
to the non-near-bipartite bricks.

\begin{problem}
\label{prb:Pfaffian-restriction}
Characterize Pfaffian nonplanar non-near-bipartite bricks.
\end{problem}

Now we turn our attention to the following
fundamental result of Lov{\'a}sz \cite{lova83} ---
that has nothing to do with Pfaffian orientations.

\begin{theorem}
Every nonbipartite matching covered graph
is either $K_4$-based, or is $\overline{C_6}$-based, or both.
\end{theorem}

This gives rise to two natural problems.

\begin{problem}
\label{prb:K4-free}
Characterize $K_4$-free {\mcg}s.
\end{problem}
\begin{problem}
\label{prb:C6bar-free}
Characterize $\overline{C_6}$-free {\mcg}s.
\end{problem}

The following result of Kothari and Murty \cite{komu16}
shows that, for both problems stated above,
one may restrict their attention to bricks.
\begin{theorem}
\label{thm:J-free-reduction-to-bricks}
Let $J$ denote any cubic brick.
A matching covered graph~$G$ is $J$-free if and only if each of its
bricks is $J$-free.
\end{theorem}

In the same work~\cite{komu16}, they solved
both problems when the graph under
consideration is a planar brick.
However, these problems remain unsolved for nonplanar bricks.

\begin{problem}
\label{prb:K4-free-restriction}
Characterize $K_4$-free nonplanar bricks.
\end{problem}
\begin{problem}
\label{prb:C6bar-free-restriction}
Characterize $\overline{C_6}$-free nonplanar bricks.
\end{problem}

It should be noted that, unlike the Pfaffian graph recognition problem,
Problems \ref{prb:K4-free-restriction} and \ref{prb:C6bar-free-restriction}
are unsolved even for the restricted case of near-bipartite graphs.

\smallskip
Another important and related problem
is that of deciding whether or not a given \mcg~$G$ is solid.
Carvalho, Lucchesi and Murty \cite{clm04} showed that
$G$ is solid if and only if each of its bricks is solid.
Thus, as usual, one may restrict attention to bricks.
The same authors, in \cite{clm06}, proved that the only
planar solid bricks are the odd wheels.
Thus it remains to solve the problem for nonplanar bricks.

\begin{problem}
\label{prb:solid-restriction}
Characterize solid nonplanar bricks.
\end{problem}

Recently, Lucchesi, Carvalho, Kothari and Murty \cite{lckm18}
showed that Problems~\ref{prb:C6bar-free-restriction}
and~\ref{prb:solid-restriction} are in fact equivalent.

\begin{theorem}
\label{thm:equivalence-of-solid-and-C6bar-free}
A nonplanar brick~$G$ is solid if and only if it is $\overline{C_6}$-free
\underline{unless}
$G$ is the Petersen graph (up to multiple edges).
\end{theorem}

The Petersen graph is $\overline{C_6}$-free but it is not solid.

\section{Proof of the Main Theorem}
\label{sec:main-theorem-proof}
Our goal is to prove the
Main Theorem~(\ref{thm:efeccb-two-qbinv-at-vertex}).
We adopt the following notation.
\begin{notation}
\label{Not:v-claw}
Let $G$ be an \efeccb, let $v \in V(G)$, and let $e_1, e_2, e_3$ be the three edges incident with $v$.
For each $i \in \{1,2,3\}$, we let $e_i := vu_i$, and we let $s_i$ and $t_i$ denote the two neighbours of $u_i$ that are distinct from $v$.
\end{notation}

For $i \in \{1,2\}$, we assume that $e_i$ is \qbinv,
and we adopt the notation and conventions introduced
in Theorem~\ref{thm:efeccb-qbinv-edge} and
Notation~\ref{Not:8-edges} ---
with the only difference being that all of the notation (except for the vertex $v$)
is decorated with subscript $i$.
For instance, the two (nontrivial) barriers
of $G-e_i$ will be denoted as $B_i$ and $B'_i$
with the convention that the two neigbours
of $v$, that are distinct from $u_i$, lie in the barrier $B_i$.
Thus $s_i,t_i \in B'_i$.

\smallskip
Note that $u_2,u_3 \in B_1$,
and that $L_1$ and $L'_1$ are the two components of $G-X_1-X'_1$.
It follows from Theorem~\ref{thm:efeccb-qbinv-edge}(\ref{itm:L-Lp})
that each of $u_2$ and $u_3$ has at most one neighbour in $V(L_1)$,
and likewise, at most one neighbour in $V(L'_1)$.
Analogously, each of $u_1$ and $u_3$ has at most one neighbour
in $V(L_2)$, and likewise, at most one neighbour in $V(L'_2)$.
These observations, and associated notational conventions, are stated below.
\begin{notation}
\label{Not:s-and-t-vertices}
For each $i \in \{1,2\}$, and for each $j \in \{1,2,3\}$, where $i \neq j$,
each of the sets $V(L_i) \cap \{s_j,t_j\}$ and
$V(L'_i) \cap \{s_j,t_j\}$ is either empty or a singleton;
furthermore, if the former set is a singleton we let $s_j$ denote its
unique element, and if the latter set is a singleton we let $t_j$
denote its unique element.
\end{notation}

For instance, if the set $\{s_1,t_1\}$ has nonempty intersection
with each of $V(L_2)$ and $V(L'_2)$ then, as per above convention,
$s_1 \in V(L_2)$ and $t_1 \in V(L'_2)$.

\subsection{The barriers $B_1$ and $B_2$}

As a first step, we proceed to show that at least one of the two sets $B_1$ and $B_2$ is
a doubleton. (Eventually, we will establish that each of them is a doubleton.)

\begin{proposition}
\label{prp:at-least-one-of-B1-B2-doubleton}
At least one of the two sets $B_1$ and $B_2$ is a doubleton.
\end{proposition}
\begin{proof}
We begin by noting that $u_3 \in B_1 \cap B_2$ and that $N(u_3) = \{v,s_3,t_3\}$.
By Lemma~\ref{lem:intersection-of-barriers}, we infer that $B_1 \cap B_2 = \{u_3\}$,
and this immediately implies $I_1 \cap I_2 = \{v\}$.

\smallskip
Now suppose that each of $B_1$ and $B_2$ has cardinality at least three.
Since $|B_1| \geq 3$, the vertex $u_3$ must have a neighbour in $I_1-v$,
otherwise $B_1-u_3$ is a nontrivial barrier in $G$. Likewise, $u_3$ must
have a neighbour in $I_2-v$. These facts imply that the set $\{s_3,t_3\}$
meets each of $I_1$ and $I_2$. Since $I_1 \cap I_2 = \{v\}$,
exactly one of $s_3$ and $t_3$ lies in $I_1$ and the other lies in $I_2$.
Adjust notation so that $s_3 \in I_1$ and $t_3 \in I_2$.
Consequently, as per Notation~\ref{Not:s-and-t-vertices},
$t_3 \in V(L'_1)$.
By Theorem~\ref{thm:efeccb-qbinv-edge}(\ref{itm:4-cut-matching}),
$t_3$ has a neighbour, say $p$, in $L'_1$.
Since $t_3 \in I_2$, we infer that $p \in B_2$.
As $B_2$ is a barrier of $G-e_2$, it follows that 
the graph $G-e_2-u_3-p$ is not matchable. In the next two paragraphs,
we will contradict this fact by constructing a perfect matching
of $G-e_2-u_3-p$.

\smallskip
Let $M'$ be a perfect matching of $J'_1-p-x_1$.
Note that $M'$ contains exactly one of $g'_1$ and $h'_1$.
Adjust notation so that $g'_1 \in M'$, and let $z'$
denote the end of $g'_1$ in $B'_1$.
By Theorem~\ref{thm:efeccb-qbinv-edge}(\ref{itm:L-Lp}),
$d'_1$ and $f'_1$ are nonadjacent, whence at least one of them
is not incident with~$z'$. Adjust notation so that
$f'_1 \notin \partial(z')$, and let $y'$ denote the end of $f'_1$ in $B'_1$.
Thus $y' \neq z'$. Now, let $M$ be a perfect matching of $J_1$
such that $f'_1 \in M$. Exactly one of $d_1$ and $f_1$ lies in $M$.
Adjust notation so that $f_1 \in M$ and let $y$ denote the end of $f_1$
in $B_1$. Note that two neighbours of $u_3$ lie in $I_1$,
whereas the third neighbour lies in $L'_1$. Thus $u_3 \neq y$.

\smallskip
Now we invoke Lemma~\ref{lem:bipartite-shore-properties} twice.
Let $N$ be a perfect matching of $G[X_1-v-u_3-y]$, and
let $N'$ be a perfect matching of $G[X'_1-u_1-y'-z']$.
Observe that $M \cup M' \cup N \cup N' \cup \{e_1\}$
is indeed a perfect matching of $G-e_2-u_3-p$. This contradicts
what we have established earlier.
We thus conclude that at least one of $B_1$ and $B_2$ has
cardinality precisely two, and this completes the proof of
Proposition~\ref{prp:at-least-one-of-B1-B2-doubleton}.
\end{proof}

\begin{figure}[!htb]
    \centering
    \begin{tikzpicture}[scale=0.65]

      %edges d and f
      \draw (-0.5,4) to [out=180,in=90] (-4,0.5)node{};
      \draw (-4.5,0.5)node[nodelabel]{$s_2$};
      \draw (-3.3,3.3)node[nodelabel]{$d_1$};
      \draw (0.5,4) to [out=220,in=90] (-3,0.5)node{};
      \draw (-2.5,0.5)node[nodelabel]{$s_3$};
      \draw (-1.7,2.1)node[nodelabel]{$f_1$};

      %edges g and h
      \draw (0.5,4) to [out=0,in=90] (4,0.5)node{};
      \draw (4.5,0.5)node[nodelabel]{$t_3$};
      \draw (3.3,3.3)node[nodelabel]{$h_1$};
      \draw (-0.5,4) to [out=320,in=90] (3,0.5)node{};
      \draw (2.5,0.5)node[nodelabel]{$t_2$};
      \draw (1.7,2.1)node[nodelabel]{$g_1$};

      %edge e=uv
      \draw (0,2) -- (0,-2);
      \draw (0.4,0)node[nodelabel]{$e_1$};

      %other edges incident at v
      \draw (0,2) -- (-0.5,4)node{};
      \draw (0,2) -- (0.5,4)node{};

      %other edges incident at u
      \draw (0,-2) -- (-0.5,-4)node{};
      \draw (-0.9,-4)node[nodelabel]{$s_1$};
      \draw (0,-2) -- (0.5,-4)node{};
      \draw (0.9,-4)node[nodelabel]{$t_1$};

      %barrier B
      %\draw (-1,3.5) -- (1,3.5) -- (1,4.5) -- (-1,4.5) -- (-1,3.5);
      \draw (-0.5,4.8)node[nodelabel]{$u_2$};
      \draw (0.5,4.8)node[nodelabel]{$u_3$};

      %barrier B'
      \draw (-3,-3.5) -- (3,-3.5) -- (3,-4.5) -- (-3,-4.5) -- (-3,-3.5);
      \draw (-3.5,-4)node[nodelabel]{$B'_1$};

      %isolated set I'-u
      \draw (-2,-5.5) -- (2,-5.5) -- (2,-6.5) -- (-2,-6.5) -- (-2,-5.5);
      \draw (-3,-6)node[nodelabel]{$I'_1-u_1$};

      %component L
      \draw (-2,1.5) -- (-5,1.5) -- (-5,-1.5) -- (-2,-1.5) -- (-2,1.5);
      \draw (-5.4,0)node[nodelabel]{$L_1$};

      %component L'
      \draw (2,1.5) -- (5,1.5) -- (5,-1.5) -- (2,-1.5) -- (2,1.5);
      \draw (5.5,0)node[nodelabel]{$L'_1$};

      %edges d' and f'
      \draw (-2.8,-3.2) -- (-3.9,-1.8);
      \draw (-3.9,-2.5)node[nodelabel]{$d'_1$};
      \draw (-2,-3.2) -- (-3.1,-1.8);
      \draw (-2,-2.5)node[nodelabel]{$f'_1$};

      %edges g' and h'
      \draw (2.8,-3.2) -- (3.9,-1.8);
      \draw (3.9,-2.5)node[nodelabel]{$h'_1$};
      \draw (2,-3.2) -- (3.1,-1.8);
      \draw (2,-2.5)node[nodelabel]{$g'_1$};

      %ends u and v of edge e
      \draw (0,2)node{};
      \draw (0.3,1.8)node[nodelabel]{$v$};
      \draw (0,-2)node{};
      \draw (0.4,-1.8)node[nodelabel]{$u_1$};

      %label G
      %\draw (0,-9)node[nodelabel]{$G$};

    \end{tikzpicture}
    \caption{The set $B_1$ is a doubleton}
    \label{fig:B1-doubleton}
  \end{figure}
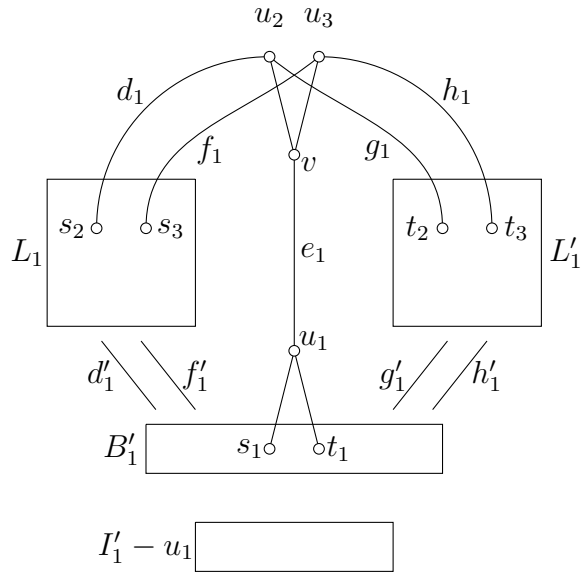

We adjust notation so that $B_1=\{u_2,u_3\}$, whence $I_1=\{v\}$.
Consequently, as per the conventions stated in Notation~\ref{Not:s-and-t-vertices},
vertices $s_2$ and $s_3$ lie in $L_1$, and $t_2$ and $t_3$ lie in $L'_1$.
We adjust notation so that $d_1,g_1 \in \partial(u_2)$
and $f_1,h_1 \in \partial(u_3)$.
See Figure~\ref{fig:B1-doubleton}.
Note that, among other things,
Theorem~\ref{thm:efeccb-qbinv-edge}(\ref{itm:L-Lp})
implies that
all of the six vertices $s_1, t_1, s_2, t_2, s_3, t_3$ are pairwise distinct,
whence $G$ has order at least $10$.

\subsection{Bricks isomorphic to $K_4$}

Let $\mathcal{J}:=\{J_1,J'_1,J_2,J'_2\}$,
and let $\mathcal{L}:=\{L_1,L'_1,L_2,L'_2\}$. Observe that a brick $J \in \mathcal{J}$
is isomorphic to $K_4$ (up to multiple edges) if and only if the corresponding (connected)
graph $L \in \mathcal{L}$
is isomorphic to~$K_2$.
Our next goal is to prove that each member of $\mathcal{J}$ is isomorphic to $K_4$
(up to multiple edges)\footnote{From now on, we shall abuse terminology slightly and just write `isomorphic to $K_4$' instead of  writing `isomorphic to $K_4$ (up to multiple edges)'.}.
However, to do so, we will require several auxiliary lemmas.

\smallskip
By assuming that the set $B_1$ is a doubleton, we have lost the symmetry
between edges $e_1$~and~$e_2$. We will thus find it convenient to first prove
that each of $J_1$ and $J'_1$ is isomorphic to~$K_4$. The following lemma
shows why doing so is in fact sufficient.

\begin{lemma}
\label{lem:J1-K4-implies-J2-K4}
If $J_1 \iso K_4$ then $J_2 \iso K_4$ and $V(L_1) \cap V(L_2)=\{s_3\}$.
Likewise,
if $J'_1 \iso K_4$ then $J'_2 \iso K_4$ and $V(L'_1) \cap V(L'_2) = \{t_3\}$.
\end{lemma}
\begin{proof}
By symmetry, it suffices to prove the first statement. Assume that $J_1 \iso K_4$,
whence $E(L_1)=\{s_2s_3\}$. Observe that $s_3$ is a common neighbour
of $u_3 \in B_2$ and $s_2 \in B'_2$. Consequently, $s_3 \in V(L_2) \cup V(L'_2)$,
and as per Notation~\ref{Not:s-and-t-vertices}, $s_3 \in V(L_2)$.
Now, since $\partial(V(L_2))$ is not a matching,
Theorem~\ref{thm:efeccb-qbinv-edge}(\ref{itm:4-cut-matching})
implies that $L_2 \iso K_2$, whence $J_2 \iso K_4$.
Note that $V(L_1) \cap V(L_2) = \{s_3\}$.
This completes the proof of Lemma~\ref{lem:J1-K4-implies-J2-K4}.
\end{proof}

\begin{corollary}
\label{cor:B'1-doubleton-Petersen}
If at least one of $J_1$ and $J'_1$ is isomorphic to $K_4$,
and if the set $B'_1$ is a doubleton,
then $G$ is the Petersen graph.
\end{corollary}
\begin{proof}
Assume that $J_1$ is isomorphic to $K_4$ and
that $|B'_1|=2$.
Since $J_1 \iso K_4$,
Lemma~\ref{lem:J1-K4-implies-J2-K4} implies that
$J_2 \iso K_4$ and $V(L_1) \cap V(L_2) = \{s_3\}$.
There are two edges joining the sets
$B'_1=\{s_1,t_1\}$ and $V(L_1)=\{s_2,s_3\}$.
Consequently, one of $s_1$ and $t_1$ lies in $V(L_2)$,
and as per Notation~\ref{Not:s-and-t-vertices},
$E(L_2) = \{s_1s_3\}$. Thus, $t_1s_2 \in E(G)$,
whence $t_1$ is a common neighbour of $u_1 \in B_2$ and $s_2 \in B'_2$.
It follows that $t_1 \in V(L'_2)$, and since $\partial(V(L'_2))$ is not a matching,
Theorem~\ref{thm:efeccb-qbinv-edge}(\ref{itm:4-cut-matching})
implies that $L'_2 \iso K_2$.
Now, we observe that $u_1$ has no neighbours in $I_2 - v$,
and that $s_2$ has no neighbours in $I'_2-u_2$.
By Lemma~\ref{lem:bipartite-shore-properties},
$B_2 = \{u_1,u_3\}$ and $B'_2 = \{s_2,t_2\}$.
Consequently, $t_2s_1 \in E(G)$, and each of $u_3$ and $t_2$
is adjacent with the unique vertex in $V(L'_2)-t_1$.
As per Notation~\ref{Not:v-claw}, $V(L'_2) = \{t_1,t_3\}$;
whence $t_1t_3, t_2t_3 \in E(G)$.
We now have a cubic subgraph of $G$ that is isomorphic to the
Petersen graph. Thus $G$ is indeed the Petersen graph.
\end{proof}

\begin{corollary}
\label{cor:B'2-doubleton-Petersen}
If at least one of $J_1$ and $J'_1$ is isomorphic to $K_4$,
and if the set $B'_2$ is a doubleton,
then $G$ is the Petersen graph.
\end{corollary}
\begin{proof}
Assume that $J_1$ is isomorphic to $K_4$ and that $|B'_2|=2$.
Since $J_1 \iso K_4$, Lemma~\ref{lem:J1-K4-implies-J2-K4}
implies that $J_2 \iso K_4$ and $V(L_1) \cap V(L_2) = \{s_3\}$.
Since $B'_2=\{s_2,t_2\}$, and $s_2s_3 \in E(G)$,
the unique vertex in $V(L_2)-s_3$ is a common neighbour of $t_2$ and $s_3$.
Consequently, $\partial(V(L'_1))$ is not a matching;
by Theorem~\ref{thm:efeccb-qbinv-edge}(\ref{itm:4-cut-matching}),
$E(L'_1) = \{t_2t_3\}$ and $J'_1 \iso K_4$.
Lemma~\ref{lem:J1-K4-implies-J2-K4} implies
that $J'_2 \iso K_4$ and $V(L'_1) \cap V(L'_2) = \{t_3\}$.
Now, observe that $u_3$ has no neighbours in $I_2-v$,
whence it follows from Lemma~\ref{lem:bipartite-shore-properties}
that $B_2=\{u_1,u_3\}$.
Consequently, $s_1,t_1 \in V(L_2) \cup V(L'_2)$, and as
per Notation~\ref{Not:s-and-t-vertices}, $s_1 \in V(L_2)$ and $t_1 \in V(L'_2)$.
In particular, $s_1s_3,t_1t_3 \in E(G)$. Finally, note that $s_2t_1 \in E(G)$.
We now have a cubic subgraph of $G$ that is isomorphic to the
Petersen graph. Thus $G$ is indeed the Petersen graph.
\end{proof}

We now embark on the arduous journey of proving that each of the bricks $J_1$
and $J'_1$ is indeed isomorphic to $K_4$.

\begin{lemma}
\label{lem:intersection-singleton-implies-K4}
If $V(L_1) \cap (B_2 \cup B'_2) = \{s_2\}$ then $J_1 \iso K_4$.
Likewise, if $V(L'_1) \cap (B_2 \cup B'_2) = \{t_2\}$ then $J'_1 \iso K_4$.
\end{lemma}
\begin{proof}
By symmetry, it suffices to prove the first statement.
Assume, to the contrary, that $V(L_1) \cap (B_2 \cup B'_2) = \{s_2\}$
and that $J_1$ is not isomorphic to~$K_4$.
By Theorem~\ref{thm:efeccb-qbinv-edge}(\ref{itm:4-cut-matching}),
$L_1$ is a $2$-connected graph,
whence each of its vertices has at least two distinct neighbours.
This fact, along with the hypothesis $V(L_1) \cap (B_2 \cup B'_2) = \{s_2\}$,
implies that $V(L_1) \cap (I_2 \cup I'_2) = \emptyset$.
Thus $V(L_1) \cap (X_2 \cup X'_2) = \{s_2\}$.
Consequently, $L_1-s_2$ is a connected subgraph of $G-X_2-X'_2$,
whence it is either a subgraph of~$L_2$ or of~$L'_2$.
Since $s_3 \in V(L_1)$, as per Notation~\ref{Not:s-and-t-vertices},
$L_1-s_2$ is a subgraph of $L_2$.

\smallskip
Now, since $s_2$ has two distinct neighbours in $L_1$,
there exist two edges joining $s_2 \in X'_2$ and $L_2$,
a contradiction to
Theorem~\ref{thm:efeccb-qbinv-edge}(\ref{itm:L-Lp}).
This proves Lemma~\ref{lem:intersection-singleton-implies-K4}.
\end{proof}

For $i \in \{1,2\}$, we say that the brick $J_i$ is {\it flexible}
if none of the four edges in the set~$\partial(V(L_i))$ depends on any of the
other three edges; otherwise, we say that the brick $J_i$ is {\it inflexible}.
Analogous definitions apply to the brick $J'_i$
--- with $L'_i$ playing the role of $L_i$.

\smallskip
Observe that, as per the above definitions:
any member of $\mathcal{J}$, that is isomorphic to $K_4$,
is in fact inflexible.
As a first step, we will prove that $J_1$ and $J'_1$ are both inflexible;
to do so, we will need two lemmas.

\begin{lemma}
\label{lem:if-either-flexible-then-the-other-K4}
If $J_1$ is flexible then $J'_1 \iso K_4$.
Likewise, if $J'_1$ is flexible then $J_1 \iso K_4$.
\end{lemma}
\begin{proof}
By symmetry, it suffices to prove the first statement.
Assume that $J_1$ is flexible.
Observe that $t_2 \in V(L'_1) \cap B'_2$.
Our plan is to leverage the flexibility of $J_1$ to
deduce that $V(L'_1) \cap (B_2 \cup B'_2) = \{t_2\}$,
and then invoke Lemma~\ref{lem:intersection-singleton-implies-K4}.
\begin{assertion}
\label{sta-1:if-either-flexible-then-the-other-K4}
For each $p \in V(L'_1)-t_2$, the graph $G-e_2-p-t_2$ is matchable;
consequently, $p \notin B'_2$.
\end{assertion}
\begin{proof}[Proof of Assertion~$\ref{sta-1:if-either-flexible-then-the-other-K4}$]
We let $M'$ denote a perfect matching of $J'_1 - p - t_2$.
First suppose that $x_1x'_1 \in M'$. Consequently, $M'-x_1x'_1$
is a perfect matching of $L'_1 - p - t_2$.
By Lemma~\ref{lem:L-p-q-matchable-e'-admissible-in-G-p-q},
with $e_3$ playing the role of $e'$,
the graph $G-e_2-p-t_2$
is matchable.

\smallskip
Now suppose that $x_1x'_1 \notin M'$. Thus $M'$
contains $h_1=u_3t_3$, and it also contains exactly one of $g'_1$~and~$h'_1$.
Adjust notation so
that $h'_1 \in M'$.
By Theorem~\ref{thm:efeccb-qbinv-edge}(\ref{itm:L-Lp}),
at least one of $d'_1$~and~$f'_1$ is not adjacent with $h'_1$.
Adjust notation so that $d'_1$ and $h'_1$ are nonadjacent.
Now, we utilize the flexibility hypothesis,
to choose a perfect matching $M$ of $J_1$
that contains both edges $d'_1$ and $d_1=u_2s_2$.
Observe that $M'$ and $M$ satisfy the hypotheses
of Lemma~\ref{lem:extend-M-union-Mp-union-e},
and thus $M' \cup M \cup \{e_1\}$
may be extended to a perfect matching of $G-p-t_2$.
Consequently, $G-e_2-p-t_2$ is matchable.
This proves Assertion~\ref{sta-1:if-either-flexible-then-the-other-K4}.
\end{proof}
\begin{assertion}
\label{sta-2:if-either-flexible-then-the-other-K4}
For each $p \in V(L'_1)$, the graph $G-e_2-p-u_3$ is matchable;
consequently, $p \notin B_2$.
\end{assertion}
\begin{proof}[Proof of Assertion~$\ref{sta-2:if-either-flexible-then-the-other-K4}$]
We let $M'$ denote a perfect matching of $J'_1 - p - x_1$.
Adjust notation so that $h'_1 \in M'$, and let $z'$ denote the end of $h'_1$ in $B'_1$.
Adjust notation so that $d'_1$ and $h'_1$ are nonadjacent, and let $y'$
denote the end of $d'_1$ in $B'_1$. (Thus $y' \neq z'$.) As before,
we make use of the flexibility hypothesis, to choose a perfect matching $M$ of $J_1$
that contains both edges $d'_1$ and $d_1=u_2s_2$.
Now, we invoke Lemma~\ref{lem:bipartite-shore-properties}
to choose a perfect matching $N'$ of $G[X'_1-u_1-y'-z']$.
Observe that $M \cup M' \cup N' \cup \{e_1\}$ is a perfect matching
of $G-e_2-p-u_3$. This proves Assertion~\ref{sta-2:if-either-flexible-then-the-other-K4}.
\end{proof}

It follows from
Assertions~\ref{sta-1:if-either-flexible-then-the-other-K4}~and~\ref{sta-2:if-either-flexible-then-the-other-K4}
that $V(L'_1) \cap (B_2 \cup B'_2) = \{t_2\}$.
Lemma~\ref{lem:intersection-singleton-implies-K4}
implies that $J'_1 \iso K_4$.
This completes the proof of
Lemma~\ref{lem:if-either-flexible-then-the-other-K4}.
\end{proof}

\begin{lemma}
\label{lem:both-K4-or-neither-K4}
Either each of $J_1$ and $J'_1$ is isomorphic to $K_4$,
or otherwise neither of them is isomorphic to $K_4$.
\end{lemma}
\begin{proof}
Assume that $J'_1 \iso K_4$, whence $E(L'_1) = \{t_2t_3\}$.
By Lemma~\ref{lem:J1-K4-implies-J2-K4},
$J'_2 \iso K_4$ and $V(L'_1) \cap V(L'_2)=\{t_3\}$.

\smallskip
If either of the sets $B'_1$ and $B'_2$ is a doubleton,
we invoke
Corollary~\ref{cor:B'1-doubleton-Petersen} or
Corollary~\ref{cor:B'2-doubleton-Petersen} (as applicable)
to deduce that $G$ is the Petersen graph; in particular, $J_1 \iso K_4$.

\smallskip
Now suppose that each of the sets $B'_1$ and $B'_2$ has three or more vertices.
We will consider two cases depending on whether or not the set $B_2$ is a doubleton.
(In each case, we will deduce that $J_1 \iso K_4$.)

\bigskip
\noindent
{\bf Case 1:} The set $B_2$ is a doubleton.

\smallskip
\noindent
Since $B_2=\{u_1,u_3\}$, each of $u_1$ and $u_3$ has
one neighbour in $L_2$, and one neighbour in~$L'_2$.
As per Notation~\ref{Not:s-and-t-vertices},
$L_2$ contains $s_1$ and $s_3$,
and $L'_2$ contains $t_1$ and $t_3$.
In particular, $E(L'_2)=\{t_1t_3\}$.
Note that one of $g'_1$ and $h'_1$ is the edge $t_1t_3$.
Since $|B'_1| \geq 3$, the vertex $t_1$ has a neighbour in $I'_1-u_1$,
whence $d'_1,f'_1 \notin \partial(t_1)$.
\begin{assertion}
\label{sta-1:both-K4-or-neither-K4}
For each $p \in V(L_1)-s_2$, the graph $G-e_2-p-s_2$ is matchable;
consequently,~$p \notin B'_2$.
\end{assertion}
\begin{proof}[Proof of Assertion~$\ref{sta-1:both-K4-or-neither-K4}$]
We let $M$ denote a perfect matching of $J_1-p-s_2$.
First suppose that $x_1x'_1 \in M$,
whence $M-x_1x'_1$ is a perfect matching of $L_1 - p -s_2$.
By Lemma~\ref{lem:L-p-q-matchable-e'-admissible-in-G-p-q},
with $e_3$ playing the role of $e'$,
the graph $G-e_2-p-s_2$ is matchable.

\smallskip
Now suppose that $x_1x'_1 \notin M$.
Thus $M$ contains $f_1=u_3s_3$, and it also contains exactly one of $d'_1$ and $f'_1$.
Let $M':=\{u_2t_2,t_1t_3\}$.
Since $d'_1,f'_1 \notin \partial(t_1)$,
the matchings $M$ and $M'$ satisfy the hypotheses of
Lemma~\ref{lem:extend-M-union-Mp-union-e},
whence $M \cup M' \cup \{e_1\}$ may be extended to a perfect
matching of $G-p-s_2$. Consequently, $G-e_2-p-s_2$ is matchable.
This proves Assertion~\ref{sta-1:both-K4-or-neither-K4}.
\end{proof}

It follows from Assertion~\ref{sta-1:both-K4-or-neither-K4} that $V(L_1) \cap B'_2 = \{s_2\}$.
Clearly, $V(L_1) \cap B_2 = \emptyset$.
Consequently, $V(L_1) \cap (B_2 \cup B'_2) = \{s_2\}$.
By Lemma~\ref{lem:intersection-singleton-implies-K4},
$J_1 \iso K_4$.

\bigskip
\noindent
{\bf Case 2:} The set $B_2$ has three or more vertices.

\smallskip
\noindent
Since each of $B_2$ and $B'_2$ has cardinality at least three,
by Lemma~\ref{lem:bipartite-shore-properties},
each of the bipartite graphs $G[X_2-v]$ and $G[X'_2-u_2]$ is connected.
We will investigate these two graphs, and use the fact that
$\partial(V(L_1)) = \{d_1,f_1,d'_1,f'_1\}$ is a cut of $G$, to arrive at a contradiction.

\smallskip
Note that $d_1=u_2s_2$ is neither an edge of $G[X_2-v]$ nor of $G[X'_2-u_2]$.
Since $t_3 \in V(L'_2)$, and $|B_2| \geq 3$, the vertex $s_3$ lies in $I_2-v$.
Consequently, $f_1=u_3s_3$ is an edge of $G[X_2-v]$.

\smallskip
The connected graph $G[X'_2-u_2]$ contains $s_2 \in V(L_1)$ and $t_2 \notin V(L_1)$;
in other words, it meets each shore of the cut $\{d_1,f_1,d'_1,f'_1\}$;
whence it contains at least one of these four edges. It follows from the 
preceding paragraph that at least
one of $d'_1$ and $f'_1$ is an edge of $G[X'_2-u_2]$.
Adjust notation so that $f'_1$ is an edge of $G[X'_2-u_2]$.

\smallskip
Now,
we observe that $\{f_1\}$ is a $1$-cut of $G[X_2-v]$;
furthermore, $G[X_2-v]-f_1$
has precisely two components: the isolated vertex $u_3$, and a nontrivial component,
say $Q$. Note that $Q=G[X_2-v] - u_3$.
The connected graph $Q$ contains $s_3 \in V(L_1)$
and $u_1 \notin V(L_1)$; whence $Q$
contains at least one edge from $\{d_1,f_1,d'_1,f'_1\}$.
We infer that $d'_1 \in E(Q)$.

\smallskip
Thus, the connected subgraph $G[X'_2-u_2]$ contains exactly one edge
from the cut $\{d_1,f_1,d'_1,f'_1\}$ --- namely, $f'_1$.
Consequently, $\{f'_1\}$
is a $1$-cut of $G[X'_2-u_2]$; furthermore, one of its shores
contains $s_2$, and the other shore contains $t_2$.
By Lemma~\ref{lem:bipartite-shore-properties},
each $1$-cut of $G[X'_2-u_2]$ is trivial.
This implies that $f'_1$ is incident with one of $s_2$ and $t_2$.
Since $t_2 \in V(L'_1)$, it is not an end of $f'_1$.
Hence, $f'_1 \in \partial(s_2)$. Consequently, $\partial(V(L_1))$
is not a matching. By Theorem~\ref{thm:efeccb-qbinv-edge}(\ref{itm:4-cut-matching}),
$L_1 \iso K_2$, whence $J_1 \iso K_4$.
\end{proof}

The following is an immediate consequence of
Lemmas~\ref{lem:if-either-flexible-then-the-other-K4}~
and~\ref{lem:both-K4-or-neither-K4}.

\begin{corollary}
Each of $J_1$ and $J'_1$ is inflexible. \qed
\end{corollary}

Since the brick $J_1$ is inflexible, by definition,
some edge in $\partial(V(L_1))=\{d_1,f_1,d'_1,f'_1\}$
depends on another edge in $\partial(V(L_1))$.
Since any pair of such edges must be nonadjacent,
we may adjust notation so that $d_1$ depends on $d'_1$.
By Theorem~\ref{thm:efeccb-qbinv-edge}(\ref{itm:3-cut-x-xp}),
the edge $x_1x'_1$ participates in every nontrivial $3$-cut of $J_1$.
Thus we may invoke Corollary~\ref{cor:dependence-implies-doubleton}
to infer that each of the sets $\{d_1,d'_1\}$ and $\{f_1,f'_1\}$
is a removable doubleton of $J_1$. An analogous argument applies
to~$J'_1$. These facts, and notational conventions,
are summarized below. See Figure~\ref{fig:J1-J'1-both-inflexible}.
\begin{notation}
\label{Not:J1-J'1-inflexible}
Each of the sets
$\{d_1,d'_1\}$ and $\{f_1,f'_1\}$ is a removable doubleton
of $J_1$,
and likewise,
each of the sets
$\{g_1,g'_1\}$ and $\{h_1,h'_1\}$ is a removable doubleton
of $J'_1$.
\end{notation}

\begin{proposition}
\label{prp:J1-J'1-both-K4}
Each of the bricks $J_1$ and $J'_1$ is isomorphic to $K_4$.
In particular, $E(L_1) = \{s_2s_3\}$ and $E(L'_1)=\{t_2t_3\}$.
\end{proposition}
\begin{proof}
By Lemma~\ref{lem:both-K4-or-neither-K4},
it suffices to show that one of $J_1$ and $J'_1$ is isomorphic to $K_4$.
Suppose, to the contrary, that neither $J_1$ nor $J'_1$ is isomorphic to $K_4$.
Consequently, each of $L_1$~and~$L'_1$ has four or more vertices;
by Theorem~\ref{thm:efeccb-qbinv-edge}(\ref{itm:4-cut-matching}),
each of them is $2$-connected;
furthermore,
each of the sets $\partial(V(L_1))$ and $\partial(V(L'_1))$
is a matching in $G$.

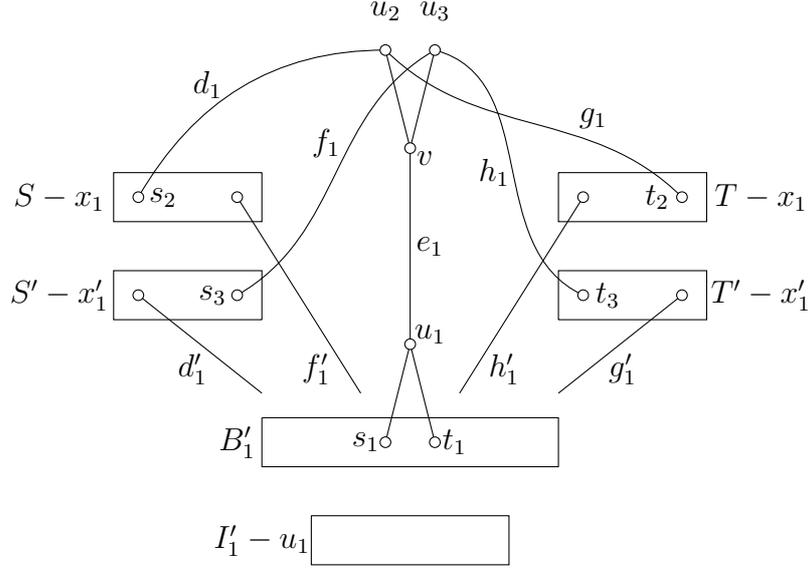
\begin{figure}[!htb]
    \centering
    \begin{tikzpicture}[scale=0.65]

      %edges d and f
      \draw (-0.5,4) to [out=180,in=60] (-5.5,1)node{};
      \draw (-5,1)node[nodelabel]{$s_2$};
      \draw (0.5,4) to [out=210,in=30] (-3.5,-1)node{};
      \draw (-4,-1)node[nodelabel]{$s_3$};
      \draw (-4.1,3.3)node[nodelabel]{$d_1$};
      \draw (-1.7,2.1)node[nodelabel]{$f_1$};

      %edges g and h
      \draw (-0.5,4) to [out=315,in=135] (5.5,1)node{};
      \draw (5,1)node[nodelabel]{$t_2$};
      \draw (0.5,4) to [out=345,in=155] (3.5,-1)node{};
      \draw (4,-1)node[nodelabel]{$t_3$};
      \draw (3.7,2.6)node[nodelabel]{$g_1$};
      \draw (1.7,1.5)node[nodelabel]{$h_1$};

      %edge e=uv
      \draw (0,2) -- (0,-2);
      \draw (0.4,0)node[nodelabel]{$e_1$};

      %other edges incident at v
      \draw (0,2) -- (-0.5,4)node{};
      \draw (0,2) -- (0.5,4)node{};

      %other edges incident at u
      \draw (0,-2) -- (-0.5,-4)node{};
      \draw (-0.9,-4)node[nodelabel]{$s_1$};
      \draw (0,-2) -- (0.5,-4)node{};
      \draw (0.9,-4)node[nodelabel]{$t_1$};

      %barrier B
      %\draw (-1,3.5) -- (1,3.5) -- (1,4.5) -- (-1,4.5) -- (-1,3.5);
      \draw (-0.5,4.8)node[nodelabel]{$u_2$};
      \draw (0.5,4.8)node[nodelabel]{$u_3$};

      %barrier B'
      \draw (-3,-3.5) -- (3,-3.5) -- (3,-4.5) -- (-3,-4.5) -- (-3,-3.5);
      \draw (-3.5,-4)node[nodelabel]{$B'_1$};

      %isolated set I'-u
      \draw (-2,-5.5) -- (2,-5.5) -- (2,-6.5) -- (-2,-6.5) -- (-2,-5.5);
      \draw (-3,-6)node[nodelabel]{$I'_1-u_1$};

      %component L
      \draw (-3,1.5) -- (-6,1.5) -- (-6,0.5) -- (-3,0.5) -- (-3,1.5);
      \draw (-3,-1.5) -- (-6,-1.5) -- (-6,-0.5) -- (-3,-0.5) -- (-3,-1.5);
      \draw (-7.1,1)node[nodelabel]{$S-x_1$};
      \draw (-7.1,-1)node[nodelabel]{$S'-x'_1$};

      %component L'
      \draw (3,1.5) -- (6,1.5) -- (6,0.5) -- (3,0.5) -- (3,1.5);
      \draw (3,-1.5) -- (6,-1.5) -- (6,-0.5) -- (3,-0.5) -- (3,-1.5);
      \draw (7.1,1)node[nodelabel]{$T-x_1$};
      \draw (7.1,-1)node[nodelabel]{$T'-x'_1$};

      %edges d' and f'
      \draw (-5.5,-1)node{} -- (-3,-3);
      \draw (-3.5,1)node{} -- (-1,-3);
      \draw (-4.4,-2.5)node[nodelabel]{$d'_1$};
      \draw (-1.9,-2.5)node[nodelabel]{$f'_1$};

      %edges g' and h'
      \draw (5.5,-1)node{} -- (3,-3);
      \draw (3.5,1)node{} -- (1,-3);
      \draw (4.3,-2.5)node[nodelabel]{$g'_1$};
      \draw (1.9,-2.5)node[nodelabel]{$h'_1$};

      %ends u and v of edge e
      \draw (0,2)node{};
      \draw (0.3,1.8)node[nodelabel]{$v$};
      \draw (0,-2)node{};
      \draw (0.4,-1.8)node[nodelabel]{$u_1$};

      %label G
      %\draw (0,-9)node[nodelabel]{$G$};

    \end{tikzpicture}
    \caption{The bricks $J_1$ and $J'_1$ are both inflexible}
    \label{fig:J1-J'1-both-inflexible}
  \end{figure}

\smallskip
We let $S$ and $S'$ denote the color classes of $J_1 -d_1-d'_1$
so that $d_1$ has both ends in $S$.
Consequently, $x_1,s_2 \in S$, and $x'_1,s_3 \in S'$.
Likewise, we let $T$ and $T'$ denote the color classes of $J'_1-g_1-g'_1$
so that $g_1$ has both ends in $T$.
Consequently, $x_1, t_2 \in T$, and $x'_1, t_3 \in T'$.
See Figure~\ref{fig:J1-J'1-both-inflexible}.

\begin{assertion}
\label{sta-1:J1-J'1-both-K4}
The edges $d'_1$ and $g'_1$ are nonadjacent.
\end{assertion}
\begin{proof}[Proof of Assertion~$\ref{sta-1:J1-J'1-both-K4}$]
We first claim that $s_3$ has a neighbour, distinct from $u_3$, that does
not lie in $B'_2$. Suppose, to the contrary, that each neighbour of $s_3$,
distinct from $u_3$, belongs to $B'_2$. Consequently,
$s_3$ has two neighbours in $B'_2$, and it has one neighbour (namely $u_3$) in $B_2$.
Thus, $s_3 \in V(L_2) \cup V(L'_2)$, and as per Notation~\ref{Not:s-and-t-vertices},
$s_3 \in V(L_2)$. Furthermore, there are two adjacent edges joining $L_2$
and $X'_2$, contrary to Theorem~\ref{thm:efeccb-qbinv-edge}(\ref{itm:L-Lp}).
Thus $s_3$ has a neighbour, say $w$, such that $w \neq u_3$ and $w \notin B'_2$.

\smallskip
By Theorem~\ref{thm:efeccb-qbinv-edge}(\ref{itm:L-Lp}),
$B'_2$ is a maximal barrier of $G-e_2$.
Since $s_2 \in B'_2$ and $w \notin B'_2$,
by Theorem~\ref{thm:canonical-partition},
these two vertices do not lie in a common barrier of $G-e_2$.
Thus, $G-e_2-w-s_2$ has a perfect matching, say $N$.
Note that, the graph $L_1$ is bipartite and has equicardinal color
classes, namely $S-x_1$ and $S'-x'_1$. Also, $w$ and $s_2$ both lie in~$S-x_1$.
Since $S-x_1$ is a stable set, we infer that $d'_1$ and $f_1=u_3s_3$
belong to $N$.

\smallskip
Observe that $g_1=u_2t_2$ belongs to $N$.
Since, $f_1 \in N$, the edge $h_1=u_3t_3$ does not lie in~$N$.
Since $L'_1$ is also bipartite and has equicardinal color classes,
namely $T-x_1$ and $T'-x'_1$, we infer that $g'_1 \in N$.

\smallskip
In particular, we have shown that $d'_1,g'_1 \in N$.
Thus they are nonadjacent.
\end{proof}

We let $y'$ and $z'$, in the set $B'_1$, denote the ends of $d'_1$ and $g'_1$,
respectively. By Assertion~\ref{sta-1:J1-J'1-both-K4}, $y'$~and~$z'$ are distinct.

\begin{assertion}
\label{sta-2:J1-J'1-both-K4}
For each $p \in V(L_1)-s_2$,
the graph $G-e_2-p-s_2$ is matchable;
consequently,~$p \notin B'_2$.
\end{assertion}
\begin{proof}[Proof of Assertion~$\ref{sta-2:J1-J'1-both-K4}$]
As noted earlier,
$L_1$ is a bipartite graph with color classes $S-x_1$ and $S'-x'_1$.
We let $p \in V(L_1)-s_2$.

\smallskip
First suppose that $p \in S'$.
By Lemma~\ref{lem:near-bipartite-J-admissibility-of-xx'},
the edge $x_1x'_1$ is admissible in the graph $J_1-d_1-d'_1-p-s_2$.
Consequently, $L_1-p-s_2$ is matchable.
By
Lemma~\ref{lem:L-p-q-matchable-e'-admissible-in-G-p-q},
with the edge $e_3$ playing the role of $e'$,
the graph $G-e_2-p-s_2$ is matchable, and we are done.

\smallskip
Now suppose that $p \in S$. Let $M$ denote a perfect matching of $J_1-p-s_2$.
Since $p$~and~$s_2$ both belong to the color class $S$ of
the bipartite graph $J_1-d_1-d'_1$,
the edge $d'_1$ lies in $M$.
Now, let $M'$ denote a perfect matching of $J'_1$
that contains the edge $g_1=u_2t_2$.
Since $\{g_1,g'_1\}$ is a removable doubleton of $J'_1$, the edge $g'_1$ lies in $M'$.

\smallskip
Since $y' \neq z'$, the matchings $M$ and $M'$ satisfy the hypotheses
of Lemma~\ref{lem:extend-M-union-Mp-union-e},
with $s_2$ playing the role of $q$.
Consequently, with $e_1$ playing the role of $e$,
we may extend
$M \cup M' \cup \{e_1\}$ to a perfect matching of $G-p-s_2$.
Thus, $G-e_2-p-s_2$ is matchable, and we are done.
\end{proof}

\begin{assertion}
\label{sta-3:J1-J'1-both-K4}
For each $p \in V(L_1)$, the graph $G-e_2-p-u_3$ is matchable;
consequently,~$p \notin B_2$.
\end{assertion}
\begin{proof}[Proof of Assertion~$\ref{sta-3:J1-J'1-both-K4}$]
We let $p \in V(L_1)$. First suppose that $p \in S$.
Let $M$ denote a perfect matching of $J_1 - p - x_1$.
Since $p, x_1 \in S$, the edge $d'_1 \in M$.
Let $M'$ denote a perfect matching of $J'_1$ that contains $g_1$,
whence $g'_1 \in M'$.
Invoking Lemma~\ref{lem:bipartite-shore-properties},
we choose a perfect matching $N'$ of $G[X'_1-u_1-y'-z']$.
Observe that $M \cup M' \cup N' \cup \{e_1\}$ is a perfect matching
of $G-e_2-p-u_3$.

\smallskip
Now suppose that $p \in S'$.
We let $M$ denote a perfect matching of $J_1-p-x'_1$.
Since $p,x'_1 \in S'$, the edge $d_1 \in M$.
We let $w' \in B'_1$ and $t' \in T$ denote the ends of $h'_1$.
Since $J'_1$ is bicritical, we may choose a perfect matching $M'$
of $J'_1 - t' - x_1$. Clearly, $g'_1 \in M'$.
Since $\partial(V(L'_1))$ is a matching in $G$, the vertices $w'$ and $z'$ are distinct.
Invoking Lemma~\ref{lem:bipartite-shore-properties},
we choose a perfect matching $N'$ of $G[X'_1-u_1-w'-z']$.
Observe that $M \cup M' \cup N' \cup \{h'_1,e_1\}$ is a perfect
matching of $G-e_2-p-u_3$.
\end{proof}

It follows from
Assertions~\ref{sta-2:J1-J'1-both-K4}~and~\ref{sta-3:J1-J'1-both-K4}
that $V(L_1) \cap (B_2 \cup B'_2) = \{s_2\}$.
By Lemma~\ref{lem:intersection-singleton-implies-K4},
$J_1 \iso K_4$, contrary to our assumption.
This completes the proof of Proposition~\ref{prp:J1-J'1-both-K4}.
\end{proof}

So far we have proved that the barrier $B_1$ of $G-e_1$ is a doubleton,
and that both bricks of $G-e_1$ are isomorphic to $K_4$.
Now we deduce, using
Lemma~\ref{lem:J1-K4-implies-J2-K4},
that analogous facts
also hold for the graph $G-e_2$.
\begin{corollary}
\label{cor:B2-doubleton-J2-J'2-K4}
The set $B_2$ is a doubleton,
and each of the bricks $J_2$ and $J'_2$ is isomorphic to $K_4$.
In particular, $E(L_2) = \{s_1s_3\}$ and $E(L'_2)=\{t_1t_3\}$.
\end{corollary}
\begin{proof}
Since each of $J_1$ and $J'_1$ is isomorphic to $K_4$,
by Lemma~\ref{lem:J1-K4-implies-J2-K4},
each of $J_2$ and $J'_2$ is also isomorphic to $K_4$;
furthermore, $s_3 \in V(L_2)$ and $t_3 \in V(L'_2)$.
Consequently, $u_3$ has no neighbours in $I_2-v$,
whence Lemma~\ref{lem:bipartite-shore-properties} implies that $B_2$ is
a doubleton. In particular, $B_2 = \{u_1,u_3\}$.
It follows that $u_1$ has a neighbour in each of $L_2$ and $L'_2$.
As per Notation~\ref{Not:s-and-t-vertices}, $s_1 \in V(L_2)$
and $t_1 \in V(L'_2)$. Hence, $E(L_2) = \{s_1s_3\}$
and $E(L'_2) = \{t_1t_3\}$.
This completes the proof of Corollary~\ref{cor:B2-doubleton-J2-J'2-K4}.
\end{proof}

\begin{figure}[!htb]
    \centering
    \begin{tikzpicture}[scale=0.7]

      %edges d' and f'
      \draw (-2,0) -- (-1,-4);
      \draw (-1.9,-2)node[nodelabel]{$d'_1$};
      \draw (-4,0) -- (-3.3,-2.5);
      \draw (-3.9,-2)node[nodelabel]{$f'_1$};

      %edges g' and h'
      \draw (1,-4) to [out=0,in=270] (4,0);
      \draw (4,-2)node[nodelabel]{$g'_1$};
      \draw (2,0) -- (1.5,-2.5);
      \draw (2,-2)node[nodelabel]{$h'_1$};

      %edges d and f
      \draw (-4,0) -- (-2,0);
      \draw (-1,4) to [out=180,in=90] (-4,0)node{};
      \draw (-4.5,0)node[nodelabel]{$s_2$};
      \draw (-3.9,2.3)node[nodelabel]{$d_1$};
      \draw (1,4) to [out=210,in=90] (-2,0)node{};
      \draw (-1.5,0)node[nodelabel]{$s_3$};
      \draw (-2,1.9)node[nodelabel]{$f_1$};

      %edges g and h
      \draw (4,0) -- (2,0);
      \draw (1,4) to [out=0,in=90] (4,0)node{};
      \draw (4.5,0)node[nodelabel]{$t_3$};
      \draw (4,2.3)node[nodelabel]{$h_1$};
      \draw (-1,4) to [out=330,in=90] (2,0)node{};
      \draw (1.5,0)node[nodelabel]{$t_2$};
      \draw (2,1.9)node[nodelabel]{$g_1$};

      %edge e=uv
      \draw (0,2) -- (0,-2);
      \draw (0.4,0)node[nodelabel]{$e_1$};

      %other edges incident at v
      \draw (0,2) -- (-1,4)node{};
      \draw (0,2) -- (1,4)node{};

      %other edges incident at u
      \draw (0,-2) -- (-1,-4)node{};
      \draw (-0.5,-4)node[nodelabel]{$s_1$};
      \draw (0,-2) -- (1,-4)node{};
      \draw (0.5,-4)node[nodelabel]{$t_1$};

      %barrier B
      %\draw (-1,3.5) -- (1,3.5) -- (1,4.5) -- (-1,4.5) -- (-1,3.5);
      \draw (-1,4.5)node[nodelabel]{$u_2$};
      \draw (1,4.5)node[nodelabel]{$u_3$};

      %barrier B'
      \draw (-3,-3.5) -- (3,-3.5) -- (3,-4.5) -- (-3,-4.5) -- (-3,-3.5);
      \draw (-3.5,-4)node[nodelabel]{$B'_1$};

      %isolated set I'-u
      \draw (-2,-5.5) -- (2,-5.5) -- (2,-6.5) -- (-2,-6.5) -- (-2,-5.5);
      \draw (-3,-6)node[nodelabel]{$I'_1-u_1$};

      %\draw (-2.8,-3.2) -- (-3.9,-1.8);
      %\draw (-3.9,-2.5)node[nodelabel]{$d'_1$};
      %\draw (-2,-3.2) -- (-3.1,-1.8);
      %\draw (-2,-2.5)node[nodelabel]{$f'_1$};

      %edges g' and h'

      %\draw (2.8,-3.2) -- (3.9,-1.8);
      %\draw (3.9,-2.5)node[nodelabel]{$h'_1$};
      %\draw (2,-3.2) -- (3.1,-1.8);
      %\draw (2,-2.5)node[nodelabel]{$g'_1$};

      %ends u and v of edge e
      \draw (0,2)node{};
      \draw (0.3,1.8)node[nodelabel]{$v$};
      \draw (0,-2)node{};
      \draw (0.4,-1.8)node[nodelabel]{$u_1$};

      %label G
      %\draw (0,-9)node[nodelabel]{$G$};

    \end{tikzpicture}
    \vline
    \begin{tikzpicture}[scale=0.7]

      %edges d' and f'
      \draw (-2,0) -- (-1,-4);
      \draw (-1.9,-2)node[nodelabel]{$d'_2$};
      \draw (-4,0) -- (-3.3,-2.5);
      \draw (-3.9,-2)node[nodelabel]{$f'_2$};

      %edges g' and h'
      \draw (1,-4) to [out=0,in=270] (4,0);
      \draw (4,-2)node[nodelabel]{$g'_2$};
      \draw (2,0) -- (1.5,-2.5);
      \draw (2,-2)node[nodelabel]{$h'_2$};

      %edges d and f
      \draw (-4,0) -- (-2,0);
      \draw (-1,4) to [out=180,in=90] (-4,0)node{};
      \draw (-4.5,0)node[nodelabel]{$s_1$};
      \draw (-3.9,2.3)node[nodelabel]{$d_2$};
      \draw (1,4) to [out=210,in=90] (-2,0)node{};
      \draw (-1.5,0)node[nodelabel]{$s_3$};
      \draw (-2,1.9)node[nodelabel]{$f_2$};

      %edges g and h
      \draw (4,0) -- (2,0);
      \draw (1,4) to [out=0,in=90] (4,0)node{};
      \draw (4.5,0)node[nodelabel]{$t_3$};
      \draw (4,2.3)node[nodelabel]{$h_2$};
      \draw (-1,4) to [out=330,in=90] (2,0)node{};
      \draw (1.5,0)node[nodelabel]{$t_1$};
      \draw (2,1.9)node[nodelabel]{$g_2$};

      %edge e=uv
      \draw (0,2) -- (0,-2);
      \draw (0.4,0)node[nodelabel]{$e_2$};

      %other edges incident at v
      \draw (0,2) -- (-1,4)node{};
      \draw (0,2) -- (1,4)node{};

      %other edges incident at u
      \draw (0,-2) -- (-1,-4)node{};
      \draw (-0.5,-4)node[nodelabel]{$s_2$};
      \draw (0,-2) -- (1,-4)node{};
      \draw (0.5,-4)node[nodelabel]{$t_2$};

      %barrier B
      %\draw (-1,3.5) -- (1,3.5) -- (1,4.5) -- (-1,4.5) -- (-1,3.5);
      \draw (-1,4.5)node[nodelabel]{$u_1$};
      \draw (1,4.5)node[nodelabel]{$u_3$};

      %barrier B'
      \draw (-3,-3.5) -- (3,-3.5) -- (3,-4.5) -- (-3,-4.5) -- (-3,-3.5);
      \draw (-3.5,-4)node[nodelabel]{$B'_2$};

      %isolated set I'-u
      \draw  (-2,-5.5)  --  (2,-5.5)   --  (2,-6.5)  --  (-2,-6.5)  --
      (-2,-5.5);
      \draw (-3,-6)node[nodelabel]{$I'_2-u_2$};

      %ends u and v of edge e
      \draw (0,2)node{};
      \draw (0.3,1.8)node[nodelabel]{$v$};
      \draw (0,-2)node{};
      \draw (0.4,-1.8)node[nodelabel]{$u_2$};

      %label G
      %\draw (0,-9)node[nodelabel]{$G$};

    \end{tikzpicture}
    \caption{Each of the sets $B_1$ and $B_2$ is a doubleton, and each
brick in $\mathcal{J}$ is isomorphic to~$K_4$}
    \label{fig:B1-B2-doubletons-and-each-brick-J-is-K4}
  \end{figure}
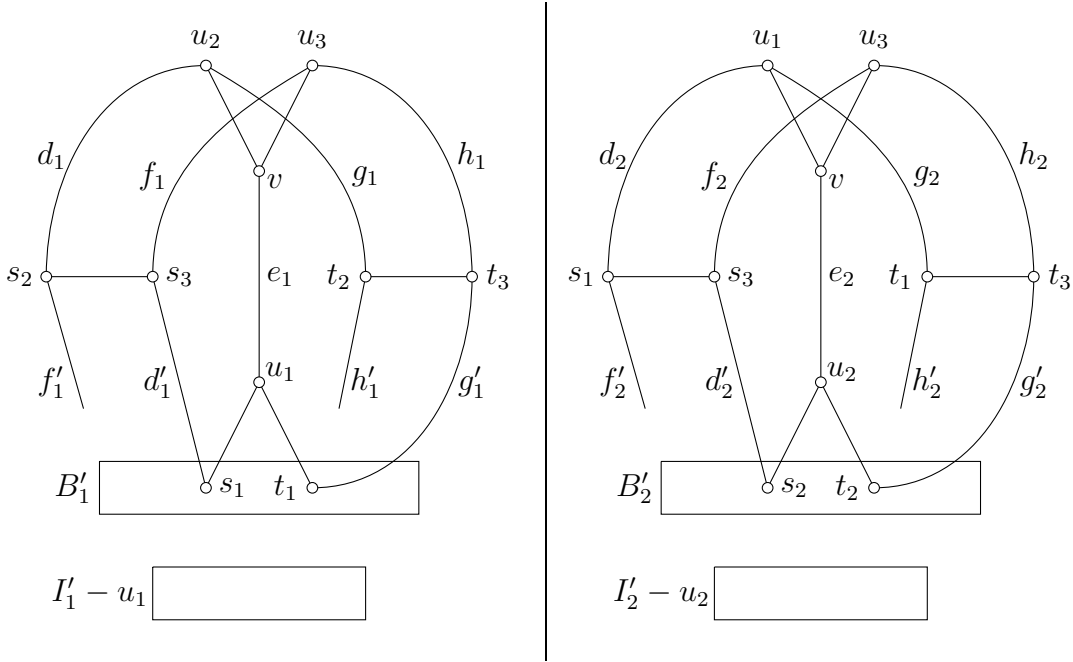

In summary, each of $B_1$ and $B_2$ is a doubleton,
and each member of $\mathcal{J}$ is isomorphic to $K_4$.
We adjust notation so that $d_2,g_2 \in \partial(u_1)$ and $f_2,h_2 \in \partial(u_3)$,
and we adopt the conventions stated below.
See Figure~\ref{fig:B1-B2-doubletons-and-each-brick-J-is-K4}.

\begin{notation}
\label{Not:J2-J'2-inflexible}
Each of the sets
$\{d_2,d'_2\}$ and $\{f_2,f'_2\}$ is a removable doubleton
of $J_2$,
and likewise,
each of the sets
$\{g_2,g'_2\}$ and $\{h_2,h'_2\}$ is a removable doubleton
of $J'_2$.
\end{notation}

For a vertex~$w$, we use $N(w)$ to denote its neighbourhood.
For instance, $N(u_3)= \{v,s_3,t_3\}$.
We let $H^*$ denote the subgraph $G-u_3-N(u_3)$.
The following is easy to see.
\begin{proposition}
The graph~$H^*$ is connected and bipartite, with color classes
$B'_1 \cup \{u_2\} = I'_2 \cup \{s_1,t_1\}$ and
$B'_2 \cup \{u_1\} = I'_1 \cup \{s_2,t_2\}$.   \qed
\end{proposition}

\begin{figure}[!htb]
\centering
\begin{tikzpicture}[scale=0.6]

%barrier B'1
\draw (-3,1.5) -- (-3,3.5) -- (-15,3.5) -- (-15,1.5) -- (-3,1.5);
\draw (-15.5,2.5)node[nodelabel]{$B'_1$};

%barrier B'2
\draw (-3,-1.5) -- (-3,-3.5) -- (-15,-3.5) -- (-15,-1.5) -- (-3,-1.5);
\draw (-15.5,-2.5)node[nodelabel]{$B'_2$};

%edges incident with s3
\draw (4,2) to [out=120,in=60] (-4,2);
\draw (4,2) to [out=315,in=60] (5,-3.5) to [out=240,in=300] (-6,-2);

%edges incident with t3
\draw (4,-2) to [out=240,in=300] (-4,-2);
\draw (4,-2) to [out=45,in=300] (5,3.5) to [out=120,in=60] (-6,2);

%edges incident with u1 and u2
\draw (-2,-2) -- (-4,2);
\draw (-2,-2) -- (-6,2);
\draw (-2,2) -- (-4,-2);
\draw (-2,2) -- (-6,-2);

%vertices s2, t2, y1, z1
%\draw (-4,-2) -- (-8,2);
%\draw (-6,-2) -- (-10,2);
\draw (-4,-2)node{}node[below,nodelabel]{$t_2$};
\draw (-6,-2)node{}node[below,nodelabel]{$s_2$};
%\draw (-8,2)node{}node[above,nodelabel]{$z_1$};
%\draw (-10,2)node{}node[above,nodelabel]{$y_1$};

%vertices s1, t1, y2, z2
%\draw (-4,2) -- (-8,-2);
%\draw (-6,2) -- (-10,-2);
\draw (-4,2)node{}node[above,nodelabel]{$s_1$};
\draw (-6,2)node{}node[above,nodelabel]{$t_1$};
%\draw (-8,-2)node{}node[below,nodelabel]{$y_2$};
%\draw (-10,-2)node{}node[below,nodelabel]{$z_2$};

%claw centered at u3
\draw (2,0) -- (0,0);
\draw (2,0) -- (4,2);
\draw (2,0) -- (4,-2);
\draw (4,2) node{}node[below,nodelabel]{$s_3$};
\draw (4,-2) node{}node[above,nodelabel]{$t_3$};
\draw (2,0)node{}node[right,nodelabel]{$u_3$};

%edges e1 and e2
\draw (0,0) -- (-2,-2);
\draw (-2,-2)node{}node[below,nodelabel]{$u_1$};
\draw (0,0) -- (-2,2)node{};
\draw (-2,2)node{}node[above,nodelabel]{$u_2$};
\draw (0,0)node{}node[left,nodelabel]{$v$};

\end{tikzpicture}
\vspace*{-0.4in}
\caption{The graph $G$ and its connected bipartite
subgraph $H^*[B'_1 \cup \{u_2\},B'_2 \cup \{u_1\}]$}
\label{fig:H*}
\end{figure}
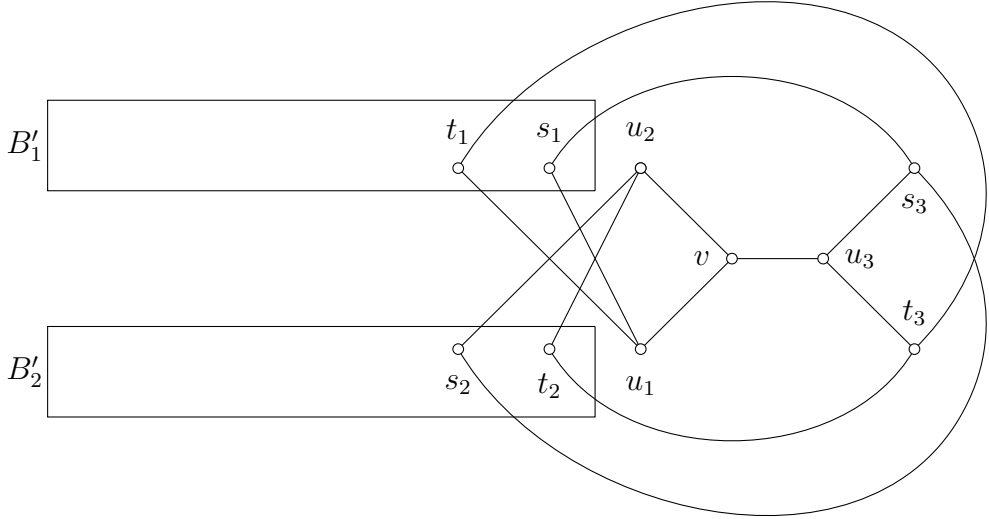

\smallskip
Figure~\ref{fig:H*} shows a different drawing of the graph $G$ that displays
the subgraph~$H^*$ clearly. The reader may easily verify the following.

\begin{proposition}
\label{prp:nonbipartite-subgraph-of-G}
Let $Q$ denote any nonbipartite subgraph of $G$.
Then $V(Q) \cap N(u_3)$ is nonempty. Furthermore, the following hold:
\begin{enumerate}[(i)]
\item If $V(Q) \cap N(u_3) = \{s_3\}$ then $s_1, s_2 \in V(Q)$.
\item If $V(Q) \cap N(u_3) = \{t_3\}$ then $t_1, t_2 \in V(Q)$.
\item if $V(Q) \cap N(u_3) = \{v\}$ then $u_1, u_2 \in V(Q)$. \qed
\end{enumerate}
\end{proposition}

%\vspace*{-0.4in}
\subsection{Nonsolid, nonplanar and non-Pfaffian}
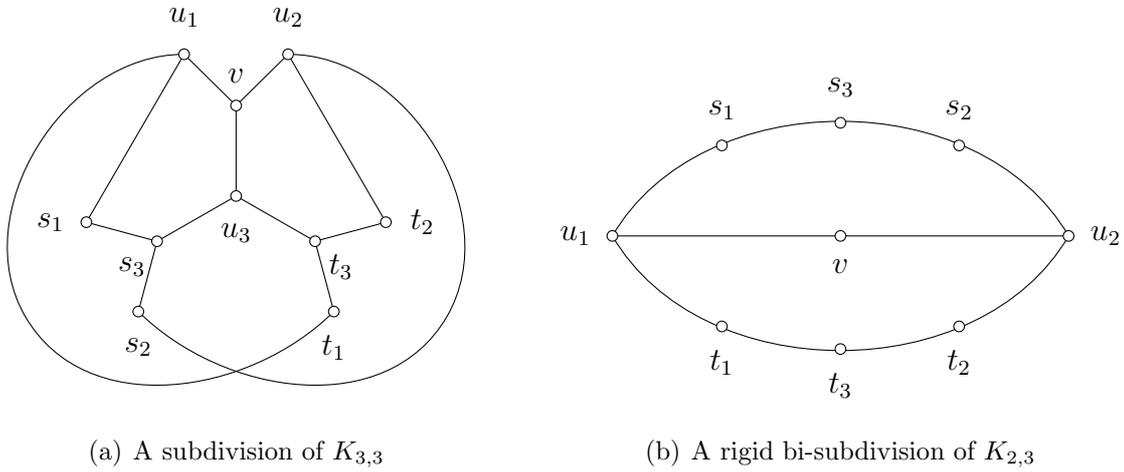
\begin{figure}[!htb]
\centering
\subfigure[A subdivision of $K_{3,3}$]
{
\begin{tikzpicture}[scale=0.8]

\draw (110:2.5) to [out=180,in=120] (210:4) to [out=300,in=225] (310:2.5);
\draw (70:2.5) to [out=0,in=60] (330:4) to [out=240,in=315] (230:2.5);

\draw (110:2.5) -- (190:2.5);
\draw (70:2.5) -- (350:2.5);

\draw (0:0) -- (90:1.5);
\draw (0:0) -- (210:1.5);
\draw (0:0) -- (330:1.5);

\draw (90:1.5) -- (70:2.5);
\draw (90:1.5) -- (110:2.5);

\draw (210:1.5) -- (230:2.5);
\draw (210:1.5) -- (190:2.5);

\draw (330:1.5) -- (350:2.5);
\draw (330:1.5) -- (310:2.5);

\draw (0:0)node{}node[below,nodelabel]{$u_3$};
\draw (90:1.5)node{}node[above,nodelabel]{$v$};
\draw (210:1.5)node{}node[below left,nodelabel]{$s_3$};
\draw (330:1.5)node{}node[below right,nodelabel]{$t_3$};
\draw (70:2.5)node{}node[above,nodelabel]{$u_2$};
\draw (110:2.5)node{}node[above,nodelabel]{$u_1$};
\draw (190:2.5)node{}node[left,nodelabel]{$s_1$};
\draw (230:2.5)node{}node[below,nodelabel]{$s_2$};
\draw (310:2.5)node{}node[below,nodelabel]{$t_1$};
\draw (350:2.5)node{}node[right,nodelabel]{$t_2$};
\end{tikzpicture}
\label{fig:subdivision-of-K33}
}
\subfigure[A rigid bi-subdivision of $K_{2,3}$]
{
\begin{tikzpicture}[scale=1.2]
\draw (0,0) -- (5,0);
\draw (0,0) to [out=60,in=120] (5,0);
\draw (0,0) to [out=300,in=240] (5,0);

\draw (0,0) node{}node[nodelabel,left]{$u_1$};
\draw (5,0) node{}node[nodelabel,right]{$u_2$};

\draw (2.5,0)node{}node[nodelabel,below]{$v$};
\draw (2.5,1.25)node{}node[nodelabel,above]{$s_3$};
\draw (2.5,-1.25)node{}node[nodelabel,below]{$t_3$};

\draw (1.2,1)node{}node[nodelabel,above]{$s_1$};
\draw (1.2,-1)node{}node[nodelabel,below]{$t_1$};

\draw (3.8,1)node{}node[nodelabel,above]{$s_2$};
\draw (3.8,-1)node{}node[nodelabel,below]{$t_2$};

\end{tikzpicture}
\label{fig:rigid-bisubdivision-of-K23}
}
\caption{The brick~$G$ is nonplanar and non-Pfaffian}
\end{figure}

Since $G$ has \qbinv\ edges, it follows from
Theorem~\ref{thm:solid-removable-implies-binv} that $G$ is nonsolid.
Figure~\ref{fig:subdivision-of-K33} shows a subgraph
that is a subdivision of $K_{3,3}$. This subgraph also clearly depicts
$4$ pentagons; we will find this useful
in the proof of Proposition~\ref{prp:Cubeplex}.

\begin{proposition}
\label{prp:G-nonplanar}
The graph $G$ is nonplanar. \qed
\end{proposition}

%\vspace*{-0.2in}

Note that, since $G$ is triangle-free, $s_1$ and $s_2$ are nonadjacent;
likewise, $t_1$ and $t_2$ are nonadjacent. The following is easy to see.

\begin{lemma}
\label{lem:Petersen-adjacencies}
The following are equivalent:
\begin{enumerate}[(i)]
\item At least one of $f'_1$ and $h'_1$ belongs to $\partial(s_1) \cup \partial(t_1)$.
\item At least one of $f'_2$ and $h'_2$ belongs to $\partial(s_2) \cup \partial(t_2)$.
\item $G$ is the Petersen graph. \qed
\end{enumerate}
\end{lemma}

\begin{proposition}
\label{prp:G-non-Pfaffian}
\footnote{By Kasteleyn's Theorem, every non-Pfaffian graph is nonplanar.
However, we have
presented separate certificates for nonplanarity and non-Pfaffian-ness,
since they are easy to observe.}
The graph~$G$ is non-Pfaffian.
\end{proposition}
\begin{proof}
The Petersen graph is non-Pfaffian.
Now suppose that $G$ is not the Petersen graph.
Let $y_2$ denote the end of $f'_1$ in $B'_1$, and
let $z_2$ denote the end of $h'_1$ in $B'_1$.
See Figure~\ref{fig:B1-B2-doubletons-and-each-brick-J-is-K4}.
Note that $y_2$ and $z_2$ need not be distinct.
However, by Lemma~\ref{lem:Petersen-adjacencies},
$\{y_2, z_2\} \cap \{s_1,t_1\} = \emptyset$.

\smallskip
Figure~\ref{fig:rigid-bisubdivision-of-K23}
shows a subgraph of~$G$ that is a bi-subdivision of $K_{2,3}$.
By Lemma~\ref{lem:rigid-bisubdivision-K23-implies-non-Pfaffian},
it suffices to show that this subgraph is a rigid bi-subdivision
of~$K_{2,3}$.

\smallskip
We let $P_1:=u_1s_1s_3s_2u_2$,
$P_2:=u_1vu_2$
and $P_3:=u_1t_1t_3t_2u_2$ denote the three edge-disjoint $u_1u_2$-paths.
For $i,j \in \{1,2,3\}$, where $i < j$, we let $C_{i,j}$
denote the cycle $P_i \cup P_j$. Thus, we only need to argue that
each of the cycles $C_{1,2}, C_{2,3}$ and $C_{1,3}$ is conformal in~$G$.
(The reader may find it useful to trace each of these cycles in
Figure~\ref{fig:B1-B2-doubletons-and-each-brick-J-is-K4}.)

\smallskip
Let us begin with the octagon $C_{1,3}$.
By Lemma~\ref{lem:bipartite-shore-properties},
$G[X'_1-u_1-s_1-t_1]$ has a perfect matching, say~$M$.
Observe that $M \cup \{e_3\}$ is a perfect matching of $G-V(C_{1,3})$.
Thus $C_{1,3}$ is conformal.

\smallskip
Now let us consider the hexagon $C_{1,2}$.
Since $s_1$ and $z_2$ are distinct,
by Lemma~\ref{lem:bipartite-shore-properties},
$G[X'_1-u_1-s_1-z_2]$ has a perfect matching, say~$M'$.
Observe that $M' \cup \{h_1, h'_1\}$
is a perfect matching of $G-V(C_{1,2})$.
Thus $C_{1,2}$ is conformal. An analogous argument
shows that $C_{2,3}$ is conformal.

\smallskip
As discussed earlier,
this completes the proof of Proposition~\ref{prp:G-non-Pfaffian}.
\end{proof}

Thus far we have proved statements (\ref{itm:bricks-J-K4})
and (\ref{itm:G-nonplanar}) of Theorem~\ref{thm:efeccb-two-qbinv-at-vertex}.
It remains to prove statements (\ref{itm:near-bipartite-implies-Cubeplex})
and (\ref{itm:any-more-qbinv-implies-Petersen}).
\subsection{The Cubeplex}
\label{sec:Cubeplex}

In this section, our goal is to prove statement (\ref{itm:near-bipartite-implies-Cubeplex})
of Theorem~\ref{thm:efeccb-two-qbinv-at-vertex}.
\begin{lemma}
\label{lem:Cubeplex}
The following are equivalent:
\begin{enumerate}[(i)]
\item The edges $f'_1$ and $h'_1$ are adjacent.
\item The edges $f'_2$ and $h'_2$ are adjacent.
\item $G$ is the Cubeplex.
\end{enumerate}
\end{lemma}
\begin{proof}
We first prove that (i) implies (ii) and (iii).
Suppose that $y_1 \in B'_1-s_1-t_1$ is a common end of $f'_1$ and $h'_1$.
We observe that $\partial(\{v,u_1,u_2,u_3,s_1,t_1,s_2,t_2,s_3,t_3,y_1\})$
is a $3$-cut of~$G$, and since $G$ is \efec, all three edges are incident
at one vertex, say $y_2$.
In particular, $y_2s_1,y_2t_1,y_2y_1 \in E(G)$, and $G$ is indeed the Cubeplex.
Also, by symmetry, (ii) implies (i) and (iii).

\smallskip
Now we prove that (iii) implies (i) and (ii). Suppose that $G$ is the Cubeplex, whence
$G$ has exactly $12$ vertices.
Consequently, $|B'_1|=3$. Let $y_1$ denote the unique vertex of $B'_1-s_1-t_1$,
and let $y_2$ denote the unique vertex of $I'_1-u_1$. Clearly, $y_2$ is incident
with each vertex of $B'_1$, and $f'_1, h'_1 \in \partial(y_1)$.
This completes the proof of Lemma~\ref{lem:Cubeplex}.
\end{proof}

\begin{proposition}
\label{prp:Cubeplex}
If $G$ is near-bipartite then $G$ is the Cubeplex; furthermore, $e_3$
participates in the (unique) removable doubleton of $G$.
\end{proposition}
\begin{proof}
Assume that $G$ is near-bipartite,
and let $R$ denote a removable doubleton of $G$.
In particular, $G-R$ is a bipartite graph, and $R$ comprises
two nonadjacent edges.

\smallskip
We observe the four $5$-cycles of $G$ that are clearly depicted
in Figure~\ref{fig:subdivision-of-K33};
view each of them as an edge set.
Let
$C_1:=$ $(vu_3, u_3s_3, s_1s_1, s_1u_1, u_1v)$,
$C_2:=$ $(vu_3, u_3s_3, s_3s_2, s_2u_2, u_2v)$,
$C_3:=$ $(vu_3, u_3t_3, t_3t_1, t_1u_1, u_1v)$,
$C_4:=$ $(vu_3, u_3t_3, t_3t_2, t_2u_2, u_2v)$,
and $\mathcal{C}:=$ $\{C_1,C_2,C_3,C_4\}$.
Since $G-R$ is bipartite, the doubleton $R$
meets each member of $\mathcal{C}$.

\smallskip
We claim that $e_3:=vu_3$ lies in $R$. Suppose not.
Observe that any three distinct members of $\mathcal{C}$
have exactly one edge in common --- namely, $e_3$.
Since $e_3 \notin R$, one edge of $R$ meets precisely two members
of $\mathcal{C}$, and the other edge of $R$ meets the remaining two members
of $\mathcal{C}$.
We note that
$C_1 \cap C_2 = \{e_3, u_3s_3\}$,
$C_1 \cap C_3 = \{e_3,vu_1\}$,
$C_1 \cap C_4 = \{e_3\}$,
$C_2 \cap C_3 = \{e_3\}$,
$C_2 \cap C_4 = \{e_3, vu_2\}$, and
$C_3 \cap C_4 = \{e_3,u_3t_3\}$.
These observations imply that either $R=\{u_3s_3,u_3t_3\}$
or $R=\{vu_1,vu_2\}$. This is absurd ---
since the two edges of $R$ are nonadjacent.

\smallskip
We have shown that $e_3 \in R$, whence $G$ has a unique removable doubleton.
Let $f_3$ denote the edge of $R$ that is distinct from $e_3$,
and let $A$~and~$B$ denote the color classes of the bipartite graph
$G-R$ such that $e_3$ has both ends in $A$.
We let $y_1, y_2 \in B$ denote the ends of $f_3$.

\smallskip
Since $e_3=vu_3$, the neighbourhood of $\{v,u_3\}$ is a subset of $B$.
Thus $u_1,u_2,s_3,t_3 \in B$.
Now, we observe that each vertex in $\{s_1,t_1,s_2,t_2\}$
has two neighbours in $\{u_1,u_2,s_3,t_3\}$;
consequently, $s_1,t_1,s_2,t_2 \in A$.
Note that each vertex in $\{u_1,u_2,s_3,t_3\}$
has three neighbours in $A$; whence each of them is distinct from $y_1$,
and from $y_2$. We have thus located six distinct vertices in $A$, and likewise,
six in $B$. Note that $f_3$ has one end in $B'_1-\{s_1,t_1\}$
and another end in $I'_1-u_1$.

\smallskip
We adjust notation so that $y_1 \in B'_1$ and $y_2 \in I'_1-u_1$.
Since $|B'_1| \geq 3$,
by Lemma~\ref{lem:bipartite-shore-properties},
the bipartite graph $G[X'_1-u_1]$ is connected;
this graph contains $f_3$ but it does not contain~$e_3$.
We claim that $G[X'_1-u_1] - f_3$ is disconnected.
Suppose not. Then its color classes are $B'_1$ and $I'_1-u_1$.
However, since $G[X'_1-u_1]-f_3$ is a subgraph of the connected
bipartite graph $G-R$,
one of its color classes is a subset of $A$, and the other is a subset of $B$.
However, since $y_1,y_2 \in B$, we have a contradiction.

\smallskip
Thus $\{f_3\}$ is a $1$-cut of $G[X'_1-u_1]$;
by Lemma~\ref{lem:bipartite-shore-properties}, it must be a trivial cut.
Hence, $y_1$ is an isolated vertex of $G[X'_1-u_1]-f_3$.
We infer that $f'_1, h'_1 \in \partial(y_1)$.
By Lemma~\ref{lem:Cubeplex}, $G$ is indeed the Cubeplex.
This proves Proposition~\ref{prp:Cubeplex}.
\end{proof}
\subsection{The Petersen graph}
\label{sec:Petersen}
In this section, our goal is to prove statement (\ref{itm:any-more-qbinv-implies-Petersen})
of Theorem~\ref{thm:efeccb-two-qbinv-at-vertex}.
\begin{lemma}
\label{lem:e3-qbinv-implies-Petersen}
If $e_3=vu_3$ is \qbinv\ then $G$ is the Petersen graph.
\end{lemma}
\begin{proof}
We assume that $e_3$ is \qbinv\ and, as usual,
we adopt the notation and conventions introduced in
Theorem~\ref{thm:efeccb-qbinv-edge} and
Notation~\ref{Not:8-edges} ---
with the only difference being that all of the notation (except for the vertex $v$)
is decorated with subscript $3$.

\smallskip
All of the preceding arguments,
pertaining to the pair of adjacent
\qbinv\ edges $e_1$~and~$e_2$,
are also applicable to the pair $e_1$~and~$e_3$.
Consequently, $B_3 = \{u_1,u_2\}$, and each of the bricks
$J_3$ and $J'_3$ is isomorphic to $K_4$, whence each of $L_3$ and $L'_3$
is isomorphic to~$K_2$. Furthermore, for $j \in \{1,2\}$,
each of the sets $V(L_3) \cap \{s_j,t_j\}$
 and $V(L'_3) \cap \{s_j,t_j\}$ is a singleton.
Adjust notation so that $s_1 \in V(L_3)$.
Since $G$ is triangle-free,
$t_2 \in V(L_3)$ and $s_1t_2 \in E(G)$.
By Lemma~\ref{lem:Petersen-adjacencies}, $G$ is indeed the Petersen graph.
This proves Lemma~\ref{lem:e3-qbinv-implies-Petersen}.
\end{proof}

\begin{proposition}
\label{prp:Petersen}
If $G$ has a \qbinv\ edge, distinct from $e_1$ and $e_2$,
then $G$ is the Petersen graph.
\end{proposition}
\begin{proof}
Assume that $e^*:=v^*u^*$ is a \qbinv\ edge of $G$, distinct from $e_1$~and~$e_2$.
If $e^* = e_3$
then the desired conclusion holds by
Lemma~\ref{lem:e3-qbinv-implies-Petersen}

\smallskip
Now suppose that $e^* \neq e_3$, whence $v \notin \{v^*,u^*\}$.
We adopt the notation and conventions introduced in
Theorem~\ref{thm:efeccb-qbinv-edge} ---
as shown in Figure~\ref{fig:e*-qbinv}.

\smallskip
We first consider the case in which $e^* \in \{u_3s_3, u_3t_3\}$.
Adjust notation so that $u_3 = v^*$ and $s_3 = u^*$,
whence $v,t_3 \in B$ and $s_1,s_2 \in B'$.
Since $u_1$ is a common neighbour of $v \in B$ and $s_1 \in B'$,
we infer that $u_1 \in V(L) \cup V(L')$.
Adjust notation so that $u_1 \in V(L)$.
Since $\partial(V(L))$ is not a matching, $L \iso K_2$ and $E(L) = \{u_1t_1\}$.
A similar argument shows that $L' \iso K_2$ and $E(L') = \{u_2t_2\}$.
Consequently, $v$ has no neighbours in $I-v^*$, whence $B$ is a doubleton.
Observe that $G[X'-u^*]$ is a subgraph of the connected bipartite graph
$G-u_3-N(u_3)$. If $|B'| \geq 3$ then, by Lemma~\ref{lem:bipartite-shore-properties},
$G[X'-u^*]$ is a connected subgraph with color classes $B'$ and $I'-u^*$;
however, this results in a contradiction since $s_1$ and $s_2$ lie in distinct
color classes of $G-u_3-N(u_3)$. Thus $B'$ is a doubleton, whence $|V(G)|=10$.
By Proposition~\ref{prp:order-10}, $G$ is indeed the Petersen graph.

\smallskip
Now suppose that $e^* \notin \{u_3s_3,u_3t_3\}$.
Thus, $u_3 \notin \{v^*,u^*\}$.

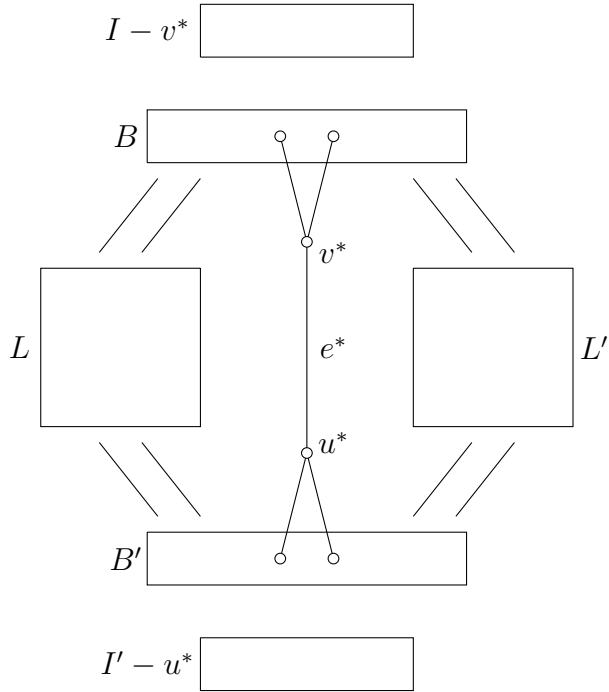
\begin{figure}[!htb]
    \centering
    \begin{tikzpicture}[scale=0.7]

      %edge e=uv
      \draw (0,2) -- (0,-2);
      \draw (0.5,0)node[nodelabel]{$e^*$};

      %other edges incident at v
      \draw (0,2) -- (-0.5,4)node{};
      \draw (0,2) -- (0.5,4)node{};

      %other edges incident at u
      \draw (0,-2) -- (-0.5,-4)node{};
      \draw (0,-2) -- (0.5,-4)node{};

      %barrier B
      \draw (-3,3.5) -- (3,3.5) -- (3,4.5) -- (-3,4.5) -- (-3,3.5);
      \draw (-3.4,4)node[nodelabel]{$B$};

      %isolated set I-v
      \draw (-2,5.5) -- (2,5.5) -- (2,6.5) -- (-2,6.5) -- (-2,5.5);
      \draw (-3,6)node[nodelabel]{$I-v^*$};

      %barrier B'
      \draw (-3,-3.5) -- (3,-3.5) -- (3,-4.5) -- (-3,-4.5) -- (-3,-3.5);
      \draw (-3.4,-4)node[nodelabel]{$B'$};

      %isolated set I'-u
      \draw (-2,-5.5) -- (2,-5.5) -- (2,-6.5) -- (-2,-6.5) -- (-2,-5.5);
      \draw (-3,-6)node[nodelabel]{$I'-u^*$};

      %component L
      \draw (-2,1.5) -- (-5,1.5) -- (-5,-1.5) -- (-2,-1.5) -- (-2,1.5);
      \draw (-5.4,0)node[nodelabel]{$L$};

      %component L'
      \draw (2,1.5) -- (5,1.5) -- (5,-1.5) -- (2,-1.5) -- (2,1.5);
      \draw (5.4,0)node[nodelabel]{$L'$};

      %edges d and f
      \draw (-2.8,3.2) -- (-3.9,1.8);
      %\draw (-3.8,2.5)node[nodelabel]{$d$};
      \draw (-2,3.2) -- (-3.1,1.8);
      %\draw (-2.2,2.5)node[nodelabel]{$f$};

      %edges g and h
      \draw (2.8,3.2) -- (3.9,1.8);
      %\draw (3.8,2.5)node[nodelabel]{$h$};
      \draw (2,3.2) -- (3.1,1.8);
      %\draw (2.1,2.5)node[nodelabel]{$g$};

      %edges d' and f'
      \draw (-2.8,-3.2) -- (-3.9,-1.8);
      %\draw (-3.8,-2.5)node[nodelabel]{$d'$};
      \draw (-2,-3.2) -- (-3.1,-1.8);
      %\draw (-2,-2.5)node[nodelabel]{$f'$};

      %edges g' and h'
      \draw (2.8,-3.2) -- (3.9,-1.8);
      %\draw (3.8,-2.5)node[nodelabel]{$h'$};
      \draw (2,-3.2) -- (3.1,-1.8);
      %\draw (2.1,-2.5)node[nodelabel]{$g'$};

      %ends u and v of edge e
      \draw (0,2)node{};
      \draw (0.5,1.8)node[nodelabel]{$v^*$};
      \draw (0,-2)node{};
      \draw (0.5,-1.8)node[nodelabel]{$u^*$};

      %label G
      %\draw (0,-9)node[nodelabel]{$G$};

    \end{tikzpicture}
\caption{Illustration for the proof of Proposition~\ref{prp:Petersen}}
\label{fig:e*-qbinv}
\end{figure}

\smallskip
We define four subgraphs as follows.
We let $G_1:=G[V(L) \cup X]$,
$G_2:=G[V(L') \cup X']$,
$G_3:=G[V(L') \cup X]$, and
$G_4:=G[V(L) \cup X']$.
By Theorem~\ref{thm:efeccb-qbinv-edge}(\ref{itm:nonbipartite-subgraphs}),
each of these four subgraphs is nonbipartite. Consequently,
Propositon~\ref{prp:nonbipartite-subgraph-of-G} applies to all of them.
In particular, each of these four subgraphs meets the set $N(u_3) = \{v,s_3,t_3\}$;
in the arguments that follow, we will use this fact implicitly.

\begin{assertion}
\label{sta:nborhood-of-u3-X-X'}
The set $N(u_3)$ meets $X \cup X'$.
\end{assertion}
\begin{proof}[Proof of Assertion~$\ref{sta:nborhood-of-u3-X-X'}$]
Suppose, to the contrary, that $N(u_3) \cap (X \cup X')$ is empty.
Thus $v,s_3,t_3 \in V(L) \cup V(L')$;
adjust notation so that two of them lie in $V(L)$,
and the third one lies in $V(L')$.
Since $u_3$ is a common neighbour, $u_3 \in X \cup X'$;
adjust notation so that $u_3 \in X$. Thus there exist
two adjacent edges joining $X$ and $L$. This
contradicts Theorem~\ref{thm:efeccb-qbinv-edge}(\ref{itm:L-Lp}).
\end{proof}

\begin{assertion}
\label{sta:nborhood-of-u3-L-L'}
The set $N(u_3)$ meets $V(L) \cup V(L')$.
\end{assertion}
\begin{proof}[Proof of Assertion~$\ref{sta:nborhood-of-u3-L-L'}$]
Suppose, to the contrary, that $N(u_3) \cap (V(L) \cup V(L'))$ is empty.
Thus $v,s_3,t_3 \in X \cup X'$;
adjust notation so that two of them lie in $X$, and the third one lies in $X'$.
Since $u_3$ is a common neighbour, we infer that, in fact, two of them lie in $B$,
the third one lies in $B'$ and $u_3 \in V(L) \cup V(L')$,
contrary to Theorem~\ref{thm:efeccb-qbinv-edge}(\ref{itm:L-Lp}).
\end{proof}

We may adjust notation so that one vertex of $N(u_3)$ lies in $X$,
another vertex of $N(u_3)$ lies in $V(L)$, and the third vertex lies in $X' \cup V(L')$.
We consider two cases depending on whether the third vertex lies in $X'$ or
in $V(L')$.

\bigskip
\noindent
{\bf Case 1:} $N(u_3) \cap V(L')$ is nonempty.

\smallskip
\noindent
We consider two cases depending on whether or not $v$ lies in $X$.

\medskip
\noindent
{\bf Case 1.1:} $v \in X$.

\smallskip
\noindent
Adjust notation so that $s_3 \in V(L)$ and $t_3 \in V(L')$.
By our assumption $v \neq v^*$, whence the common neighbour $u_3 \in B$,
and $v \in I - v^*$. Consequently, $u_1, u_2 \in B$.
Since $G_4$ is nonbipartite,
Proposition~\ref{prp:nonbipartite-subgraph-of-G} implies that $s_1, s_2 \in V(L)$.
This is absurd --- since it results in three distinct edges joining $L$ and $X$.

\medskip
\noindent
{\bf Case 1.2:} $v \notin X$.

\smallskip
\noindent
Adjust notation so that $t_3 \in X, s_3 \in V(L)$ and $v \in V(L')$.
The common neighbour $u_3 \in B$, and $t_3 \in I$.
First suppose that $t_3 \in I-v^*$, whence $t_1, t_2 \in B$.
Since $G_2$ is nonbipartite,
Proposition~\ref{prp:nonbipartite-subgraph-of-G} implies that $u_1, u_2 \in V(L')$.
This is absurd --- since it results in three distinct edges joining $L'$ and $X$.

\smallskip
Now suppose that $t_3=v^*$, whence $B$ is a doubleton.
Adjust notation so that $t_2 \in B$ and $t_1=u^*$.
Consequently, $u_1 \in B'$.
Observe that $\partial(V(L'))$ is not a matching,
whence $L' \iso K_2$ and $E(L') = \{vu_2\}$.
Now it follows that $s_2 \in B'$.
Consequently, $\partial(V(L))$ is not a matching;
thus $L \iso K_2$ and $E(L) = \{s_3s_1\}$.
This implies that $t_2s_1 \in E(G)$.
By Lemma~\ref{lem:Petersen-adjacencies},
$G$ is indeed the Petersen graph.

\bigskip
\noindent
{\bf Case 2:} $N(u_3) \cap X'$ is nonempty.

\smallskip
\noindent
We consider two cases depending on whether or not $v$ lies in $V(L)$.

\medskip
\noindent
{\bf Case 2.1:} $v \in V(L)$.

\smallskip
\noindent
Adjust notation so that $s_3 \in X$ and $t_3 \in X'$.
We infer that the common neighbour $u_3 \in V(L)$, and
that $s_3 \in B$ and $t_3 \in B'$.
Since $\partial(V(L))$ is not a matching, $L \iso K_2$ and $E(L) = \{vu_3\}$.
Consequently, one of $u_1$ and $u_2$ lies in $B$, and the other lies in $B'$.
Adjust notation so that $u_1 \in B$ and $u_2 \in B'$.
Since $s_2$ is a common neighbour of $s_3 \in B$ and $u_2 \in B'$,
we infer that $s_2 \in V(L')$.
Likewise, $t_1$ is a common neighbour of $u_1 \in B$ and $t_3 \in B'$,
whence $t_1 \in V(L')$.
Also, $\partial(V(L'))$ is not a matching; consequently, $L' \iso K_2$
and $s_2t_1 \in E(G)$.
By Lemma~\ref{lem:Petersen-adjacencies}, $G$ is indeed the Petersen graph.

\medskip
\noindent
{\bf Case 2.2:} $v \notin V(L)$.

\smallskip
\noindent
Adjust notation so that $t_3 \in V(L)$, $s_3 \in X$ and $v \in X'$.
Thus their common neighbour $u_3 \in V(L)$, and $s_3 \in B$ and $v \in B'$.
Since $\partial(V(L))$ is not a matching, $L \iso K_2$ and $E(L)=\{t_3u_3\}$.
One of $t_1$ and $t_2$ lies in $B$, and the other lies in $B'$.
Adjust notation so that $t_1 \in B$ and $t_2 \in B'$.
Since $u_1$ is a common neighbour of $t_1 \in B$ and $v \in B'$,
we infer that $u_1 \in V(L')$.
Since $\partial(V(L'))$ is not a matching, $L' \iso K_2$
and $E(L') = \{u_1s_1\}$.
We observe that $vu_2s_2s_3$ is a path that joins $v \in B'$ and $s_3 \in B$;
this implies that $u_2 = u^*$ and $v^* = s_2$.
Now $v \in B'$ has no neighbours in $I'-u^*$, whence $B'$ is a doubleton,
and $t_2s_1 \in E(G)$.
By Lemma~\ref{lem:Petersen-adjacencies}, $G$ is indeed the Petersen graph.

\bigskip
Thus, in each case, we have either arrived at a contradiction, or we have
arrived at the conclusion that $G$ is the Petersen graph.
This completes the proof of Proposition~\ref{prp:Petersen}.
\end{proof}

This completes the proof of the
Main Theorem (\ref{thm:efeccb-two-qbinv-at-vertex}). \hfill \qed

\section{Consequences of the Main Theorem}
\label{sec:main-theorem-consequences}

We let $\mathcal{G}$ denote the set of all {\efeccb}s,
except for $K_4$.
We partition $\mathcal{G}$ into two subsets
$\mathcal{G}_1$ and $\mathcal{G}_2$ as follows.
A graph $G \in \mathcal{G}$ belongs to $\mathcal{G}_2$
if and only if $G$ has a vertex that is incident with
at least two \qbinv\ edges.
In particular, the Petersen graph and the Cubeplex
are members of $\mathcal{G}_2$.
The Cubeplex is the only near-bipartite member of $\mathcal{G}_2$,
and it has precisely two \qbinv\ edges and a unique removable doubleton;
whence it has $14$ \binv\ edges. 
Now, we prove Theorem~\ref{thm:efeccnnbb-lower-bound}
which states that, if $G$ is any non-near-bipartite member of~$\mathcal{G}$, distinct
from the Petersen graph, then at least two-third of its edges are \binv.

\begin{proof}[Proof of Theorem~$\ref{thm:efeccnnbb-lower-bound}$]
First suppose that $G$ is a non-near-bipartite member of~$\mathcal{G}_2$
distinct from the Petersen graph.
By Theorem~\ref{thm:efeccb-two-qbinv-at-vertex},
$G$ has precisely two \qbinv\ edges
and $\frac{3|V(G)|}{2}-2$ \binv\ edges.

\smallskip
Now suppose that $G$ is a non-near-bipartite member of $\mathcal{G}_1$,
whence each vertex of $G$ is incident with at least two \binv\ edges;
consequently, $G$ has at least $|V(G)|$ \binv\ edges.
This completes the proof of Theorem~\ref{thm:efeccnnbb-lower-bound}.
\end{proof}

\smallskip
We now point out some other interesting facts that are immediate consequences
of Theorem~\ref{thm:efeccb-two-qbinv-at-vertex}.
Let $G$ be any member of $\mathcal{G}_2$
that is distinct from the Petersen graph.
Then $G$ has a unique vertex that $v$ that is incident with precisely
two \qbinv\ edges $e_1$~and~$e_2$. Let $e_3:=vu_3$ denote the third edge
incident with $v$.
It follows that the automorphism group of $G$ has at least
two singleton orbits: $\{v\}$ and $\{u_3\}$.
In particular,
the Petersen graph is the only vertex-transitive member of $\mathcal{G}_2$.

\smallskip
As mentioned earlier, the Cubeplex and the Twinplex are the only two
near-bipartite graphs that are minimally \mbox{non-Pfaffian}
(see \cite{fili01,clm12}). It has always intrigued us that
the Twinplex is far more symmetric than the Cubeplex; in particular,
the automorphism group of the Twinplex has no singleton orbits.
The discussion in the preceding paragraph perhaps throws some more
light on this phenomenon.

\smallskip
As per Theorem~\ref{thm:efeccb-two-qbinv-at-vertex},
each member of $\mathcal{G}_2$ is nonsolid, nonplanar and non-Pfaffian.
By Theorem~\ref{thm:equivalence-of-solid-and-C6bar-free},
each member of $\mathcal{G}_2$,
except the Petersen graph, is $\overline{C_6}$-based.

\smallskip
Now let $G \in \mathcal{G}_2$ and let $e$ denote a \qbinv\ edge.
Then both bricks of $G-e$ are isomorphic to $K_4$.
By Theorem~\ref{thm:J-free-reduction-to-bricks},
$G-e$ is $K_4$-based.
Consequently, $G$ is also $K_4$-based.

\smallskip
The following is a brief summary of the above discussion
pertaining to the bricks in~$\mathcal{G}_2$.
\begin{corollary}
Let $G$ denote any member of $\mathcal{G}_2$ that is distinct
from the Petersen graph.
Then the following statements hold:
\begin{enumerate}[(i)]
\item $G$ has exactly two \qbinv\ edges, say~$e_1$ and $e_2$, and these are adjacent.
\item If $G$ is not the Cubeplex then each edge, except $e_1$~and~$e_2$,
is \binv; in particular, $G$ is non-near-bipartite.
\item The automorphism group of $G$ has at least two singleton orbits.
\item $G$ is nonplanar and non-Pfaffian.
\item $G$ is $\overline{C_6}$-based and nonsolid.
\item $G$ is $K_4$-based. \qed
\end{enumerate}
\end{corollary}

\section{An infinite family of cubic bricks}
\label{sec:infinite-family}

Let us recall Theorem~\ref{thm:efeccb-removable-edges},
which states that, if $e$ is a removable edge of an \efeccb,
then $b(G-e) \in \{1,2\}$.

\smallskip
In this section, we will demonstrate
that the conclusion of Theorem~\ref{thm:efeccb-removable-edges}
does not hold
for cubic bricks, in general. In particular, for any integer $k \geq 3$,
we describe how one may
construct a cubic brick $G$ that has a removable edge $e$
so that $b(G-e) = k$.

\smallskip
We will start
from a cubic brace $H$ of order $2k+2$, and we will perform
a few operations in order to obtain $G$.
In particular, we will require the operation of `splicing' two graphs,
that is defined formally in \cite{lckm18}.
The following is easy to prove.
\begin{proposition}
\label{prp:splicing}
A splicing of any two {\mcg}s yields another \mcg.
A splicing of any two bicritical graphs yields another bicritical graph. \qed
\end{proposition}

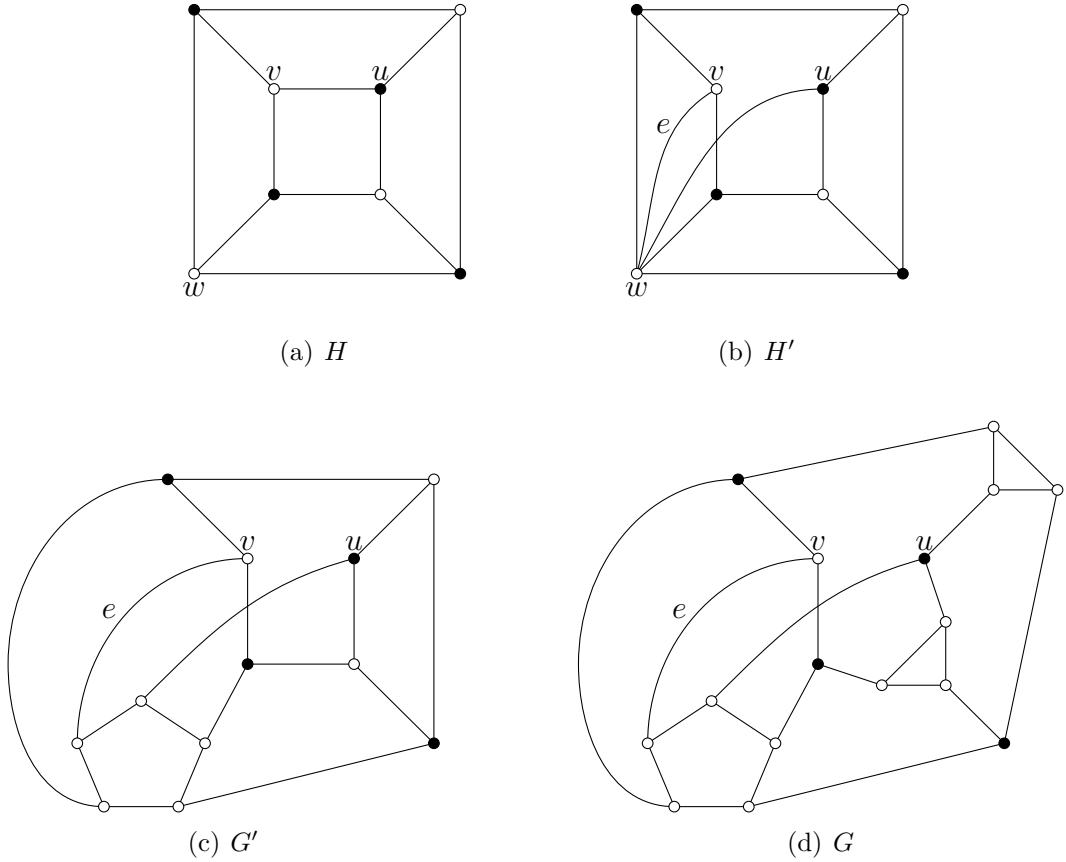
\begin{figure}[!htb]
\centering
\subfigure[$H$]
{
\begin{tikzpicture}[scale=0.7]
\draw (1,1) -- (1,-1) -- (-1,-1) -- (-1,1) -- (1,1);
\draw (2.5,2.5) -- (2.5,-2.5) -- (-2.5,-2.5) -- (-2.5,2.5) -- (2.5,2.5);

\draw (1,1) -- (2.5,2.5);
\draw (1,-1) -- (2.5,-2.5);
\draw (-1,-1) -- (-2.5,-2.5);
\draw (-1,1) -- (-2.5,2.5);

\draw (1,1)node[fill=black]{};
\draw (2.5,2.5)node{};
\draw (1,-1)node{};
\draw (2.5,-2.5)node[fill=black]{};
\draw (-1,-1)node[fill=black]{};
\draw (-2.5,-2.5)node{};
\draw (-1,1)node{};
\draw (-2.5,2.5)node[fill=black]{};

\draw (-1,1.3)node[nodelabel]{$v$};
\draw (1,1.3)node[nodelabel]{$u$};
\draw (-2.5,-2.8)node[nodelabel]{$w$};
\end{tikzpicture}
\label{fig:cube}
}
\hspace*{0.5in}
\subfigure[$H'$]
{
\begin{tikzpicture}[scale=0.7]
%two new edges
\draw (-2.5,-2.5) to [out=75,in=210] (-1,1);
\draw (-2,0.3)node[nodelabel]{$e$};
\draw (-2.5,-2.5) to [out=60,in=180] (1,1);

\draw (1,1) -- (1,-1) -- (-1,-1) -- (-1,1);
\draw (2.5,2.5) -- (2.5,-2.5) -- (-2.5,-2.5) -- (-2.5,2.5) -- (2.5,2.5);

\draw (1,1) -- (2.5,2.5);
\draw (1,-1) -- (2.5,-2.5);
\draw (-1,-1) -- (-2.5,-2.5);
\draw (-1,1) -- (-2.5,2.5);

\draw (1,1)node[fill=black]{};
\draw (2.5,2.5)node{};
\draw (1,-1)node{};
\draw (2.5,-2.5)node[fill=black]{};
\draw (-1,-1)node[fill=black]{};
\draw (-2.5,-2.5)node{};
\draw (-1,1)node{};
\draw (-2.5,2.5)node[fill=black]{};

\draw (-1,1.3)node[nodelabel]{$v$};
\draw (1,1.3)node[nodelabel]{$u$};
\draw (-2.5,-2.8)node[nodelabel]{$w$};
\end{tikzpicture}
\label{fig:cube-H'}
}

\vspace*{0.2in}
\subfigure[$G'$]
{
\begin{tikzpicture}[scale=0.7]
%spokes of W5
\draw (-3.7,-3.7) to [out=180,in=270] (-5.5,-1) to [out=90,in=180] (-2.5,2.5);
\draw (-2.3,-3.7) -- (2.5,-2.5);
\draw (-1.8,-2.5) -- (-1,-1);
\draw (-3,-1.7) to [out=45,in=195] (1,1);
\draw (-4.2,-2.5) to [out=90,in=180] (-1,1);
\draw (-3.6,0)node[nodelabel]{$e$};

%\draw (-4.2,-2.5) -- (-3.7,-3.7)node{} -- (-2.3,-3.7)node{} -- (-1.8,-2.5)node{}
%-- (-3,-1.7)node{} -- (-4.2,-2.5)node{};

%two new edges
%\draw (-2.5,-2.5) to [out=75,in=210] (-1,1);
%\draw (-2,0.3)node[nodelabel]{$e$};
%\draw (-2.5,-2.5) to [out=60,in=180] (1,1);

\draw (1,1) -- (1,-1) -- (-1,-1) -- (-1,1);
\draw (2.5,2.5) -- (2.5,-2.5);
\draw (-2.5,2.5) -- (2.5,2.5);

\draw (1,1) -- (2.5,2.5);
\draw (1,-1) -- (2.5,-2.5);
%\draw (-1,-1) -- (-2.5,-2.5);
\draw (-1,1) -- (-2.5,2.5);

\draw (1,1)node[fill=black]{};
\draw (2.5,2.5)node{};
\draw (1,-1)node{};
\draw (2.5,-2.5)node[fill=black]{};
\draw (-1,-1)node[fill=black]{};
%\draw (-2.5,-2.5)node{};
\draw (-1,1)node{};
\draw (-2.5,2.5)node[fill=black]{};

\draw (-1,1.3)node[nodelabel]{$v$};
\draw (1,1.3)node[nodelabel]{$u$};
%\draw (-2.5,-2.8)node[nodelabel]{$w$};

%vertices of W5
\draw (-4.2,-2.5) -- (-3.7,-3.7)node{} -- (-2.3,-3.7)node{} -- (-1.8,-2.5)node{}
-- (-3,-1.7)node{} -- (-4.2,-2.5)node{};

\end{tikzpicture}
\label{fig:cube-G'}
}
\hspace*{0.5in}
\subfigure[$G$]
{
\begin{tikzpicture}[scale=0.7]

%second K4 splicing
\draw (2.3,2.3) -- (1,1);
\draw (3.5,2.3) -- (2.5,-2.5);
\draw (2.3,3.5) -- (-2.5,2.5);
%\draw (2.3,2.3) -- (3.5,2.3)node{} -- (2.3,3.5)node{} -- (2.3,2.3)node{};

%first K4 splicing
\draw (1.4,-1.4) -- (2.5,-2.5);
\draw (1.4,-0.2) -- (1,1);
\draw (0.2,-1.4) -- (-1,-1);

%spokes of W5
\draw (-3.7,-3.7) to [out=180,in=270] (-5.5,-1) to [out=90,in=180] (-2.5,2.5);
\draw (-2.3,-3.7) -- (2.5,-2.5);
\draw (-1.8,-2.5) -- (-1,-1);
\draw (-3,-1.7) to [out=45,in=195] (1,1);
\draw (-4.2,-2.5) to [out=90,in=180] (-1,1);
\draw (-3.6,0)node[nodelabel]{$e$};

\draw (-1,-1) -- (-1,1);

\draw (-1,1) -- (-2.5,2.5);

\draw (1,1)node[fill=black]{};
\draw (2.5,-2.5)node[fill=black]{};
\draw (-1,-1)node[fill=black]{};
\draw (-1,1)node{};
\draw (-2.5,2.5)node[fill=black]{};

\draw (-1,1.3)node[nodelabel]{$v$};
\draw (1,1.3)node[nodelabel]{$u$};

%vertices of W5
\draw (-4.2,-2.5) -- (-3.7,-3.7)node{} -- (-2.3,-3.7)node{} -- (-1.8,-2.5)node{}
-- (-3,-1.7)node{} -- (-4.2,-2.5)node{};

\draw (1.4,-1.4) -- (1.4,-0.2)node{} -- (0.2,-1.4)node{} -- (1.4,-1.4)node{};

\draw (2.3,2.3) -- (3.5,2.3)node{} -- (2.3,3.5)node{} -- (2.3,2.3)node{};

\end{tikzpicture}
\label{fig:cube-G}
}
\caption{Constructing a cubic brick $G$ with a removable edge $e$ so that $b(G-e)=3$}
\label{fig:infinite-family-example}
\end{figure}

\smallskip
Some well-known examples of infinite families of cubic braces are:
the prisms of order~$4k$ where $k \geq 2$,
and the M{\"o}bius ladders of order $4k+2$ where $k \geq 1$.
See \cite{komu16} for definitions. The smallest prism is the cube graph,
shown in Figure~\ref{fig:cube}, and we shall use this
to illustrate the construction that we are about to describe.

\smallskip
Let $k \geq 3$ be an integer. We consider any cubic brace $H[A,B]$
of order $2k+2$. We choose an edge $uv$,
adjusting notation so that $v \in A$ and $u \in B$,
and we choose a vertex $w \in A$ such that $u$~and~$w$ are nonadjacent.
(Such a choice is possible since $H$ is of order eight or more.)
Now, let $H':=H-uv+uw+vw$.
Observe that $H'$ is a simple graph, in which vertex~$w$ has degree five,
and every other vertex is cubic.
Also, $H'$ is not \mc; in particular,
the edge $e:=vw$ is inadmissible.
See Figures~\ref{fig:cube}~and~\ref{fig:cube-H'}.

\smallskip
We let $G'$ denote a cubic graph obtained by splicing $H'$,
and the odd wheel $W_5$, at their only noncubic vertices.
(This is equivalent to `replacing' the vertex $w$, of $H'$, by a $5$-cycle,
so as to obtain a cubic graph.)
See Figures~\ref{fig:cube-H'}~and~\ref{fig:cube-G'}.

\smallskip
We let $G$ denote the cubic graph obtained by splicing $G'$ with a copy of $K_4$
at each vertex in the set~$A-v-w$.
(Splicing a cubic graph with $K_4$, at a given vertex,
is equivalent to `replacing' that vertex by a triangle
so as to obtain another cubic graph.)
See Figures~\ref{fig:cube-G'}~and~\ref{fig:cube-G}.

\smallskip
In the proof of the following, we will omit a few details.
However, we provide all of the important steps.

\begin{proposition}
\label{prp:infinite-family}
The graph $G$ is a cubic brick, and $e$ is a removable edge of $G$.
Furthermore, $b(G-e)=k$.
\end{proposition}
\begin{proof}
First of all, we argue why $G$ is a brick.
As noted in Section~\ref{sec:tight-cut-decompositions},
every cubic brace is \efec. In particular, $H$ is \efec.
Using this fact, one may deduce that $H'$ is free of nontrivial $3$-cuts.
Now, since $G'$ is obtained
by splicing $H'$ and the odd wheel $W_5$, one may infer that $G'$ is also
free of nontrivial $3$-cuts. Thus, $G'$
is an \efeccg. Since $G'$ is nonbipartite,
Corollary~\ref{cor:efeccg-brick-brace} implies that $G'$ is a brick.
In particular, by Theorem~\ref{thm:elp-brick-characterization},
$G'$ is a bicritical graph.

\smallskip
Since $G$ is obtained from the bicritical graph~$G'$
by repeatedly splicing with copies of the bicritical graph~$K_4$,
it follows from Proposition~\ref{prp:splicing}
that $G$ is also bicritical. Clearly, $G$ is cubic.
Thus, by Corollary~\ref{cor:bicritical-not-brick}, $G$ is a brick.

\smallskip
Now we argue why $e$ is a removable edge of $G$, or equivalently,
why $G-e$ is \mc.
In a brace (of order six or more), every edge is removable.
Thus $H-uv$ is a bipartite \mcg.
Proposition~\ref{prp:bip-mc-characterization} implies that
any bipartite graph, which is obtained from a bipartite \mcg\ by
adding an edge, is also \mc.
Thus, $H-uv+uw$, which is the same as $H'-e$,
is also \mc.
Observe that $G-e$ may be obtained by splicing $H'-e$ and $W_5-e$,
each of which is \mc, whence $G-e$ is also \mc.

\smallskip
Finally, we argue why $b(G-e)=k$. Note that the set $B$ is a barrier in $G-e$
and that $G-e-B$ has $k$ nontrivial components (and one trivial component,
namely the vertex~$v$). In particular, $k-1$ of the nontrivial components
are triangles, and the last one is a $5$-cycle.
This yields a laminar family of $k$ nontrivial tight cuts, and each of them
produces exactly one brick (that is in fact isomorphic to $K_4$).
Thus $b(G-e)=k$.

\smallskip
This completes the proof of Proposition~\ref{prp:infinite-family}.
\end{proof}

%%%%%%%%%%%%%%%%%%%%%%%%%%%%%%%%%%%%%%%%%%%%%%%%%%%%%%%
\subsection*{Acknowledgements}

This work commenced when the first author was a postdoc at the
University of Waterloo (Canada) and was visiting Massey University (New Zealand).
A significant portion of this work was done while the first author was a
postdoc at the University of Vienna, and also during his visit to UFMS Campo Grande (Brazil).
The work was completed during his postdoc at the University of Campinas (Brazil).
The first author extends his gratitude to all institutions involved, and especially
to Jochen K{\"o}nemann (Waterloo), and to Ilse Fischer (Vienna), for providing financial support for the
aforementioned visits.

%%%%%%%%%%%%%%%%%%%%%%%%%%%%%%%%%%%%%%%%%%%%%%%%%%%%%%%
% You do not have to use the same format for your references, but 
%    include everything in this file.  Don't use natbib please.
% If you use BibTeX to create a bibliography, copy the .bbl file into here.

%\bibliographystyle{plain}
%\bibliography{clm} 

\end{document}